\documentclass[11pt]{amsart}

\usepackage{amscd,calc,graphicx}

\newtheorem{theorem}{Theorem}[section]
\newtheorem{proposition}[theorem]{Proposition}
\newtheorem{corollary}[theorem]{Corollary}
\newtheorem{lemma}[theorem]{Lemma}

\newtheorem*{numberlesstheorem}{Theorem}
\newtheorem*{numberlesscorollary}{Corollary}

\newcommand{\set}[1]{\{ #1 \}}         % set braces  {  }

\newcommand{\Diff}{\operatorname{Diff}}
\newcommand{\Imb}{\operatorname{Imb}}

\renewcommand{\P}{\operatorname{{\mathbb P}}}
\newcommand{\R}{\operatorname{{\mathbb R}}}
\newcommand{\Z}{\operatorname{{\mathbb Z}}}

\newcommand{\sA}{\textsc{a}}
\newcommand{\sB}{\textsc{b}}

\newcommand{\sX}{\textsc{x}}
\newcommand{\sY}{\textsc{y}}

\newcommand{\diff}{\operatorname{diff}}
\newcommand{\Exp}{\operatorname{Exp}}
\newcommand{\Fr}{\operatorname{Fr}}
\newcommand{\Isom}{\operatorname{Isom}}
\newcommand{\isom}{\operatorname{isom}}
\newcommand{\vertical}{\operatorname{vert}}
\newcommand{\imb}{\operatorname{imb}}

\newcommand{\Is}{\operatorname{{\mathcal I}}}

\newcommand{\cdim}{\operatorname{cdim}}
\newcommand{\Dih}{\operatorname{Dih}}
\renewcommand{\dim}{\operatorname{dim}}
\newcommand{\Or}{\operatorname{O}}
\renewcommand{\O}{\operatorname{O}}
\newcommand{\SO}{\operatorname{SO}}
\newcommand{\Ostar}{\Or(2)^*}
\newcommand{\ttimes}{\operatorname{\,\widetilde{\times}\,}}

\newcommand{\longpage}{\enlargethispage{\baselineskip}}
\newcommand{\shortpage}{\enlargethispage{-\baselineskip}}

\begin{document}

%\frontmatter
\title[The Smale Conjecture for Lens Spaces]{The Smale Conjecture for Lens
Spaces}

\author{Sungbok Hong}    % first author
\address{\parbox{4in}{Department of Mathematics\\ 
Korea University\\
Seoul 136-701\\ Korea\medskip}}
\email{shong@semi.korea.ac.kr}
\urladdr{http://elie.korea.ac.kr/$_{\widetilde{\phantom{i}}}\,$shong/}
\thanks{Supported in part by the Korea Research Foundation (2002-070-C00012)}

\author{Darryl McCullough}    % second author
\address{\parbox{4in}{Department of Mathematics\\ 
University of Oklahoma\\ Norman, OK 73019\\ U. S. A.\medskip}}
\email{dmccullough@math.ou.edu}
\urladdr{www.math.ou.edu/$_{\widetilde{\phantom{i}}}\,$dmccullough/}
\thanks{Supported in part by the National Science Foundation.}
\author{J. H. Rubinstein}   % third author
\address{\parbox{4in}{Department of Mathematics\\
University of Melbourne\\
Parkville, Victoria 3052\\
Australia\medskip}}
\email{rubin@ms.unimelb.edu.au}
\urladdr{www.ms.unimelb.edu.au/$_{\widetilde{\phantom{n}}}$rubin/}
\thanks{Supported in part by the Australian Research Council.}

\date{\today}
\subjclass[2000]{Primary 57M99; Secondary 57M50, 57R50, 57S05, 58D05}

\maketitle

Throughout, the term \emph{lens space} will mean a $3$-dimensional lens space
$L(m,q)$ other than $L(1,0)$ (the $3$-sphere), $L(0,1)$ (the product
$S^2\times S^1$), and \textit{other than} $L(2,1)$ (the real projective
$3$-space $\R\P^3$). In addition, we always select $q$ so that $1\leq q<m/2$.

Lens spaces are elliptic $3$-manifolds. That is, they may be regarded as
the quotient of the standard round $3$-sphere $S^3$ by a finite subgroup of
the group $\SO(4)$ of orientation-preserving isometries of $S^3$. Then,
they inherit a Riemannian metric of constant positive curvature. For this
metric, the isometry group $\Isom(L(m,q))$ is a Lie group of dimension
either $2$ or $4$. These groups are given in Table~\ref{tab:lens spaces}.

S. Smale \cite{Smale} proved that for the standard round $2$-sphere $S^2$,
the inclusion of the isometry group $\O(3)$ into the diffeomorphism group
$\Diff(S^2)$ is a homotopy equivalence. He conjectured that the analogous
result holds true for the $3$-sphere, that is, that $\O(4)\to \Diff(S^3)$
is a homotopy equivalence. J.~Cerf~\cite{Cerf} proved that the inclusion
induces a bijection on path components, and the full conjecture was proven
by A.~Hatcher \cite{H1}. This is a result of fundamental importance in the
theory of $3$-manifolds, because it shows that smooth structures on
$3$-manifolds are unique up to diffeomorphism, and that for many purposes,
there is no essential difference between the group of homeomorphisms and
the group of diffeomorphisms of a $3$-manifold~\cite{Cerf1}.

A natural extension of Smale's conjecture is that if $M$ is any elliptic
$3$-manifold, then $\Isom(M)\to \Diff(M)$ is a homotopy equivalence.  This
has been proven for some cases \cite{I1, I2, MR}, among them the lens
spaces $L(4n,2n-1)$, $n\geq 2$. Our main result extends this to all
lens spaces:

\begin{numberlesstheorem}[Smale Conjecture for Lens Spaces]
For any lens space $L$, the inclusion $\Isom(L)\to \Diff(L)$ is a homotopy
equivalence.
\end{numberlesstheorem}

One consequence of this is the determination of the homeomorphism type of
$\Diff(L)$. Recall that a Fr\'echet space is a locally convex complete
metrizable linear space (see for example \cite[Proposition 6.4]{BP}, or
\cite{Leslie}). If $M$ is a closed smooth manifold, then with the
$C^\infty$-topology, $\Diff(M)$ is a second countable infinite-dimensional
manifold locally modeled on the Fr\'echet space of smooth vector fields on
$M$ (see for example \cite[section 1.2]{Banyaga} for the local structure,
and for the second countability see chapter~2 of \cite{Hirsch}, especially
section~2.4).  By the Anderson-Kadec Theorem \cite[Corollary VI.5.2]{BP},
every infinite-dimensional separable Fr\'echet space is homeomorphic to
$\R^\infty$, the countable product of lines. A theorem of Henderson and
Schori (\cite[Theorem~IX.7.3]{BP}, originally announced in \cite{HS}) shows
that if $Y$ is any locally convex space with $Y$ homeomorphic to
$Y^\infty$, then manifolds locally modeled on $Y$ are homeomorphic whenever
they have the same homotopy type. Therefore the Smale Conjecture for Lens
Spaces gives immediately the homeomorphism type of $\Diff(L)$:

\begin{numberlesscorollary} For any lens space $L$,
$\Diff(L)$ is homeomorphic to $\Isom(L)\times \R^\infty$.\par
\end{numberlesscorollary}

\noindent
This combines with the known calculations of $\Isom(L)$ in
Table~\ref{tab:lens spaces} to give a complete classification of $\Diff(L)$
for lens spaces into four homeomorphism types. In the following corollary,
we assume as usual that $L(m,q)$ is written with $1\leq q<m/2$, and we
write $P_n$ for the discrete space with $n$ points.
\begin{numberlesscorollary} For a lens space $L(m,q)$, the homeomorphism
type of $\Diff(L)$ is as follows:
\begin{enumerate}
\item For $m$ odd, $\Diff(L(m,1))\approx S^1\times S^3\times P_2$.
\item For $m$ even, $\Diff(L(m,1))\approx S^1\times \SO(3)\times P_2$.
\item For $q>1$ and $q^2\not\equiv \pm 1\pmod{m}$,
$\Diff(L(m,q))\approx S^1\times S^1\times P_2$.
\item For $q>1$ and $q^2\equiv \pm 1\pmod{m}$,
$\Diff(L(m,q))\approx S^1\times S^1\times P_4$.
\end{enumerate}
\end{numberlesscorollary}

Our homeomorphism classification contrasts with the fact that the
\emph{isomorphism type} of $\Diff(L)$ determines $L$. In fact, every
abstract group isomorphism between the diffeomorphism groups of two smooth
manifolds without boundary is induced by a diffeomorphism between the
manifolds \cite{Banyaga, Filip, Takens}.

The authors gratefully acknowledge support from many sources during the
pursuance of this work, including the Australian Research Council, the
Basic Science Research Center of Korea University, the U. S. National
Science Foundation, the University of Oklahoma College of Arts and
Sciences, and the University of Oklahoma Vice President for Research.

\newpage
\section[Outline of the proof]
{Outline of the proof}
\label{sec:outline}

In this section, we will outline the proof of the Smale Conjecture for Lens
Spaces. 

Some initial reductions, detailed in section~\ref{sec:reduction}, reduce
the Smale Conjecture for Lens Spaces to showing that the inclusion
$\diff_f(L)\to \diff(L)$ is an isomorphism on homotopy groups. Here,
$\diff(L)$ is the connected component of the identity in $\Diff(L)$, and
$\diff_f(L)$ is the connected component of the identity in the group of
diffeomorphisms that are fiber-preserving with respect to a Seifert
fibering of $L$ induced from the Hopf fibering of its universal cover,
$S^3$. To simplify the exposition, most of the paper is devoted just to
proving that $\diff_f(L)\to \diff(L)$ is surjective on homotopy groups,
that is, that a map from $S^d$ to $\diff(L)$ is homotopic to a map into
$\diff_f(L)$. The injectivity is obtained in section~\ref{sec:Dk} by a
combination of tricks and minor adaptations of the main program.

Of course, a major difficulty in working with elliptic $3$-manifolds is
their lack of incompressible surfaces. In their place, we use another
structure which has a certain degree of essentiality, called a
\emph{sweepout.} This means a structure on $L$ as a quotient of $P\times
I$, where $P$ is a torus, in which $P\times \{0\}$ and $P\times \{1\}$ are
collapsed to core circles of the solid tori of a genus~1 Heegaard splitting
of $L$. For $0<t<1$, $P\times\{t\}$ becomes a Heegaard torus in $L$,
denoted by $P_t$ and called a \emph{level.} The sweepout is chosen so that
each $P_t$ is a union of fibers. Sweepouts are examined in
section~\ref{sec:RSgraphic}.

Start with a parameterized family of diffeomorphisms $f\colon L\times
S^d\to L$, and for $u\in S^d$ denote by $f_u$ the restriction of $f$ to
$L\times\{u\}$.  The procedure that deforms $f$ to make each $f_u$
fiber-preserving has three major steps.

Step~1 (``finding good levels'') is to perturb $f$ so that for each $u$,
there is some pair $(s,t)$ so that $f_u(P_s)$ intersects $P_t$
transversely, in a collection of circles each of which is either essential
in both $f_u(P_s)$ and $P_t$ (a \emph{biessential} intersection), or
inessential in both (a \emph{discal} intersection), and at least one
intersection circle is biessential. This pair is said to intersect in
\emph{good position,} and if none of the intersections is discal, in
\emph{very good position.}  These concepts are developed in
section~\ref{sec:very good position}, after a preliminary examination of
annuli in solid tori in section~\ref{sec:preliminaries}.

To accomplish Step~1, the methodology of Rubinstein and Scharlemann in
\cite{RS} is adapted. This is reviewed in
section~\ref{sec:Rubinstein-Scharlemann}. First, one perturbs $f$ to be in
``general position,'' as defined in section~\ref{sec:generalposition}. The
intersections of the $f_u(P_s)$ and $P_t$ are then sufficiently
well-controlled to define a \emph{graphic} in the square $I^2$. That is,
the pairs $(s,t)$ for which $f_u(P_s)$ and $P_t$ do not intersect
transversely form a graph imbedded in the square. The complementary regions
of this graph in $I^2$ are labeled according to a procedure in \cite{RS},
and in section~\ref{sec:goodregions} we show that the properties of general
position salvage enough of the combinatorics of these labels developed in
\cite{RS} to deduce that at least one of the complementary regions consists
of pairs in good position.

Perhaps the hardest work of the paper, and certainly the part that takes us
farthest from the usual confines of low-dimensional topology, is the
verification that sufficient ``general position'' can be achieved. Since we
use parameterized families, we must allow $f_u(P_s)$ and $P_t$ to have
large numbers of tangencies, some of which may be of high order. It turns
out that to make the combinatorics of \cite{RS} go through, we must achieve
that at each parameter \textit{there are at most finitely many pairs
$(s,t)$ where $f_u(P_s)$ and $P_t$ have multiple or high-order tangencies}
(at least, for pairs not extremely close to the boundary of the
square). The need for this requirement is illustrated by examples in
section~\ref{sec:examples}, where we construct pairs of sweepouts with all
tangencies of Morse type, but having no pair of levels intersecting in good
position. To achieve the necessary degree of general position, we use
results of a number of people, notably J. W. Bruce \cite{Bruce} and
F. Sergeraert \cite{Sergeraert}.

Step~2 (``from good to very good'') is to deform $f$ to eliminate the
discal intersections of $f_u(P_s)$ and $P_t$, for certain pairs in good
position that have been found in Step~1, so that they intersect in very
good position. This is an application of Hatcher's parameterization
methods~\cite{H}. One must be careful here, since an isotopy that
eliminates a discal intersection can also eliminate a biessential
intersection, and if all biessential intersections were eliminated by the
procedure, the resulting pair would no longer be in very good position.
Lemma~\ref{lem:no biessential elimination} ensures that not all biessential
intersections will be eliminated.

Step~3 (``from very good to fiber-preserving'') is to use the pairs in very
good position to deform the family so that each $f_u$ is fiber-preserving.
This is carried out in sections~\ref{sec:lemmas}
and~\ref{sec:fiber-preserving families}. The basic idea is first to use the
biessential intersections to deform the $f_u$ so that $f_u(P_s)$ actually
equals $P_t$ (for certain $(s,t)$ pairs that orginially intersected in good
position), then use known results about the diffeomorphism groups of
surfaces and Haken $3$-manifolds to make the $f_u$ fiber-preserving on
$P_s$ and then on its complementary solid tori. This process is technically
complicated for two reasons. First, although a biessential intersection is
essential in both tori, it can be contractible in one of the complementary
solid tori of $P_t$, and $f_u(P_s)$ can meet that complementary solid torus
in annuli that are not parallel into $P_t$. So one may be able to push the
annuli out from only one side of $P_t$. Secondly, the fitting together of
these isotopies requires one to work with not just one level but many
levels at a single parameter.

Two natural questions are whether Bonahon's original method for determining
the mapping class group $\pi_0(\Diff(L))$ \cite{Bonahon} can be adapted to
the parameterized setting, and whether our methodology can be used to
recover his results. Concerning the first question, we have had no success
with this approach, as we see no way to perturb the family to the point
where the method can be started at each parameter. For the second, the
answer is yes. In fact, the key geometric step of \cite{Bonahon} is the
isotopy uniqueness of genus-one Heegaard surfaces in $L$, which was already
reproven in Rubinstein and Scharlemann's original
work~\cite[Corollary~6.3]{RS}.

\newpage
\section[Reductions]
{Reductions}
\label{sec:reduction}

In this section, we carry out initial reductions. The Conjecture
will be reduced to a purely topological problem of deforming
parameterized families of diffeomorphisms to families of
diffeomorphisms that preserve a certain Seifert fibering of $L$, which we
call the Hopf fibering.

The paper \cite{M} contains a calculation of the isometry groups of all
elliptic $3$-manifolds (calculations for lens spaces were also given in
\cite{HR} and \cite{KM1}). Among the elliptic $3$-manifolds, the lens
spaces have the most complicated isometry groups, given in
Table~\ref{tab:lens spaces}.  In the table, $\isom(L(m,q))$ is the path
component of the identity map in the isometry group $\Isom(L(m,q))$, and
$\Is(L(m,q))$ is the group of path components of $\Isom(L(m,q))$.  The
orthogonal groups are denoted by $\Or(4)$, $\SO(3)$ and $\Or(2)$, $C_k$ is
the cyclic group of order $k$, and $\Dih(S^1\times S^1)$ is the semidirect
product $(S^1\times S^1)\circ C_2$, where $C_2$ acts by complex conjugation
in both factors. Also, $\Ostar$ is the nontrivial central extension of
$\Or(2)$ by $C_2$, that is, the preimage of $\Or(2)\subset \SO(3)$ under
the $2$-fold covering map $S^3\to \SO(3)$. If $H_1$ and $H_2$ are groups,
each containing $-1$ as a central involution, then the quotient $(H_1\times
H_2)/\langle (-1,-1)\rangle$ is denoted by $H_1\ttimes H_2$. In particular,
$\SO(4)$ itself is $S^3\ttimes S^3$, and contains the subgroups
$\Ostar\ttimes S^3$ and $S^1\ttimes S^1$. The latter is isomorphic to
$S^1\times S^1$.\par
\begin{table}
\begin{small}
\renewcommand{\arraystretch}{1.5}
\newlength{\minipagewidth}%
\setlength{\tabcolsep}{1.5 ex}
\setlength{\fboxsep}{0pt}
\fbox{%
\begin{tabular}{l|l|l}
$m$, $q$&$\isom(L(m,q))$&$\Is(L(m,q))$\\
\hline
\hline
$m=1$&$\Or(4)$&$C_2$\\
\hline
$m=2$&$\SO(3)\times \SO(3)$&$C_2$\\
\hline
\settowidth{\minipagewidth}{$m>2$, $m$ even}%
\begin{minipage}{\minipagewidth}%
\noindent $m>2$, $m$ odd,\par\end{minipage} $q=1$&$\Ostar\ttimes S^3$&$C_2$\\  
\hline
$m>2$, $m$ even, $q=1$&$\Or(2)\times \SO(3)$&$C_2$\\
\hline
$m>2$, $1<q<m/2$, $q^2\not\equiv\pm1\bmod{m}$&$\Dih(S^1\times S^1)$&$C_2$\\
\hline
$m>2$, $1<q<m/2$, $q^2\equiv-1\bmod{m}$&$S^1\ttimes S^1$&$C_4$\\
\hline
$m>2$,
%\settowidth{\minipagewidth}{$1<q<m/2$, $q^2\equiv1\pod{m}$,}%
\settowidth{\minipagewidth}{$\gcd(m,q+1)\gcd(m,q-1)=m$}%
\begin{minipage}[t]{\minipagewidth}%
\noindent $1<q<m/2$, $q^2\equiv1\bmod{m}$,\par
\noindent $\gcd(m,q+1)\gcd(m,q-1)=m$\rule[-1.2 ex]{0mm}{0mm}\par%
\end{minipage}%
&$\Or(2)\ttimes \Or(2)$&$C_2\times C_2$\\
\hline
$m>2$,
%\settowidth{\minipagewidth}{$1<q<m/2$, $q^2\equiv1\pod{m}$,}%
\settowidth{\minipagewidth}{$\gcd(m,q+1)\gcd(m,q-1)=2m$}%
\begin{minipage}[t]{\minipagewidth}%
\noindent $1<q<m/2$, $q^2\equiv1\bmod{m}$,\par
\noindent $\gcd(m,q+1)\gcd(m,q-1)=2m$\rule[-1.2 ex]{0mm}{0mm}\par%
\end{minipage}%
&$\Or(2)\times \Or(2)$&$C_2\times C_2$
\end{tabular}}
\end{small}
\bigskip
\caption{Isometry groups of $L(m,q)$}
\label{tab:lens spaces}
\end{table}

From Table~\ref{tab:lens spaces}, one sees that $\isom(L(m,1))$ is
homeomorphic to $S^1\times S^3$ for $m$ odd, and to $S^1\times \SO(3)$ for
$m$ even, while for $q>1$, $\isom(L(m,q))$ is homeomorphic to $S^1\times
S^1$.  These observations were used in the second corollary stated in the
introduction.

The following theorem from \cite{M} is the ``$\pi_0$-part'' of the Smale
Conjecture for elliptic $3$-manifolds.
\begin{theorem}
Let $M$ be an elliptic $3$-manifold. Then the inclusion of $\Isom(M)$ into
$\Diff(M)$ is a bijection on path components.
\label{thm:pi0 Smale}
\end{theorem}
\noindent Consequently, to prove the Smale Conjecture for a lens space $L$,
it is sufficient to prove that the inclusion of the connected components of
the identity map $\isom(L)\to\diff(L)$ is a homotopy equivalence. 

Since $\Diff(M)$ is an infinite-dimensional manifold locally modeled on
$\R^\infty$, Corollary~IX.7.1 of \cite{BP} (originally theorem~4 of
\cite{Henderson}) shows that $\Diff(M)$ admits an open imbedding into
$\R^\infty$. Theorems II.6.2 and II.6.3 of \cite{BP} then show that
$\Diff(M)$ has the homotopy type of a CW-complex (as far as we know, this
fact is due originally to Palais \cite{Palais}). So $\diff(M)$ has the
homotopy type of a CW-complex, and the same is true for $\isom(M)$, since
it is a manifold. Therefore it is sufficient to prove that
$\isom(L)\to\diff(L)$ is a weak homotopy equivalence.

Section~1.4 of \cite{M} gives a certain way to imbed $\pi_1(L)$ into
$\SO(4)$ so that its action on $S^3$ is fiber-preserving for the
fibers of the Hopf bundle structure of $S^3$. Consequently, this bundle
structure descends to a Seifert fibering of $L$, which we call the
\emph{Hopf fibering} of $L$. If $q=1$, this Hopf fibering is actually an
$S^1$-bundle structure, while if $q>1$, it has two exceptional fibers with
invariants of the form $(k,q_1)$, $(k,q_2)$ where $k=m/\gcd(q-1,m)$ (see
Table~4 of \cite{M}). We will always use the Hopf fibering as the
Seifert-fibered structure of $L$.

A diffeomorphism from $L$ to $L$ is called \emph{fiber-preserving} if the
image of each fiber is a fiber, and \emph{vertical} if it preserves each
fiber. By $\diff_f(L)$ we denote the connected component of the identity
map in the space of fiber-preserving diffeomorphisms.  Theorem~2.1 of
\cite{M} shows that (since $m>2$) every orientation-preserving isometry of
$L$ preserves the Hopf fibering on $L$. In particular, $\isom(L)\subset
\diff_f(L)$, so there are inclusions
\[\isom(L)\to \diff_f(L)\to \diff(L)\ .\]

\begin{theorem}
The inclusion $\isom(L)\to \diff_f(L)$ is a weak homotopy equivalence.
\label{thm:reduce to fiber-preserving}
\end{theorem}

\begin{proof}
The argument is similar to the latter part of the proof of theorem~4.2 of
\cite{MR}, so we only give a sketch. There is a diagram
\begin{equation*}
\begin{CD}
S^1 @>>> \isom(L) @>>>  \isom(L_0)\\ 
@VVV @VVV @VVV { }\\
\vertical(L) @>>>  \diff_f(L) @>>>  \diff_{orb}(L_0)
\end{CD}
\end{equation*}
\noindent
where $L_0$ is the quotient orbifold and $\diff_{orb}(L_0)$ is the group of
orbifold diffeomorphisms of $L_0$, and $\vertical(L)$ is the group of
vertical diffeomorphisms. The first row is a fibration, in fact an
$S^1$-bundle, and the second row is a fibration by theorem~8.3 of
\cite{KM}. The vertical arrows are inclusions. When $q=1$, $L_0$ is the
$2$-sphere and the right-hand vertical arrow is the inclusion of $\SO(3)$
into $\diff(S^2)$, which is a homotopy equivalence by \cite{Smale}. When
$q>1$, $L_0$ is a $2$-sphere with two cone points, $\isom(L_0)$ is
homeomorphic to $S^1$, and $\diff_{orb}(L_0)$ is essentially the connected
component of the identity in the diffeomorphism group of the annulus. Again
the right-hand vertical arrow is a weak homotopy equivalence. The left-hand
vertical arrow is a weak homotopy equivalence in both cases, so the middle
arrow is as well.
\end{proof}

Theorem~\ref{thm:reduce to fiber-preserving} reduces the Smale Conjecture
for Lens Spaces to proving that the inclusion $\diff_f(L)\to \diff(L)$ is a
weak homotopy equivalence. For this it is sufficient to prove that for all
$d\geq 1$, any map $f\colon (D^d,S^{d-1})\to (\diff(L),\diff_f(L))$ is
homotopic, through maps taking $S^{d-1}$ to $\diff_f(L)$, to a map from
$D^d$ into $\diff_f(L)$. To simplify the exposition, we work
until the final section with a map $f\colon S^d\to \diff(L)$ and show
that it is homotopic to a map into $\diff_f(L)$. In the final section, we
give a trick that enables the entire procedure to be adapted to maps
$f\colon (D^d,S^{d-1})\to (\diff(L),\diff_f(L))$, completing the proof.

\newpage
\section[Annuli in solid tori]
{Annuli in solid tori}
\label{sec:preliminaries}

Annuli in solid tori will appear frequently in our work. Incompressible
annuli present little difficulty, but we will also need to examine
compressible annuli, whose behavior is somewhat more complicated. In this
section, we provide some basic definitions and lemmas.

A loop $\alpha$ in a solid torus $V$ is called a \emph{longitude} if its
homotopy class is a generator of the infinite cyclic group $\pi_1(V)$.  If
in addition there is a product structure $V=S^1\times D^2$ for which
$\alpha=S^1\times\set{0}$, then $\alpha$ is called a \emph{core circle} of
$V$. A subset of a solid torus $V$ is called a \emph{core region} when it
contains a core circle of $V$. An imbedded circle in $\partial V$ which is
essential in $\partial V$ and contractible in $V$ is called a
\emph{meridian} of $V$; a properly imbedded disk in $V$ whose boundary is a
meridian is called a \emph{meridian disk} of $V$.

Annuli in solid tori will always be assumed to be properly imbedded, which
for us includes the property of being transverse to the boundary, unless
they are actually contained in the boundary. The next three results are
elementary topological facts, and we do not include proofs.
\begin{proposition} 
Let $A$ be a boundary-parallel annulus in a solid torus $V$, which
separates $V$ into $V_0$ and $V_1$, and for $i=0,1,$ let $A_i=V_i\cap
\partial V$.  Then $A$ is parallel to $A_i$ if and only if $V_{1-i}$ is a
core region.
\label{prop:boundary-parallel annuli}
\end{proposition}
\noindent 

\begin{proposition} 
Let $A$ be a properly imbedded annulus in a solid torus $V$, which
separates $V$ into $V_0$ and $V_1$, and let $A_i=V_i\cap \partial V$. The
following are equivalent:
\begin{enumerate}
\item
$A$ contains a longitude of $V$.
\item
$A$ contains a core circle of $V$.
\item
$A$ is parallel to both $A_0$ and $A_1$.
\item
Both $V_0$ and $V_1$ contain longitudes of $V$.
\item
Both $V_0$ and $V_1$ are core regions of $V$.
\end{enumerate}
\label{prop:longitudinal annuli}
\end{proposition}
\noindent An annulus satisfying the conditions in
proposition~\ref{prop:longitudinal annuli} is said to be
\emph{longitudinal.} A longitudinal annulus must be incompressible.

\begin{proposition} Let $V$ be a solid torus and let $\cup A_i$ be a union 
of disjoint boundary-parallel annuli in $V$. Let $C$ be a core circle of
$V$ that is disjoint from $\cup A_i$. For each $A_i$, let $V_i$ be the
closure of the complementary component of $A_i$ that does not contain $C$,
and let $B_i=V_i\cap \partial V$. Then $A_i$ is parallel to $B_i$.
Furthermore, either
\begin{enumerate}
\item
no $A_i$ is longitudinal, and exactly one component of $V-\cup A_i$
is a core region, or
\item
every $A_i$ is longitudinal, and every component of $V-\cup A_i$ is a core
region.
\end{enumerate}
\label{prop:core circles}
\end{proposition}

\begin{figure}
\includegraphics[width=40 ex]{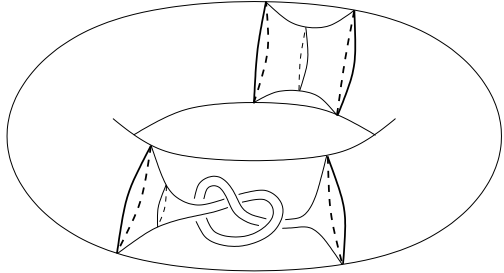}
\caption{Meridional annuli in a solid torus.}
\label{fig:meridional annuli}
\end{figure}

There are various kinds of compressible annuli in solid tori. For example,
there are boundaries of regular neighborhoods of properly imbedded arcs,
possibly knotted. Also, there are annuli with one boundary circle a
meridian and the other a contractible circle in the boundary torus. When
both boundary circles are meridians, we call the annulus \emph{meridional.}
As shown in figure~\ref{fig:meridional annuli}, meridional annuli are not
necessarily boundary-parallel.

Although meridional annuli need not be boundary-parallel, they behave
homologically as though they were, and as a consequence any family of
meridional annuli misses some longitude.
\begin{lemma}
Let $A_1,\ldots\,$, $A_n$ be disjoint meridional annuli in a solid torus
$V$. Then:
\begin{enumerate}
\item
Each $A_i$ separates $V$ into two components, $V_{i,0}$ and $V_{i,1}$, for
which $A_i$ is incompressible in $V_{i,0}$ and compressible in $V_{i,1}$.
\item $V_{i,1}$ contains a meridian disk of $V$.
\item $\pi_1(V_{i,0})\to \pi_1(V)$ is the zero homomorphism.
\item The intersection of the $V_{i,1}$ is the unique
component of the complement of $\cup A_i$ that contains a longitude of $V$.
\end{enumerate}
\label{lem:meridional annuli}
\end{lemma}

\begin{proof}
For each $i$, every loop in $V$ has even algebraic intersection with $A_i$,
since every loop in $\partial V$ does, so $A_i$ separates $V$. Since $A_i$
is not incompressible, it must be compressible in one of its complementary
components, $V_{i,1}$, and since $V$ is irreducible, $A_i$ must be
incompressible in the other complementary component, $V_{i,0}$.

Notice that $V_{i,1}$ must contain a meridian disk of $V$. Indeed, if $K$ is
the union of $A_i$ with a compressing disk in $V_{i,1}$, then two of the
components of the frontier of a regular neighborhood of $K$ in $V$ are
meridian disks of $V_{i,1}$. Consequently, $\pi_1(V_{i,0})\to \pi_1(V)$ is
the zero homomorphism. The Mayer-Vietoris sequence shows that $H_1(A)\to
H_1(V_{i,0})$ and $H_1(V_{i,1})\to H_1(V)$ are isomorphisms.

Let $V_1$ be the intersection of the $V_{i,1}$, and let $V_0$ be the union
of the $V_{i,0}$. The Mayer-Vietoris sequence shows that $V_1$ is
connected, and that $H_1(V_1)\to H_1(V)$ is an isomorphism, so $V_1$
contains a longitude of $V$. No other complementary component of $\cup A_i$
contains a longitude, since each such component lies in one of the
$V_{i,0}$ and all of its loops must be contractible in~$V$.
\end{proof}

\newpage
\section[Heegaard tori in very good position]
{Heegaard tori in very good position}
\label{sec:very good position}

A \emph{Heegaard torus} in a lens space $L$ is a torus that separates $L$
into two solid tori.  In this section we will develop some properties of
Heegaard tori.  Also, we introduce the concepts of discal and biessential
intersection circles, good position, and very good position, which will be
used extensively in later sections.

When $P$ is a Heegaard torus bounding solid tori $V$ and $W$, and $Q$ is a
Heegaard torus contained in the interior of $V$, $Q$ need not be parallel
to $\partial V$. For example, start with a core circle in $V$, move a small
portion of it to $\partial V$, then pass it across a meridian disk of $W$
and back into $V$.  This moves the core circle to its band-connected sum in
$V$ with an $(m,q)$-curve in $\partial V$. By varying the choice of band---
for example, by twisting it or tying knots in it--- and by iterating this
construction, one can construct complicated knotted circles in $V$ which
are isotopic in $L$ to a core circle of $V$. The boundary of a regular
neighborhood of such a circle is a Heegaard torus of $L$. But here is one
restriction on Heegaard tori:

\begin{proposition}
Let $P$ be a Heegaard torus in a lens space $L$, bounding solid tori $V$
and $W$. If a loop $\ell$ imbedded in $P$ is a core circle for a solid
torus of some genus-$1$ Heegaard splitting of $L$, then $\ell$ is a
longitude for either $V$ or $W$.
\label{prop:Heegaard loops}
\end{proposition}

\begin{proof}
Since $L$ is not simply-connected, $\ell$ is not a meridian for either $V$
or $W$, consequently $\pi_1(\ell)\to \pi_1(V)$ and $\pi_1(\ell)\to
\pi_1(W)$ are injective. So $P-\ell$ is an open annulus separating
$L-\ell$, making $\pi_1(L-\ell)$ a free product with amalgamation
$\Z*_{\Z}\Z$. Since $\ell$ is a core circle, $\pi_1(L-\ell)$ is infinite
cyclic, so at least one of the inclusions of the amalgamating subgroup to
the infinite cyclic factors is surjective.
\end{proof}

Let $F_1$ and $F_2$ be transversely intersecting imbedded surfaces in the
interior of a $3$-manifold $M$. A component of $F_1\cap F_2$ is called
\emph{discal} when it is contractible in both $F_1$ and $F_2$, and
\emph{biessential} when it is essential in both. We say that $F_1$ and
$F_2$ are \emph{in good position} when every component of their
intersection is either discal or biessential, and at least one is
biessential, and we say that they are \emph{in very good position} when
they are in good position and every component of their intersection is
biessential.

\begin{figure}
\includegraphics[width=45 ex]{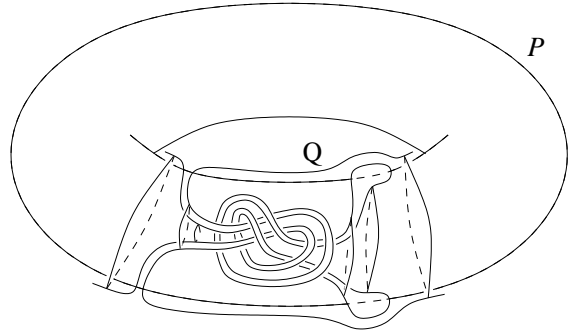}
\caption{Heegard tori in very good position, with non-boundary-parallel
meridional annuli.}
\label{fig:bad torus}
\end{figure}

Later, we will go to considerable effort to obtain pairs Heegaard tori for
lens spaces that intersect in very good position. Even then, the
configuration can be complicated. Consider a Heegaard torus $P$ bounding
solid tori $V$ and $W$, and another Heegaard torus $Q$ that meets $P$ in
very good position. When the intersection circles are not meridians for
either $V$ or $W$, the components of $Q\cap V$ and $Q\cap W$ are annuli
that are incompressible in $V$ and $W$, and must be as described in
proposition~\ref{prop:core circles}. But if the intersection circles are
meridians for one of the solid tori, say $V$, then $Q\cap V$ consists of
meridional annuli, and as shown in figure~\ref{fig:bad torus}, they need
not be boundary-parallel. To obtain that configuration, one starts with a
torus $Q$ parallel to $P$ and outside $P$, and changes $Q$ only by an
isotopy on a regular neighborhood of a meridian $c$ of $Q$.  First, $c$
passes across a meridian in $P$, then shrinks down to a small circle which
traces around a knot. Then, it expands out to another meridian in $P$ and
pushes across. The resulting torus meets $P$ in four circles which are
meridians for $V$, and meets $V$ in two annuli, both isotopic to the
non-boundary-parallel annulus in figure~\ref{fig:meridional annuli}. The
next lemma gives a small but important restriction on meridional annuli of
$Q\cap V$.

\begin{lemma}
Let $P$ be a Heegaard torus which separates a lens space into two solid
tori $V$ and $W$.  Let $Q$ be another Heegaard torus whose intersection
with $V$ consists of a single meridional annulus $A$. Then $A$ is
boundary-parallel in $V$.
\label{lem:one annulus}
\end{lemma}

\begin{proof}
From lemma~\ref{lem:meridional annuli}, $A$ separates $V$ into two
components $V_0$ and $V_1$, such that $A$ is compressible in $V_1$ and
$V_1$ contains a longitude of $V$. Suppose that $A$ is not
boundary-parallel in~$V$.

Let $A_0=V_0\cap \partial V$. Of the two solid tori in $L$ bounded by $Q$,
let $X$ be the one that contains $A_0$, and let $Y$ be the other one. Since
$Q\cap V$ consists only of $A$, $Y$ contains $V_1$, and in particular
contains a compressing disk for $A$ in $V_1$ and a longitude for $V$.

Suppose that $A_0$ were incompressible in $X$. Since $A_0$ is not
parallel to $A$, it would be parallel to $\overline{\partial
X-A}$. So $V_0$ would contain a core circle of $X$. Since $\pi_1(V_0)\to
\pi_1(V)$ is the zero homomorphism, this implies that $L$ is
simply-connected, a contradiction. So $A_0$ is compressible in $X$. A
compressing disk for $A_0$ in $X$ is part of a $2$-sphere that meets $Y$
only in a compressing disk of $A$ in $V_1$. This $2$-sphere has algebraic
intersection $\pm 1$ with the longitude of $V$ in $V_1$,
contradicting the irreducibility of $L$.
\end{proof}

Regarding $D^2$ as the unit disk in the plane, for $0<r<1$ let $rD^2$
denote $\{(x,y)\;|\;x^2+y^2\leq r^2\}$.  A solid torus $X$ imbedded in a
solid torus $V$ is called \emph{concentric in $V$} if there is some product
structure $V=D^2\times S^1$ such that $X=rD^2\times S^1$. Equivalently, $X$
is in the interior of $V$ and some (hence every) core circle of $X$ is a
core circle of $V$.

The next lemma shows how we will use Heegard tori that meet in very good
position.
\begin{lemma}
Let $P$ be a Heegaard torus which separates a lens space into two solid
tori $V$ and $W$.  Let $Q$ be another Heegaard torus, that meets $P$ in
very good position, and assume that the annuli of $Q\cap V$ are
incompressible in $V$. Then at least one component $C$ of $V-(Q\cap V)$
satisfies both of the following:
\begin{enumerate}
\item
$C$ is a core region for $V$.
\item
Suppose that $Q$ is moved by isotopy to a torus $Q_1$ in $W$, by pushing
the annuli of $Q\cap V$ one-by-one out of $V$ using isotopies that move
them across regions of $V-C$, and let $X$ be the solid torus bounded by
$Q_1$ that contains $V$. Then $V$ is concentric in $X$.
\item
After all but one of the annuli have been pushed out of $V$, the
image $Q_0$ of $Q$ is isotopic to $P$ relative to $Q_0\cap P$.
\end{enumerate}
\label{lem:pushout}
\end{lemma}

\begin{proof}
Assume first that $Q\cap V$ has only one component $A$. Then $\partial A$
separates $P$ into two annuli, $A_1$ and $A_2$.  Since $A$ is
incompressible in $V$, it is parallel in $V$ to one of the $A_i$, say
$A_1$. Let $A'=Q\cap W$.

If $A'$ is longitudinal, then $A'$ is parallel in $W$ to $A_2$. So pushing
$A$ across $A_1$ moves $Q$ to a torus in $W$ parallel to $P$, and the lemma
holds, with $C$ being the region between $A$ and $A_2$. An isotopy from $Q$
to $P$ can be carried out relative to $Q\cap P$, giving the last statement
of the lemma. Suppose that $A'$ is not longitudinal. If $A'$ is
incompressible, then it is boundary parallel in $W$. If $A'$ is not
incompressible, then since $P$ and $Q$ meet in very good position, $A'$ is
meridional, and by lemma~\ref{lem:one annulus} it is again
boundary-parallel in $W$. If $A'$ is parallel to $A_2$, then we are
finished as before. If $A'$ is parallel to $A_1$, but not to $A_2$, then
there is an isotopy moving $Q$ to a regular neighborhood of a core circle
of $A_1$. By proposition~\ref{prop:Heegaard loops}, $A$ is longitudinal, so
must also be parallel in $V$ to $A_2$. In this case, we take $C$ to be the
region between $A$ and~$A_1$.

Suppose now that $Q\cap V$ and hence also $Q\cap W$ consist of $n$ annuli,
where $n>1$. By isotopies pushing outermost annuli in $V$ across $P$, we
obtain $Q_0$ with $Q_0\cap V$ consisting of one annulus $A$. At least one
of its complementary components, call it $C$, satisfies the lemma. Let $Z$
be the union of the regions across which the $n-1$ annuli were
pushed. Since $C$ is a core region, $C\cap (V-Z)$ is also a core region
(since a core circle of $V$ in $C$ can be moved, by the reverse of the
pushout isotopies, to a core circle of $V$ in $C\cap (V-Z)$).  So $C\cap
(V-Z)$ satisfies the conclusion of the lemma.
\end{proof}

Here is a first consequence of lemma~\ref{lem:pushout}.
\begin{corollary}
Let $P$ be a Heegaard torus which separates a lens space into two solid
tori $V$ and $W$, and let $Q$ be another Heegaard torus separating it into
$X$ and $Y$. Assume that $Q$ meets $P$ in very good position. If the
circles of $P\cap Q$ are meridians (respectively, longitudes) in $X$ or in
$Y$, then they are meridians (longitudes) in $V$ or in~$W$.\par
\label{coro:comeridian}
\end{corollary}

\begin{proof}
We may choose notation so that the annuli of $Q\cap V$ are incompressible
in $V$.  Use lemma~\ref{lem:pushout} to move $Q$ out of $V$. After all but
one annulus has been pushed out, the image $Q_0$ of $Q$ is isotopic to $P$
relative to $Q_0\cap P$. That is, the original $Q$ is isotopic to $P$ by an
isotopy relative to $Q_0\cap P$. If the circles of $Q\cap P$ were
originally meridians of $X$ or $Y$, then in particular those of $Q_0\cap P$
are meridians of $X$ or $Y$ after the isotopy, that is, of $V$ or $W$.  The
case of longitudes is similar.
\end{proof}

\newpage
\section[Sweepouts, and levels in very good positions]
{Sweepouts, and levels in very good position}
\label{sec:RSgraphic}

In this section we will define sweepouts and related structures. Also, we
will prove an important technical lemma concerning pairs of sweepouts
having levels that meet in very good position.

By a \emph{sweepout} of a closed orientable $3$-manifold, we mean
smooth map $\tau\colon P\times [0,1]\to M$, where $P$ is a closed
orientable surface, such that
\begin{enumerate}
\item
$T_0=\tau(P\times\set{0})$ and $T_1=\tau(P\times \set{1})$ are disjoint
graphs with each vertex of valence $3$.
\item 
Each $T_i$ is a union of a collection of smoothly imbedded arcs and circles
in $M$.
\item
$\tau\vert_{P\times(0,1)}\colon P\times (0,1)\to M$ is a diffeomorphism
onto $M-(T_0\cup T_1)$.
\item
Near $P\times\partial I$, $\tau$ gives a mapping cylinder neighborhod of
$T_0\cup T_1$.
\end{enumerate}
Associated to any $t$ with $0<t<1$, there is a Heegaard splitting
$M=V_t\cup W_t$, where $V_t=\tau(P\times [0,t])$ and $W_t=\tau(P\times
[t,1])$. For each $t$, $T_0$ is a deformation retract of $V_t$ and $T_1$ is
a deformation retract of $W_t$.  We denote $\tau(P\times \set{t})$ by
$P_t$, and call it a \emph{level} of $\tau$. Also, for
$0<s< t<1$ we denote the closure of the region between $P_s$ and $P_t$
(that is, $\tau(P\times [s,t])$) by $R(s,t)$.  Note that any genus-$1$
Heegaard splitting of $L$ provides sweepouts with $T_0$ and $T_1$ as core
circles of the two solid tori, and the Heegaard torus as one of the levels.

A sweepout $\tau\colon P\times [0,1]\to M$ induces a continuous projection
function $\pi\colon M\to [0,1]$ by the rule $\pi(\tau(x,t))=t$.  By
composing this with a smooth bijection from $[0,1]$ to $[0,1]$ all of whose
derivatives vanish at $0$ and at $1$, we may reparameterize $\tau$ to
ensure that $\pi$ is a smooth map. We always assume that $\tau$ has been
selected to have this property.

The next lemma gives very strong restrictions on levels of two different
sweepouts of a lens space that intersect in very good position.  For its
proof, recall that a \emph{spine} for a connected surface $P$ is a
$1$-dimensional cell complex in $P$ whose complement consists of open
disks.

\begin{figure}
\includegraphics[width=45 ex]{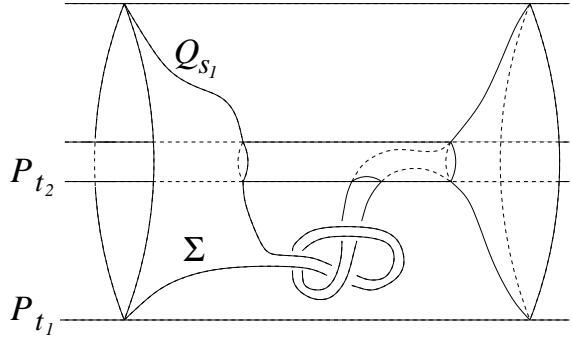}
\caption{Case (2) of lemma \ref{lem:hitting levels}}
\label{fig:sigma}
\end{figure}

\begin{lemma}
Let $L$ be a lens space. Let $\tau \colon T\times[0,1]\to L$ be a sweepout
as above, where $T$ is a torus. Let $\sigma\colon T\times[0,1]\to L$ be
another sweepout, with levels $Q_s=\sigma(T\times\set{s})$.  Suppose that
for $t_1<t_2$, $s_1\neq s_2$, and $i=1,2$, $Q_{s_i}$ and $P_{t_i}$
intersect in very good position, and that $Q_{s_1}$ has no discal
intersections with $P_{t_2}$. If $Q_{s_1}$ has nonempty intersection with
$P_{t_2}$, then either
\begin{enumerate}
\item
every intersection circle of $Q_{s_1}$ with $P_{t_2}$ is biessential, and
consequently $Q_{s_1}\cap R(t_1,t_2)$ contains an annulus with one boundary
circle essential in $P_{t_1}$ and the other essential in $P_{t_2}$, or
\item
for $i=1,2$, $Q_{s_i}\cap P_{t_i}$ consists of meridians of $W_{t_i}$, and
$Q_{s_1}\cap R(t_1,t_2)$ contains a surface $\Sigma$ which is a homology
from a circle of $Q_{s_1}\cap P_{t_1}$ to a union of circles in $P_{t_2}$.
\end{enumerate}
\label{lem:hitting levels}
\end{lemma}
\noindent Figure~\ref{fig:sigma} illustrates case (2) of
lemma~\ref{lem:hitting levels}.

We mention that that to apply lemma~\ref{lem:hitting levels} when
$t_1>t_2$, we interchange the roles of $V_{t_i}$ and $W_{t_i}$. The
intersection circles in case (2) are then meridians of the $V_{t_i}$
rather than the $W_{t_i}$.

\begin{proof}[Proof of lemma~\ref{lem:hitting levels}]
Assume for now that the circles of $Q_{s_2}\cap P_{t_2}$ are not meridians
of $W_{t_2}$.

We first rule out the possibility that there exists a circle of
$Q_{s_1}\cap P_{t_2}$ that is inessential in $Q_{s_1}$. If so, there would
be a circle $C$ of $Q_{s_1}\cap P_{t_2}$, bounding a disk $D$ in $Q_{s_1}$
with interior disjoint from $P_{t_2}$. Since $Q_{s_1}$ and $P_{t_2}$ have
no discal intersections, $C$ is essential in $P_{t_2}$, so $D$ is a
meridian disk for $V_{t_2}$ or $W_{t_2}$. It cannot be a meridian disk for
$V_{t_2}$, for then some circle of $D\cap P_t$ would be a meridian of
$V_{t_1}$, contradicting the fact that $Q_{s_1}$ and $P_{t_1}$ meet in very
good position. But $D$ cannot be a meridian disk for $W_{t_2}$, since $D$
is disjoint from $Q_{s_2}$ and the circles of $Q_{s_2}\cap P_{t_2}$ are not
meridians of~$W_{t_2}$.

We now rule out the possiblity that there exists a circle of $Q_{s_1}\cap
P_{t_2}$ that is essential in $Q_{s_1}$ and inessential in $P_{t_2}$.
There is at least one biessential intersection circle of $Q_{s_1}$ with
$P_{t_1}$, hence also an annulus $A$ in $Q_{s_1}$ with one boundary circle
inessential in $P_{t_2}$ and the other essential in either $P_{t_1}$ or
$P_{t_2}$, with no intersection circle of the interior of $A$ with
$P_{t_1}\cup P_{t_2}$ essential in $A$.  The interior of $A$ must be
disjoint from $P_{t_1}$, since $Q_{s_1}$ meets $P_{t_1}$ in very good
position. It must also be disjoint from $P_{t_2}$, by the previous
paragraph. So, since $A$ has at least one boundary circle in $P_{t_2}$, it
is properly imbedded either in $R(t_1,t_2)$ or in $W_{t_2}$. It cannot be
in $R(t_1,t_2)$, since it has one boundary circle inessential in $P_{t_2}$
and the other essential in $P_{t_1}\cup P_{t_2}$. So $A$ is in $W_{t_2}$,
and since one boundary circle is inessential in $P_{t_2}$, the other must
be a meridian, contradicting the assumption that no circle of $Q_{s_2}\cap
P_{t_2}$ is a meridian of $W_{t_2}$. Thus conclusion~(1) holds when cirlces
of $Q_{s_2}\cap P_{t_2}$ are not meridians of~$W_{t_2}$.

Assume now that the circles of $Q_{s_2}\cap P_{t_2}$ are meridians
of~$W_{t_2}$. We will achieve conclusion~(2).

Suppose first that some circle of $Q_{s_1}\cap P_{t_2}$ is essential in
$Q_{s_1}$. Then there is an annulus $A$ in $Q_{s_1}$ with one boundary
circle essential in $P_{t_1}$, the other essential in $P_{t_2}$, and all
intersections of the interior of $A$ with $P_{t_1}\cup P_{t_2}$ inessential
in $A$. Since $Q_{s_1}$ meets $P_{t_1}$ in very good position, the interior
of $A$ must be disjoint from $P_{t_1}$. So $A\cap R(t_1,t_2)$ contains a
planar surface $\Sigma$ with one boundary component a circle of
$Q_{s_1}\cap P_{t_1}$ and the other boundary components circles in
$P_{t_2}$ which are meridians in $W_{t_2}$, giving the conclusion~(2) of
the lemma.

Suppose now that every circle of $Q_{s_1}\cap P_{t_2}$ is contractible in
$Q_{s_1}$. We will show that this case is impossible. An intersection
circle innermost on $Q_{s_1}$ bounds a disk $D$ in $Q_{s_1}$ which is a
meridian disk for $W_{t_2}$, since $\partial D$ is essential in $P_{t_2}$
and disjoint from $Q_{s_2}\cap P_{t_2}$.  Now, use lemma~\ref{lem:pushout}
to push $Q_{s_2}\cap V_{t_2}$ out of $V_{t_2}$ by an ambient isotopy of
$L$.  Suppose for contradiction that one of these pushouts, say, pushing an
annulus $A_0$ in $Q_{s_2}$ across an annulus in $P_{t_2}$, also eliminates a
circle of $Q_{s_1}\cap P_{t_1}$. Let $Z$ be the region of parallelism
across which $A_0$ is pushed. Since $Z$ contains an essential loop of
$Q_{s_1}$, and each circle of $Q_{s_1}\cap P_{t_2}$ is contractible in
$Q_{s_1}$, $Z$ contains a spine of $Q_{s_1}$. This spine is isotopic in $Z$
into a neighborhood of a boundary circle of $A_0$. Since this boundary
circle is a meridian of $W_{t_2}$, every circle in the spine is
contractible in $L$. This contradicts the fact that $Q_{s_1}$ is a Heegaard
torus. So the pushouts do not eliminate intersections of $Q_{s_1}$ with
$P_{t_1}$, and after the pushouts are completed, the image of $Q_{s_1}$
still meets $P_{t_1}$.

During the pushouts, some of the intersection circles of $Q_{s_1}$ with
$P_{t_2}$ may disappear, but not all of them, since the pushouts only move
points into $W_{t_2}$. So after the pushouts, there is a circle of
$Q_{s_1}\cap P_{t_2}$ that bounds a innermost disk in $Q_{s_1}$ (since all
the original intersection circles of $Q_{s_1}$ with $P_{t_2}$ bound disks
in $Q_{s_1}$, and the new intersection circles are a subset of the old
ones). Since the boundary of this disk is a meridian of $W_{t_2}$, the disk
it bounds in $Q_{s_1}$ must be a meridian disk of $W_{t_2}$. The image of
$Q_{s_2}$ lies in $W_{t_2}$ and misses this meridian disk, contradicting
the fact that $Q_{s_2}$ is a Heegaard torus.
\end{proof}

\newpage
\section[Morse general position and the Rubinstein-Scharlemann graphic]
{The Rubinstein-Scharlemann graphic}
\label{sec:Rubinstein-Scharlemann}

The purpose of this section is to present a number of definitions, and to
sketch the proof of theorem~\ref{thm:RS} below, originally from~\cite{RS}.
It requires the hypothesis that two sweepouts meet in general position in a
strong sense that we call Morse general position. In
section~\ref{sec:goodregions}, this proof will be adapted to the weaker
concept of general position developed in section~\ref{sec:generalposition}.

Consider a smooth function $f\colon (\R^2,0)\to (\R,0)$. A critical point
of $f$ is \emph{stable} when it is locally equivalent under smooth change
of coordinates of the domain and range to $f(x,y)=x^2+y^2$ or
$f(x,y)=x^2-y^2$. The first type is called a \emph{center,} and the second
a \emph{saddle.} An unstable critical point is called a \emph{birth-death}
point if it is locally $f(x,y)=x^2+y^3$.

Let $\tau\colon P\times [0,1]\to M$ be sweepouts as in
section~\ref{sec:RSgraphic}. As in that section, we denote
$\tau(P\times\set{0,1})$ by $T$, $\tau(P\times\set{t})$ by $P_t$,
$\tau(P\times\set[0,t])$ by $V_t$, and $\tau(P\times[t,1])$ by $W_t$.  For
a second sweepout $\sigma\colon Q\times [0,1]\to M$, we denote
$\sigma(Q\times\set{0,1})$ by $S$, $\sigma(Q\times\set{s})$ by $Q_s$,
$\sigma(Q\times[0,s])$ by $X_s$, and $\sigma(Q\times[s,1])$ by $Y_s$.  We
call $Q_s$ a \emph{$\sigma$-level} and $P_t$ a~\emph{$\tau$-level.}

A tangency of $Q_s$ and $P_t$ at a point $w$ is said to be \emph{of Morse
type} at $w$ if in some local $xyz$-coordinates with origin at $w$, $P_t$
is the $xy$-plane and $Q_s$ is the graph of a function which has a stable
critical point or a birth-death point at the origin.  Note that this
condition is symmetric in $Q_s$ and $P_t$. We may refer to a tangency as
stable or unstable, and as a center, saddle, or birth-death point.

A tangency of $S$ with a $\tau$-level is said to be \emph{stable} if there
are local $xyz$-coordinates in which the $\tau$-levels are the planes
$\R^2\times\{z\}$ and $S$ is the graph of $z=x^2$ in the $xz$-plane. In
particular, the tangency is isolated and cannot occur at a vertex of
$S$. There is an analogous definition of stable tangency of $T$ with a
$\sigma$-level.

We will say that $\sigma$ and $\tau$ are in \emph{Morse general position}
when the following hold:
%\newcounter{MGPcounter}
%\begin{list}{(MP\arabic{MGPcounter})}
%{\usecounter{MGPcounter}
%\setlength{\labelwidth}{7 ex}
%\setlength{\leftmargin}{9 ex}
%\setlength{\labelsep}{2 ex}
%\setlength{\parsep}{5pt}
%}
\begin{enumerate}
\item
$S$ is disjoint from $T$, 
\item 
all tangencies of $S$ with $\tau$-levels and of $T$ with $\sigma$-levels
are stable,
\item 
all tangencies of $\sigma$-levels with $\tau$-levels are of Morse type, and
only finitely many are birth-death points,
\item
each pair consisting of a $\sigma$-level and a $\tau$-level has at most two
tangencies, and
\item
there are only finitely many pairs consisting of
a $\sigma$-level and a $\tau$-level with two tangencies, and for each of
these pairs both tangencies are stable.
\end{enumerate}

%\end{list}

The following concept due to A. Casson and C. McA.\ Gordon \cite{CG} is a
crucial ingredient in \cite{RS}. A Heegaard splitting $M=V\cup_P W$ is
called \emph{strongly irreducible} when every compressing disk for $V$
meets every compressing disk for $W$.

Suppose that $P$ is a Heegaard surface in $M$, bounding a handlebody
$V$. We define a \emph{precompression} or \emph{precompressing disk} for
$P$ in $V$ to be an imbedded disk $D$ in $M$ such that
\begin{enumerate}
\item $\partial D$ is an essential loop in $P$, 
\item $D$ meets $P$ transversely at $\partial D$,
and $V$ contains a neighborhood of $\partial D$,
\item the interior of $D$ is transverse to $P$, and its intersections with
$P$ are discal.
\end{enumerate}
Provided that $M$ is irreducible, a precompression for $P$ in $V$ is
isotopic relative to a neighborhood of $\partial D$ to a compressing disk
for $P$ in $V$. In particular, if the Heegaard splitting is strongly
irreducible, then the boundaries of a precompression for $P$ in $V$ and a
precompression for $P$ in $\overline{M-V}$ must intersect.

A sweepout is called \emph{strongly irreducible} when the associated
Heegaard splittings are strongly irreducible. We can now state the main
technical result of \cite{RS}.
\begin{theorem}[Rubinstein-Scharlemann]
Let $M\neq S^3$ be a closed orientable $3$-manifold, and let
$\sigma,\tau\colon F\times[0,1]\to M$ be strongly irreducible sweepouts of
$M$ which are in Morse general position. Then there exists $(s,t)\in
(0,1)\times (0,1)$ such that $Q_s$ and $P_t$ meet in good position.
\label{thm:RS}
\end{theorem}

The closure in $I^2$ of the set $(s,t)$ for which $Q_s$ and $P_t$ have a
tangency is a graph $\Gamma$. On $\partial I^2$, it can have valence-$1$
vertices corresponding to valence-$3$ vertices of $S$ or $T$, and
valence-$2$ vertices corresponding to points of tangency of $S$ with a
$\tau$-level or $T$ with a $\sigma$-level (see p.~1008 of \cite{RS}, see
also \cite{KS} for an exposition with examples). In the interior of $I^2$,
it can have valence-$4$ vertices which correspond to a pair of levels which
have two stable tangencies, and valence-$2$ vertices which correspond to
pairs of levels having a birth-death tangency.

The components of the complement of $\Gamma$ in the interior of $I^2$ are
called \emph{regions.}  Each region is either unlabeled or bears a label
consisting of up to four letters. The labels are determined by the
following conditions on $Q_s$ and $P_t$, which by transversality hold
either for every $(s,t)$ or for no $(s,t)$ in a region.
\begin{enumerate}
\item
If $Q_s$ contains a precompression for $P_t$ in $V_t$ (respectively, in
$W_t$), the region receives the letter $A$ (respectively, $B$).
\label{item:AB}
\item
If $P_t$ contains a precompression for $Q_s$ in $X_s$ (respectively, in
$Y_s$), the region receives the letter $X$ (respectively, $Y$).
\label{Item:XY}
\item
If the region has neither an $A$-label nor a $B$-label, and $V_t$
(respectively, $W_t$), contains a spine of $Q_s$, the region receives the
letter $b$ (respectively, $a$).
\item
If the region has neither an $X$-label nor a $Y$-label, and $X_s$
(respectively, $Y_s$), contains a spine of $P_t$, the region receives the
letter $y$ (respectively, $x$).
\end{enumerate}

With these conventions, $Q_s$ and $P_t$ are in good position if and only if
the region containing $(s,t)$ is unlabeled.  To check this, assume first
that they are in good position.  Since all intersections are biessential or
discal, neither surface can contain a precompressing disk for the other,
and since there is a biessential intersection circle, the complement of one
surface cannot contain a spine for the other. For the converse, an
intersection circle which is not biessential or discal leads to a
precompression as in~(1) or~(2), so assume that all intersections are
discal.  Then the complement of the intersection circles in $Q_s$ contains
a spine, so the region has either an $a$- or $b$-label, and by the same
reasoning applied to $P_t$ the region has either an $x$- or $y$-label. This
verifies the assertion, as well as the following lemma.

\begin{lemma}
If the label of a region contains the letter $a$ or $b$, then it must also
contain either $x$ or $y$. Similarly, if it contains $x$ or $y$, then it
must also contain $a$ or $b$.
\label{lem:no lone small letters}
\end{lemma}

We call the data consisting of the graph $\Gamma\subset I^2$ and the
labeling of a subset of its regions the \emph{Rubinstein-Scharlemann
graphic} associated to the sweepouts. Regions of the graphic are called
\emph{adjacent} if there is an edge of $\Gamma$ which is contained in both
of their closures.

At this point, we begin to make use of the fact that the sweepouts are
strongly irreducible. The labels will then have the following properties,
where $\sA$ stands for either of $A$ and $a$, and $\sB$, $\sX$, and $\sY$
are defined similarly.  
\newcounter{RScounter}
\begin{list}{(RS\arabic{RScounter})}
{\usecounter{RScounter}
\setlength{\labelwidth}{7 ex}
\setlength{\leftmargin}{10 ex}
\setlength{\labelsep}{2 ex}
\setlength{\parsep}{5pt}
}
\item
A label cannot contain both an $\sA$ and a $\sB$, or both an
$\sX$ and a $\sY$ (direct from the labeling rules and the
definition of strong irreducibility).
\label{item:noAB}
\item
If the label of a region contains $\sA$, then the label of any
adjacent region cannot contain $\sB$. Similarly for $\sX$ and
$\sY$ (Corollary~5.5 of~\cite{RS}).
\label{item:edge labeling}
\item
If all four letters $\sA$, $\sB$, $\sX$, and $\sY$ appear in the labels of
the regions that meet at a valence-$4$ vertex of $\Gamma$, then two
opposite regions must be unlabeled (Lemma~5.7 of \cite{RS}).
\label{item:2-cell labeling}
\end{list}

Property (RS\ref{item:edge labeling}) warrants special comment, since it
will play a major role in our later work. The analysis of labels of
adjacent regions given in section~5 of \cite{RS} uses only the fact that
for the points $(s,t)$ in an open edge of $\Gamma$, the corresponding $Q_s$
and $P_t$ have a single stable tangency. The open edges of the more general
graphics we will use for the diffeomorphisms in parameterized families in
general position will still have this property, so the labels of their
graphics will still satisfy property~(RS\ref{item:edge labeling}). They
will not satisfy property property~(RS\ref{item:2-cell labeling}), indeed
the $\Gamma$ for their graphics will have vertices of high valence, so
property~(RS\ref{item:2-cell labeling}) will not even be meaningful.

We now analyze the labels of regions whose closures meet $\partial I^2$, as
on p.~1012 of \cite{RS}. Consider first a region whose closure meets the
side $s=0$ (we consider $s$ to be the horizontal coordinate, so this is the
left-hand side of the square).  The region must contains points $(s,t)$
with $s$ arbitrarily close to~$0$. These correspond to $Q_s$ which are
extremely close to $S_0$. For almost all $t$, $S_0$ is transverse to $P_t$,
and for sufficiently small $s$ any intersection of such a $P_t$ with $Q_s$
must be an essential circle of $Q_s$ bounding a disk in $P_t$ that lies in
$X_s$, in which case the region must have an $X$-label. If $P_t$ is
disjoint from $Q_s$, then $P_t$ lies in $Y_s$ so the region has an
$x$-label. That is, all such regions have an $\sX$-label. Similarly, the
label of any region whose closure meets the edge $t=0$ (respectively,
$s=1$, $t=1$) contains $\sA$ (respectively, $\sY$, $\sB$).

\begin{figure}
\includegraphics[width=20 ex]{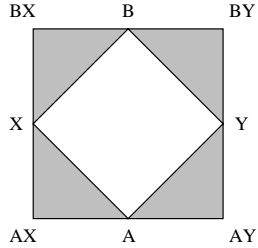}
\caption{The Diagram.}
\label{fig:Diagram}
\end{figure}

We will set up some of the remaining steps a bit differently from those of
\cite{RS}, so that their adaptation to our later arguments will be more
transparent. We have seen that it is sufficient to prove that there exists
an unlabeled region in the graphic defined by the sweepouts. To accomplish
this, Rubinstein and Scharlemann use the shaded subset of the square shown
in figure~\ref{fig:Diagram}. It is a simplicial complex in which each of
the four triangles is a $2$-simplex. Henceforth we will refer to it as
\emph{the Diagram.}

Suppose for contradiction that every region in the
Rubin\-stein-Schar\-le\-mann graphic is labeled.  Let $\Delta$ be a
triangulation of $I^2$ such that each vertex of $\Gamma$ and each corner of
$I^2$ is a $0$-simplex, and each edge of $\Gamma$ is a union of
$1$-simplices. Let $K$ be $I^2$ with the structure of a regular $2$-complex
dual to $\Delta$. We observe the following properties of $K$:
\newcounter{Kcounter}
\begin{list}{(K\arabic{Kcounter})}
{\usecounter{Kcounter}
\setlength{\labelwidth}{7 ex}
\setlength{\leftmargin}{10 ex}
\setlength{\labelsep}{2 ex}
\setlength{\parsep}{5pt}
}
\item
Each $0$-cell of $K$ lies in the interior of a side of
$\partial I^2$ or in a region.\par
\label{item:K1}
\item
Each $1$-cell of $K$ either lies in $\partial I^2$, or is disjoint from
$\Gamma$, or crosses one edge of $\Gamma$ transversely in one point.
\label{item:K2}
\item
Each $2$-cell of $K$ either contains no vertex of $\Gamma$, in which case
all of its $0$-cell faces that are not in $\partial I^2$ lie in one region
or in two adjacent regions, or contains one vertex of $\Gamma$, in which
case all of its $0$-cell faces which do not lie in $\partial I^2$ lie in
the union of the regions whose closures contain that vertex.
\label{item:K3}
\end{list}

We now construct a map from $K$ to the Diagram. First, each $0$-cell in
$\partial K$ is sent to one of the single-letter $0$-simplices of the
diagram: if it lies in the side $s=0$ (respectively, $t=0$, $s=1$, $t=1$)
then it is sent to the $0$-simplex labeled $\sX$ (respectively, $sA$, $sY$,
$sB$).  Similarly, any $1$-cell in a side of $\partial K$ is sent to the
$0$-simplex that is the image of its endpoints, and the four $1$-cells in
$\partial K$ dual to the original corners are send to the $1$-simplex whose
endpoints are the images of the endpoints of the $1$-cell. Notice that
$\partial K$ maps essentially onto the circle consisting of the four
diagonal $1$-simplices of the Diagram.

We will now show that if there is no unlabeled region, this map extends to
$K$, a contradiction. Since an unlabled region produces pairs $Q_s$ and
$P_t$ that meet in good position, this will complete the proof sketch of
theorem~\ref{thm:RS}.

Now we consider cells of $K$ that do not lie entirely in $\partial K$. Each
$0$-cell in the interior of $K$ lies in a region.  By (RS\ref{item:noAB}),
the label of each $0$-cell has a form associated to one of the
$0$-simplices of the Diagram, and we send the $0$-cell to that $0$-simplex.

Consider a $1$-cell of $K$ that does not lie in $\partial K$. Suppose it
has one endpoint in $\partial K$, say in the side $s=0$ (the other cases
are similar). The other endpoint lies in a region whose closure meets the
side $s=0$, so its label contains $\sX$. Therefore the images of the
endpoints of the $1$-cell both contain $\sX$, so lie either in a
$0$-simplex or a $1$-simplex of the Diagram. We extend the map to the
$1$-cell by sending it into that $0$- or $1$-simplex.  Suppose the $1$-cell
lies in the interior of $K$.  Its endpoints lie either in one region or in
two adjacent regions. If the former, or the latter and the labels of the
regions are equal, we send the $1$-cell to the $0$-simplex for that
label. If the latter and the labels of the regions are different, then
property (RS\ref{item:edge labeling}) shows that the labels span a unique
$1$-simplex of the Diagram, in which case we send the $1$-cell to that
$1$-simplex.

Assuming that the map has been extended to the $1$-cells in this way,
consider a $2$-cell of $K$. Suppose first that it has a face that lies in
the side $s=0$ (the other cases are similar). Then each of its $0$-cell
faces lies in one of the sides $s=0$, $t=0$, or $t=1$, or in a region whose
closure meets $s=0$. In the latter case, we have seen that the label of the
region must contain $\sX$, so it cannot contain $\sY$, and in particular it
cannot be a single letter $\sY$. In no case does the $0$-cell map to the
vertex $\sY$ of the Diagram, so the image of the boundary of the $2$-cell
maps into the complement of that vertex in the Diagram. Since that
complement is contractible, the map extends over the $2$-cell.

Suppose now that the $2$-cell lies entirely in the interior of $K$.  If it
is dual to a $0$-simplex of $\Delta$ that lies in a region or in the
interior of an edge of $\Gamma$, then all its $0$-cell faces lie in a
region or in two adjacent regions. In this case, all of its $1$-dimensional
faces map into some $1$-simplex of the Diagram, so the map on the faces
extends to a map of the $2$-cell into that $1$-simplex. Suppose the
$2$-cell is dual to a vertex of $\Gamma$. Its faces lie in the union of
regions whose closures contain the vertex. If the vertex has valence $2$,
then all $0$-cell faces lie in two adjacent regions (actually, in this
case, the regions must have the same label) and the map extends to the
$2$-cell as before. If the vertex has valence $4$, then by
(RS\ref{item:2-cell labeling}), the labels of the four regions whose
closures contain the vertex must all avoid at least one of the four
letters. This implies that the boundary of the $2$-cell of $K$ maps into a
contractible subset of the Diagram. So again the map can be extended over
the $2$-cell, giving us the desired contradiction.

We emphasize that the map from $K$ to the Diagram carries each $1$-cell of
$K$ to a $0$-simplex or a $1$-simplex of the Diagram, principally due to
property~(RS\ref{item:edge labeling}).

\newpage
\section[Graphics having no unlabeled region]
{Graphics having no unlabeled region}
\label{sec:examples}

One cannot hope to perturb a parameterized family of sweepouts to be in
Morse general position. One must allow for the possibility of levels having
tangencies of high order, and having more than two tangencies. We will see
in section~\ref{sec:generalposition} that all such phenomena can be
isolated at the vertices of the graph $\Gamma$ in the graphic. In
particular, the $(s,t)$ that lie on the open edges of $\Gamma$ will still
correspond to pairs of levels that have a single stable tangency, and
therefore their associated graphics will still have property
(RS\ref{item:edge labeling}). Achieving this property for the edges of
$\Gamma$ will require considerable effort, so before beginning the task, we
will show that the hard work really is necessary. We will give here
examples of pairs of sweepouts which have a graphic with no unlabeled
region. It will be clear that what goes wrong is the existence of edges of
$\Gamma$ that consist of pairs having multiple tangencies, and the
corresponding failure of the graphic to have property~(RS\ref{item:edge
labeling}).

This section is not part of the proof of the Smale Conjecture for Lens
Spaces, and can be read independently from the rest of the paper (provided
that one is familiar with Rubinstein-Scharlemann graphics and their
labeling scheme).

We will first construct examples in $S^2\times S^1$, then show how to
further modify them to obtain examples in any $L(m,q)$.

\begin{figure}
\includegraphics[width=0.95\textwidth]{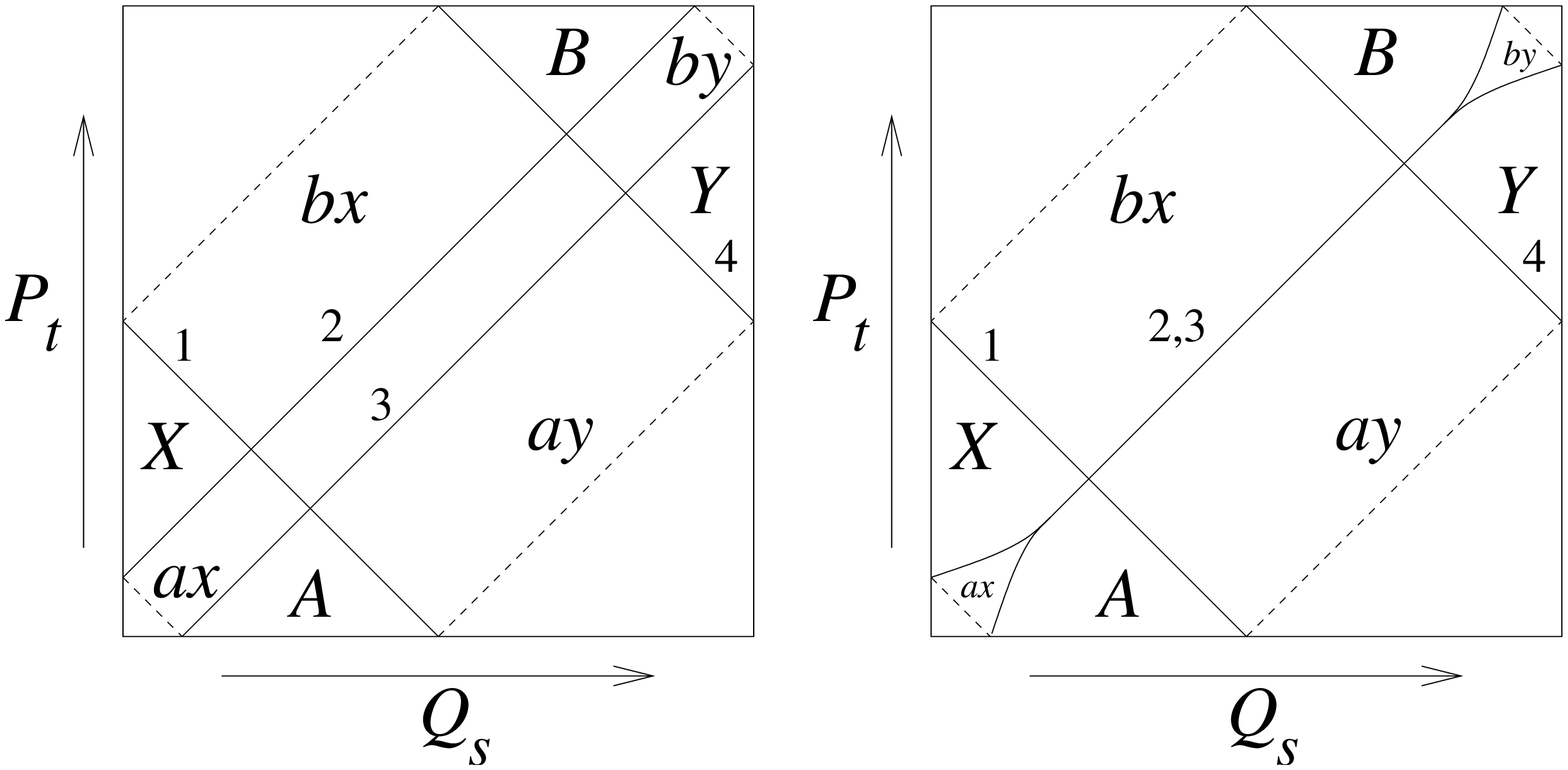}
\caption{Graphics before and after deformation.}
\label{fig:graphics}
\end{figure}

The first step is to construct a pair of sweepouts of $S^2\times S^1$, with
the graphic shown on the left in figure~\ref{fig:graphics}.  In
figure~\ref{fig:graphics}, the edges of pairs for which the corresponding
levels have a single center tangency are shown as dotted. The four corner
regions are not labeled, since their labels are the same as the regions
that are adjacent to them along an edge of centers.

After constructing the sweepouts that produce the first graphic, we will
see how to move one of the sweepouts by isotopy to ``collapse'' the
unlabeled region. Two edges of the first graphic are moved to coincide,
producing the graphic on the right in figure~\ref{fig:graphics}. The three
open edges that lie on the diagonal $y=x$ consist of pairs of levels which
have two saddle tangencies. The two vertices where the edges labeled $1$
and $4$ cross the diagonal correspond to pairs having three saddle
tangencies.

\begin{figure}[ht]
\includegraphics[width=66 ex]{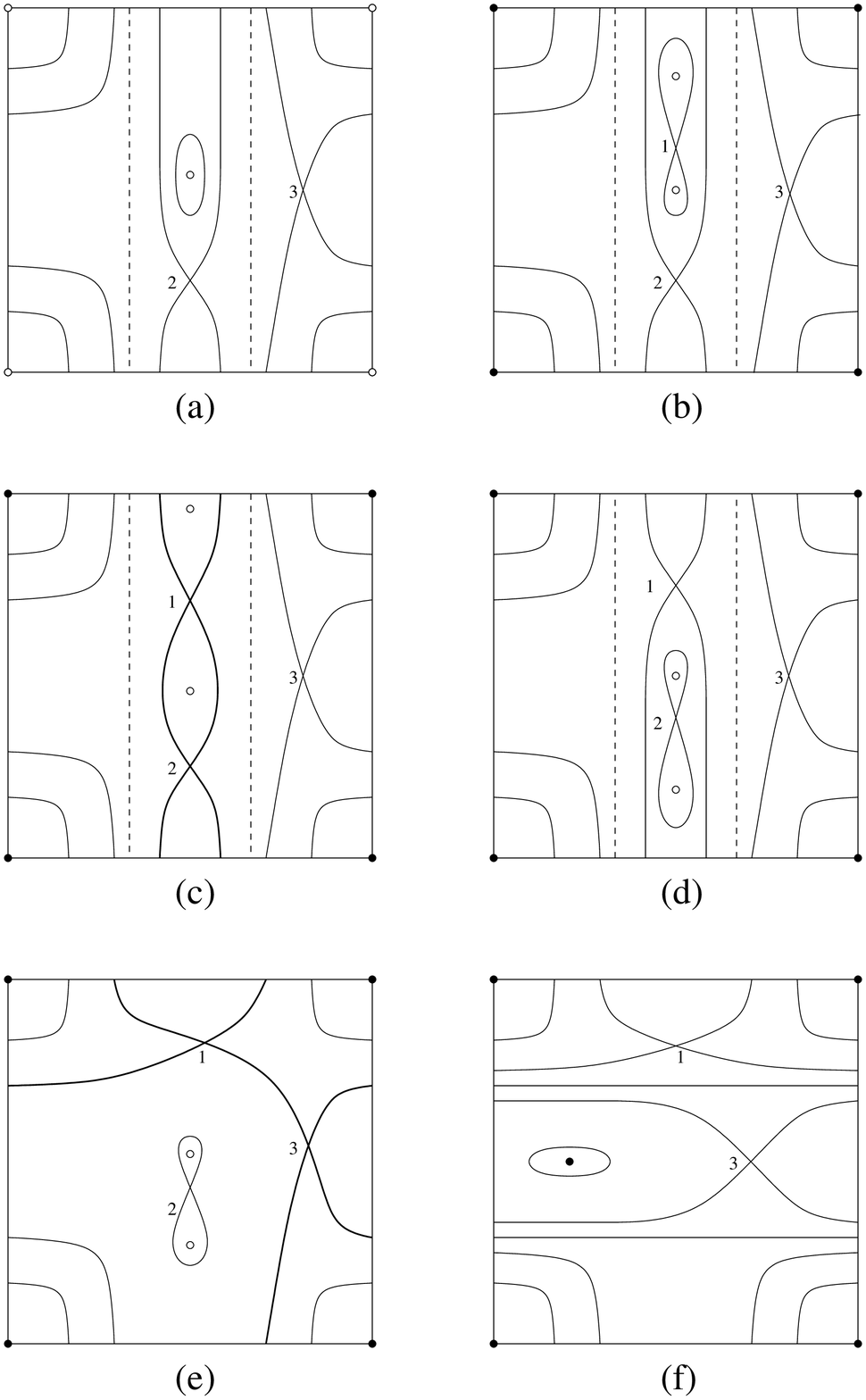}
\caption[sweepout2]{Intersections of the $Q_s$ with fixed $P_t$ as $t$
decreases from $1/2$ to $0$, for the sweepouts with an unlabeled region.
\begin{enumerate}
\item[(a)] $P_{1/2}$.
\item[(b)] $P_t$ where $s(1)<s(2)<s(3)$.
\item[(c)] $P_t$ where $s(1)=s(2)$.
\item[(d)] $P_t$ where $s(2)<s(1)<s(3)$.
\item[(e)] $P_t$ where $s(1)=s(3)$.
\item[(f)] $P_t$ where $s(3)<s(1)$, and after saddle 2 changes to a center.
\end{enumerate}}
\label{fig:Morse}
\end{figure}

\begin{figure}[ht]
\includegraphics[width=66 ex]{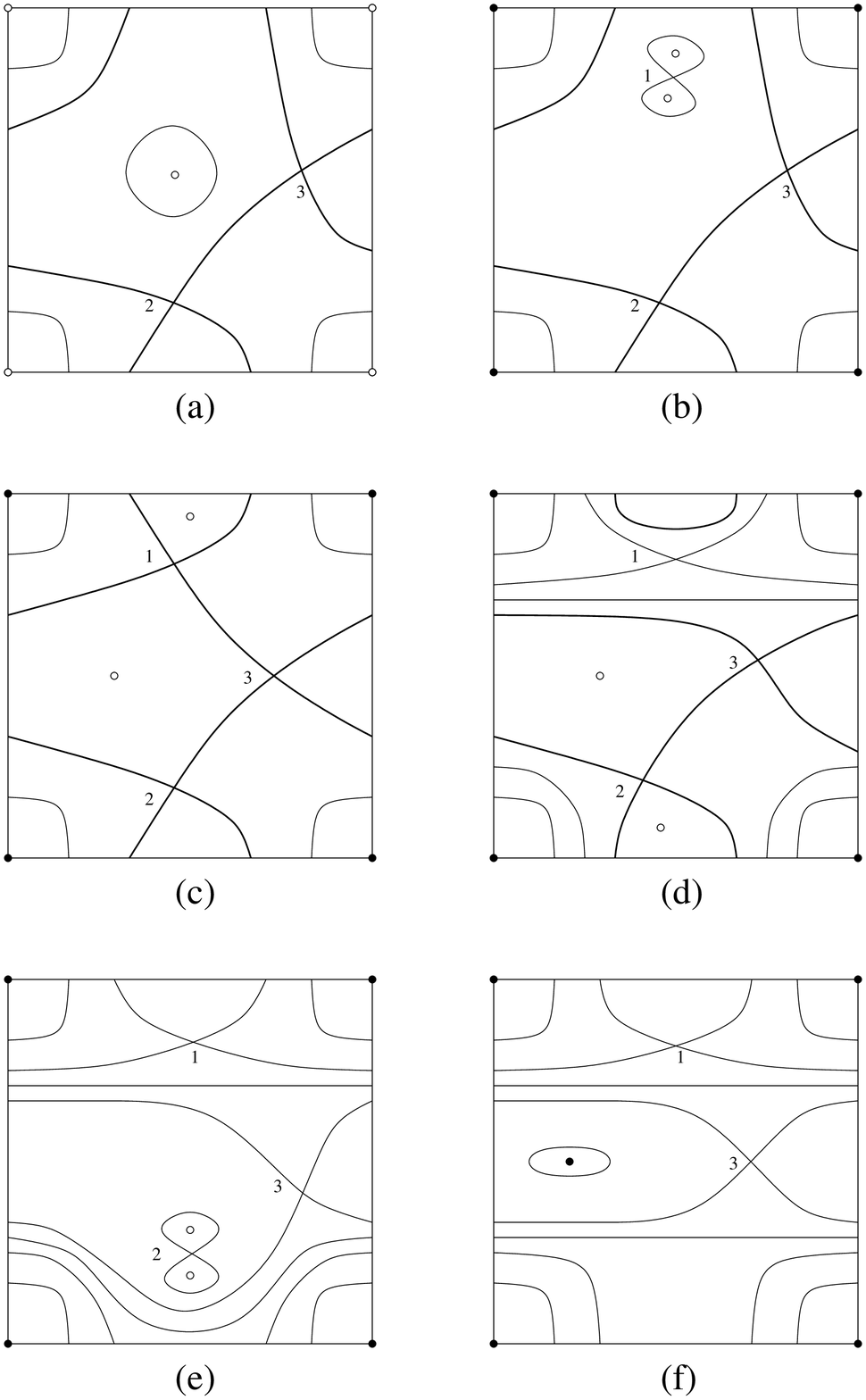}
\caption[sweepout2]{Intersections of the $Q_s$ with fixed $P_t$ as $t$
decreases from $1/2$ to $0$, for the sweepouts with no unlabeled region.
\begin{enumerate}
\item[(a)] $P_{1/2}$.
\item[(b)] $P_t$ where $s(1)<s(2)=s(3)$.
\item[(c)] $P_t$ where $s(1)=s(2)=s(3)$.
\item[(d)] $P_t$ where $s(2)=s(3)<s(1)$.
\item[(e)] $P_t$ where $s(2)<s(3)<s(1)$.
\item[(f)] $P_t$ where $s(3)<s(1)$, and after saddle 2 changes to a center.
\end{enumerate}}
\label{fig:nogood}
\end{figure}

\begin{figure}[ht]
\includegraphics[width=66 ex]{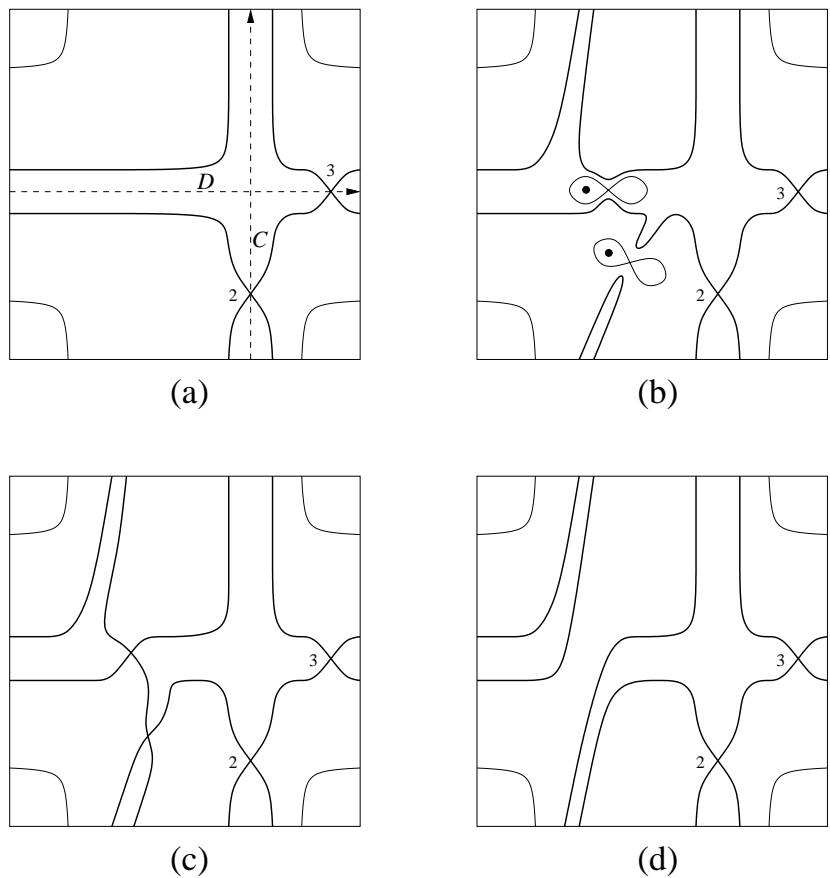}
\caption[sweepout2]{Various $P_t$-levels illustrating a bowtie
elaboration.}
\begin{enumerate}
\item[(a)] $P_t$ for $t$ below the bowtie elaboration.
\item[(b)] $P_t$ above the two birth-death points, but below $t_0$.
\item[(c)] $P_{t_0}$.
\item[(d)] $P_t$ for $t>t_0$.
\end{enumerate}
\label{fig:bowtie_levels}
\end{figure}

As it is rather difficult to visualize the sweepouts directly, we describe
them by level pictures for various $P_t$.  The $Q_s$ appear as level curves
in each $P_t$. Here are some general conventions:
\begin{enumerate}
\item[(i)] A solid dot is a center tangency.
\item[(ii)] An open dot (i.~e.~a tiny circle) is a point in one of the singular
circles $S_i$ of the $Q_s$-sweepout.
\item[(iii)] Double-thickness lines are intersections with a $Q_s$ that
have more than one tangency.
\item[(iv)] In figures~\ref{fig:Morse} and~\ref{fig:nogood}, dashed lines
are biessential intersection circles (in figure~\ref{fig:bowtie_levels},
they have a different meaning).
\end{enumerate}

In a picture of a $P_t$, the level curves $P_t\cap Q_s$ that contain
saddles appear as curves with self-crossings, and we label the crossings
with $1$, $2$, $3$, or $4$ to indicate which edge of the graphic in
figure~\ref{fig:graphics} contains that $(s,t)$-pair. For a fixed $t$,
$s(n)$ will denote the $s$-level of saddle~$n$. That is, in the graphic the
edge of $\Gamma$ labeled $n$ contains the point $(s(n),t)$.

Figure~\ref{fig:Morse} shows some $P_t$ with $t\leq 1/2$, for a sweepout of
$S^2\times S^1$ whose graphic is the one shown in the left of
figure~\ref{fig:graphics}. Here are some notes on figure~\ref{fig:Morse}.
\begin{enumerate}
\item In (a)-(f), the circles $x=\text{constant}$ are longitudes of $V_t$,
and the circles $y=\text{constant}$ are meridians.
\item The point represented by the four corners is the point of $P_t$ with
largest $s$-level. In (a) it is a tangency of $P_{1/2}$ with $S_1$, and in
(b)-(f) it is a center tangency of $P_t$ with $Q_{t+1/2}$.
\item The open dots in the interior of the squares are intersections of
$P_t$ with $S_0$. In (a) it is a tangency of $P_{1/2}$ with $S_0$, in
(b)-(e) they are transverse intersections. In (f), $P_t$ is disjoint from
$S_0$.
\item In (b), saddle~1 has appeared. Circles of $Q_s\cap P_t$ with $s<s(1)$
are essential in $Q_s$, and these $(s,t)$ lie in the region labeled $X$ in
the graphic. Circles of $Q_s\cap P_t$ with $s(1)<s<s(2)$ enclose the
figure-$8$ in (b), which is $P_t\cap Q_{s(1)}$. They are inessential in
both $Q_s$ and $P_t$, and these $(s,t)$ lie in the region labeled $bx$. The
vertical dotted lines are biessential intersections corresponding to a pair
in the unlabeled region. Finally, one crosses $Q_{s(3)}$, and eventually
reaches the center tangency.
\item The horizontal level curves shown in (f) are meridians of $V_t$ that
bound disks in the $Q_s$ that contain them. This $(s,t)$ lies in the region
labeled $A$ in the graphic.
\end{enumerate}

For $t>1/2$, the intersection pattern of $P_t$ with the $Q_s$ is isomorphic
to the pattern for $P_{1-t}$, by an isomorphism for which $Q_s$ corresponds
to $Q_{1-s}$. As one starts $t$ at $1/2$ and moves upward through
$t$-levels, saddle $4$ appears inside the component of $P_t-Q_{s(3)}$ that
is an open disk, and expands until the level where $s(3)=s(4)$. The
biessential intersection circles in (a)-(d) are again longitudes in $V_t$
and in $W_t$, and the horizontal intersection circles in (f) are meridians
of $W_t$. These $(s,t)$ lie in the region labeled $B$ in the graphic. This
completes the description of the sweepouts in Morse general position.

Figure~\ref{fig:nogood} shows some $P_t$ for a sweepout of $S^2\times S^1$
whose graphic is the one shown in the right of figure~\ref{fig:graphics}.
This sweepout is obtained from the previous one by an isotopy that moves
parts of the $Q_s$ levels down (to lower $t$-levels) near saddle $2$ and
up near saddle $3$. Again, the portion that is shown fits together with a
similar portion for $1/2\leq t\leq 1$. As $t$ increases past $1/2$, saddle
$4$ appears in the component of $P_t-S_{s(2)}$ that contains the point
which appears as the four corners.

\newpage
\newpage

\begin{figure}
\includegraphics[width=24 ex]{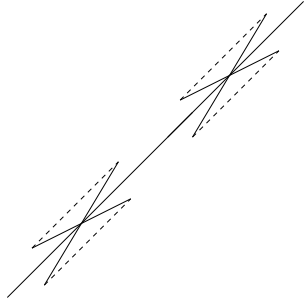}
\caption[sweepout2]{Bowtie elaborations of the $2,3$-edge}
\label{fig:bowtie}
\end{figure}

\newpage
We will now explain how to modify this construction to obtain a pair of
sweepouts with no unlabeled region for any $L(m,q)$. The graphic will be
the same as the one on the right in figure~\ref{fig:graphics}, except that
near the point $(1/2,1/2)$, a small portion of the $2,3$-edge will have a
sequence of elaborations, called \emph{bowties,} two of which are shown in
figure~\ref{fig:bowtie}. We remark that for any $t$ near $1/2$, one has
$s(2)=s(3)=t$, since the $2,3$-edge is the diagonal of the graphic. That
is, $P_t\cap Q_t$ contains saddles~$2$ and~$3$ for $t$ near~$1/2$.

Figure~\ref{fig:bowtie_levels} shows various $P_t$ for a bowtie
elaboration. Consider the lower-left bowtie in figure~\ref{fig:bowtie}. Let
$(t_0,t_0)$ be the point where its two saddle edges cross the diagonal
$2,3$-edge. Figure~\ref{fig:bowtie_levels}(a) shows $P_t\cap Q_t$ in $P_t$
for $t$ near $1/2$ but with $t$ below the level where the bowtie
elaboration begins.  This could be the lower endpoint of the portion of the
$2,3$-edge shown in figure~\ref{fig:bowtie}.  We have drawn the levels in
$P_t$ a bit differently from the picture of $P_{1/2}$ in
figure~\ref{fig:nogood}, but this picture is isotopic to
figure~\ref{fig:nogood}(a).  In figure~\ref{fig:bowtie_levels}(a), the
closure $X$ of one of the components of $P_t-(P_t\cap Q_t)$ deformation
retracts to a figure-$8$ $C\cup D$, where the circle $C$ passes through
saddle $2$ and is a longitude $L$ of $V_t$, and the circle $D$
passes through saddle $3$ and is a meridian $M$ of~$V_t$.

As one moves up in $t$-levels, one passes through two birth-death points,
producing two center-saddle pairs. One birth-death point is at an $s$-level
with $s<t$, and the other is at an $s$-level with $s>t$, so they lie in
different components of $P_t-Q_t$.  Figure~\ref{fig:bowtie_levels}(b) shows
the two new center and saddle tangencies, along with $P_t\cap Q_t$, in such
a level with $t<t_0$. Figure~\ref{fig:bowtie_levels}(c) shows
$P_{t_0}$. The two new saddles are then in $Q_{t_0}$, along with saddles
$2$ and $3$. In figure~\ref{fig:bowtie_levels}(d), we see $P_t\cap Q_t$ in
a $P_t$ for $t>t_0$. The effect is to reposition the component $X$ so that
$C$ still represents $L$, but $D$ represents $M+L$.

We remark that figure~\ref{fig:bowtie_levels} is rather schematic. Since
the bowtie elaborations lie very close to the diagonal, the new centers and
saddles actually lie very close to $Q_t$. The elongations on $P_t\cap Q_t$
that reach toward the saddle points in figure~\ref{fig:bowtie_levels}(c)
would actually follow along very close to the $P_t\cap Q_t$ of
figure~\ref{fig:bowtie_levels}(a).

By a sequence of such bowtie elaborations, one can change $P_t\cap Q_t$ so
that $C$ and $D$ represent any pair of generators of $H_1(P_t)$. In
particular, we may move them so that $D$ represents $mL+qM$, the meridian
of $W_t$ in $L(m,q)$. Then, $C$ is a longitude of $W_t$, and the portion in
$W_{1/2}$ of the sweepouts of $S^2\times S^1$ with no unlabeled region can
be placed into this $W_t$ in $L(m,q)$, producing a pair of sweepouts whose
graphic is obtained from the one on the right in figure~\ref{fig:graphics}
by bowtie elaborations.  Since the bowtie elaborations produce no unlabeled
regions, this graphic has no unlabeled region.

\newpage
\section[Graphics for parameterized families]
{Graphics for parameterized families}
\label{sec:generalposition}

In this section we prove that a parameterized family of sweepouts can be
perturbed so that a suitable graphic exists at each parameter. As discussed
in section~\ref{sec:examples}, in a parameterized family one must allow for
the possibility of levels having tangencies of high order, and having more
than two tangencies.

Additional complications arise because one cannot avoid having parameters
where the singular sets of the sweepouts intersect, or where the singular
sets have high-order tangencies with levels. We sidestep these
complications by working only with sweepout parameters that lie in an
interval $[\epsilon,1-\epsilon]$. The graphic is only considered to exist
on the square $[\epsilon,1-\epsilon]\times[\epsilon,1-\epsilon]$, which we
call $I^2_\epsilon$. The number $\epsilon$ is chosen so that the labels of
regions whose closure meets a side of $I^2_\epsilon$ will be known to
include certain letters. Just as before, this will ensure that the map to
the Diagram be essential on the boundary of the dual complex~$K$.

These considerations, and the examples in section~\ref{sec:examples},
motivate our definition of a general position family of diffeomorphisms.
As usual, let $M$ be a closed orientable $3$-manifold and $\tau\colon
P\times [0,1]\to M$ a sweepout with singular set $T=T_0\cup T_1$ and level
surfaces $P_t$ bounding handlebodies $V_t$ and $W_t$. Let $f\colon M\times
W\to M$ be a parameterized family of diffeomorphisms, where $W$ is a
compact manifold. For $u\in W$ we denote the restriction of $f$ to
$M\times\set{u}$ by $f_u$. When a choice of parameter $u$ has been fixed,
we denote $f_u(P_s)$ by $Q_s$, and $f_u(V_s)$ and $f_u(W_s)$ by $X_s$ and
$Y_s$ respectively. When $Q_s$ meets $P_t$ transversely, a label is
assigned to $(s,t)$ as in section~\ref{sec:Rubinstein-Scharlemann}.

A preliminary definition will be needed. We say that a positive number
$\epsilon$ \emph{gives border label control} for $f$ if the following hold
at each parameter $u$:
\begin{enumerate}
\item
If $t\leq 2\epsilon$, then there exists $r$ such that $Q_r$ meets $P_t$
transversely and contains a compressing disk of $V_t$.
\item
If $t\geq 1-2\epsilon$, then there exists $r$ such that $Q_r$ meets $P_t$
transversely and contains a compressing disk of $W_t$.
\item
If $s\leq 2\epsilon$, then there exists $r$ such that $P_r$ meets $Q_s$
transversely and contains a compressing disk of $X_s$.
\item
If $s\geq 1-2\epsilon$, then there exists $r$ such that $P_r$ meets $Q_s$
transversely and contains a compressing disk of $Y_s$.
\end{enumerate}

Throughout this section, a \textit{graph} is a compact space which is a
disjoint union of a CW-complex of dimension $\leq 1$ and circles. The
circles, if any, are considered to be open edges of the graph.

We say that $f$ is \textit{in general position} (with respect to the
sweepout $\tau$) if there exists $\epsilon >0$ such that $\epsilon$ gives
border label control for $f$ and such that the following hold for each
parameter $u\in W$.  
\newcounter{GPcounter}
\begin{list}{(GP\arabic{GPcounter})}
{\usecounter{GPcounter}
\setlength{\labelwidth}{7 ex}
\setlength{\leftmargin}{9 ex}
\setlength{\labelsep}{2 ex}
\setlength{\parsep}{5pt}
}
\item 
For each $(s,t)$ in $I^2_\epsilon$, $Q_s\cap P_t$ is a graph. At each
point in an open edge of this graph, $Q_s$ meets $P_t$ transversely. At
each vertex, they are tangent.
\item
The $(s,t)\in I^2_\epsilon$ for which $Q_s$ has a tangency with $P_t$ form
a graph $\Gamma_u$ in $I^2_\epsilon$.
\item
If $(s,t)$ lies in an open edge of $\Gamma_u$, then $Q_s$ and $P_t$ have a
single stable tangency.
\end{list}

Here is the main result of this section.
\begin{theorem}
Let $f\colon M\times W\to M$ be a parameterized family of diffeomorphisms.
Then by an arbitrarily small deformation, $f$ can be put into general
position with respect to $\tau$.
\label{thm:generalposition}
\end{theorem}
\noindent The proof of theorem~\ref{thm:generalposition} will constitute
the remainder of this section. Since the argument is rather long, we will
break it into subsections. Until subsection~\ref{sec:borderlabel}, $M$ can
be a closed manifold of arbitrary dimension~$m$.

\subsection{The Parameterized Extension Principle}
\label{subsec:PEP}

Let $f\colon M\times W\to M$ be a parameterized family of diffeomorphisms.
Here, $W$ and $M$ are smooth manifolds, with $M$ closed and $W$ compact,
and $f$ is continuous as a map from $W$ to $\Diff(M)$, where $\Diff(M)$ has
as always the $C^\infty$-topology.

We first recall that $\Diff(M)$ is locally convex in the following strong
sense. Fix a Riemannian metric on $M$ for which $\partial M$ is totally
geodesic, and let $f'\colon M\times W\to M$ be a parameterized family of
smooth maps. Assume that $f'$ is close enough to $f$, in the compact-open
topology on maps from $W$ to the space of smooth maps $C^\infty(M,M)$, so
that for each $u\in W$ and each $x\in M$, there is a unique short vector
$v_{x,u}$ at $f(x,u)$ such that $\Exp(v_{x,u})=f'(x,u)$. Putting
$F_t(x,u)=\Exp(tv_{x,u})$ defines a parameterized family $F_t$ of
homotopies from $f_u$ to $f_u'$. The diffeomorphisms form an open subset of
the smooth maps from $M$ to $M$, so when $f'$ is sufficiently close to $f$,
each $(F_t)_u$ and in particular $f'$ will be a parameterized family of
diffeomorphisms. Consequently, \textit{if a modification of a parameterized
family $f$ of diffeomorphisms can be achieved by taking a family which can
be selected to be arbitrarily close in the $C^\infty$-topology, then it can
be achieved by a deformation of $f$ through families of diffeomorphisms.}

This observation works just as well when $M$ is an open manifold, provided
that we use the spaces $\Diff_c(M)$ and $C^\infty_c(M,M)$ of
diffeomorphisms and maps with compact support (those which agree with the
identity outside of a compact subset of $M$).

By very similar considerations, if $f\colon M'\times W\to M$ is a
parameterized family of imbeddings of a submanifold $M'$ of $M$ (possibly
of codimension $0$) into $M$, then any map $f'\colon M'\times W\to M$
sufficiently close to $f$ will also consist of imbeddings, and is homotopic
to $f$ through parameterized families of imbeddings.

We now state a powerful extension theorem for isotopies of submanifolds,
due to R. Palais \cite{Palais}. In the theorem, $N$ is a not necessarily
compact manifold, all spaces of maps have the strong $C^r$-topology, $1\leq
r\leq \infty$, $\Imb(X,N)$ denotes the space of smooth imbeddings of the
submanifold $X$ (with $X$ a closed manifold)
of $N$ into $N$, and $\Diff_c(N)$ denotes the space
of diffeomorphisms of $N$ with compact support.
\begin{theorem}[Palais Extension Theorem] 
Let $N$ be a smooth manifold without boundary, and $X$ and $Y$
submanifolds of $N$, with $Y\subseteq X$. Then the restriction maps
$\Imb(X,N)\to \Imb(Y,N)$ and $\Diff_c(N)\to\Imb(Y,N)$ are locally trivial
fibrations.
\label{thm:Palais}
\end{theorem}
\noindent
In most of our applications of the Palais Extension Theorem, we will need
considerable control. This control is present in Palais' setup, but not
explicit in the statement of the Palais Extension Theorem, so we will
rephrase Palais' method to prove the precise statement that will be needed.
\begin{theorem}[Parameterized Extension Principle]
Let $M$ and $W$ be compact smooth manifolds, let $M_0$ be a
submanifold of $M$ of positive codimension, and let $U$ be an open subset
of the interior of $M$ with $M_0\subset U$. Suppose that $F\colon M\times
W\to M$ is a parameterized family of diffeomorphisms of $M$. If $g\in
C^\infty(M_0\times W,M)$ is sufficiently close to $F\vert_{M_0\times W}$,
then there is a deformation $G$ of $F$ such that $G\vert_{M_0\times W}=g$,
and $G=F$ on $(M-U)\times W$. By selecting $g$ sufficiently close to
$F\vert_{M_0\times W}$, $G$ may be selected arbitrarily close to~$F$.
\end{theorem}
\begin{proof}
The key step in the proof of the Palais Extension Theorem is a method of
extending imbeddings to diffeomorphisms, given as Lemmas~c and~d
in~\cite{P}. Choose any Riemannian metric on $M$. Fixing an imbedding
$i\colon M_0\to U$ and an imbedding $j$ sufficiently close to $i$, a
section of the tangent bundle of $M$ is defined over $i(M_0)$ by choosing
at each $i(x)$ the unique short vector $w_x$ such that the exponential
function at $i(x)$ sends $w_x$ to $j(x)$. Using a construction involving
parallel translation along paths in the fibers of a tubular neighborhood of
$i(M_0)$, the section over $i(M_0)$ is extended to a vector field $w$ on
$M$, with compact support in $U$. The map $J\colon M\to M$ that carries
each $p$ to $\Exp_p(w_p)$ sends each $i(x)$ to $j(x)$.  When $j$ is close
to $i$, the vector field $w$ is close to the zero vector field, so $J$ is
close to the identity in the $C^\infty$-topology.  Since in the
$C^\infty$-topology the diffeomorphisms form an open subset of the smooth
maps, $J$ will be a diffeomorphism when $j$ is sufficiently close to
$i$. If $g$ is sufficiently close to $F\vert_{M_0\times W}$ so that each
$g_u$ is an imbedding, then this process can be applied at each parameter
$u$ to the imbeddings $i=F\vert_{M_0\times\set{u}}$ and $j=g_u$.  The
tubular neighborhoods must be selected to vary continuously, so that the
resulting $J_u$ vary continuously in $u$. The family $G$ is defined
by~$G_u=J_u\circ F_u$.
\end{proof}

%\noindent The support of $G$ on $U$ 
%in the Parameterized Extension Principle is possible
%because the Palais Extension Theorem uses the space of diffeomorphisms with
%compact support, and can be applied with $N=U$.

%There are various places where we need to apply results about closed or
%open manifold to the case of manifolds with boundary. Our general rule is
%that a map on a manifold $M$ with boundary is considered to be smooth if it
%extends to an open neighborhood of $M$ in the manifold obtained from $M$
%and $\partial M\times [0,\infty)$ by identifying $\partial M$ with
%$\partial M\times \{ 0\}$. This condition allows one to extend the map to a
%map with domain the double of $M$ along its boundary. Perturbations of this
%map restrict to perturbations of the original map on $M$. Such procedures
%are sufficient for all of our needs, and we omit these details. It is not
%difficult also to develop a version of the Parameterized Extension
%Principle for the case of manifolds with boundary, compact submanifolds
%with boundary, and proper imbeddings, but we will not need this. [check???]

\subsection{Weak transversality}
Although individual maps may be put transverse to a submanifold of the
range, it is not possible to perturb a parameterized family so that each
individual member of the family is transverse. But a very nice result of
J. W. Bruce, Theorem~1.1 of~\cite{Bruce}, allows one to simultaneously
improve the members of a family.
\begin{theorem}[J. W. Bruce]
Let $A$, $B$ and $U$ be smooth manifolds and $C\subset B$ a submanifold.
There is a residual family of mappings $F\in C^\infty(A\times U,B)$ such
that:
\begin{enumerate}
\item[(a)] For each $u\in U$, the restriction $F_u=F|_{A\times \{u\}}\colon
A\to B$ is transverse to $C$ except possibly on a discrete set of points.
\item[(b)] For each $u\in U$, the set $F_u^{-1}(C)$ is a smooth submanifold
of codimension equal to the codimension of $C$ in $B$, except possibly at a
discrete set of points. At each of these exceptional points $F_u^{-1}(C)$
is locally diffeomorphic to the germ of an algebraic variety, with the
exceptional point corresponding to an isolated singular point of the
variety.
\end{enumerate}
\label{thm:Bruce1}
\end{theorem}
\noindent
That is, $F_u^{-1}(C)$ is smooth except at isolated points where it has
topologically a nice cone-like structure.  It is not assumed that any of
the manifolds involved is compact.

Theorem 1.3 of \cite{Bruce} is a version of theorem~\ref{thm:Bruce1} in
which $C$ is replaced by a bundle $\phi\colon B\to D$. The statement is:
\begin{theorem}[J. W. Bruce]
For a residual family of mappings $F\in C^\infty(A\times U,B)$, the
conclusions of theorem~\ref{thm:Bruce1} hold for all submanifolds
$C=\phi^{-1}(d)$, $d\in D$.
\label{thm:Bruce2}
\end{theorem}

We should comment on the significance of the residual subset in these two
theorems. The method of proof of these theorems is to define, in an
appropriate jet space, a locally algebraic subset which contains the jets
of all the maps that fail these weak transversality conditions. These
subsets have increasing codimension as higher-order jets are taken. A
variant of Thom transversality (lemma~1.6 of \cite{Bruce}) allows one to
perturb a parameterized family of maps so that these jets are avoided and
the conclusion holds.  When $A$ and $W$ are compact, the image of $A\times
W$ will lie in the open complement of the locally algebraic sets of
sufficiently high codimension. Consequently, any map sufficiently close to
the perturbed map will also satisfy the conclusions of the theorems. In all
of our applications, the spaces involved will be compact, and \textit{we
tacitly assume that the result of any procedure holds on an open
neighborhood of the perturbed map.}

We now adapt the methodology of Bruce to prove a version of
theorem~\ref{thm:Bruce1} in which the submanifold $C$ is replaced by the
zero set of a nontrivial polynomial. We will prove it only for the case
when $A=I$, although a more general version should be possible.

\begin{proposition}
Let $P\colon \R^n\to \R$ be a nonzero polynomial and put $V=P^{-1}(0)$.
Let $W$ be compact. Then for all $G$ in an open dense subset of
$C^\infty(I\times W,\R^n)$, each $G_u^{-1}(V)$ is finite.
\label{prop:polynomial_Bruce}
\end{proposition}

\begin{proof}
Let $J_0^k(1,n)$ be the space of germs of degree-$k$ polynomials from
$(\R,0)$ to $\R^n$; an element of $J_0^k(1,n)$ can be written as
$(a_{1,0}+a_{1,1}t+\cdots+a_{1,k}t^k,\ldots,a_{n,0}+a_{n,1}t+\cdots+a_{n,k}t^k)$,
so that $J_0^k(1,n)$ can be identified with $\R^{(k+1)n}$. Note that the
jet space $J^k(I,\R^n)$ can be regarded as $I\times J_0^k(1,n)$, by
identifying the jet of $\alpha\colon I\to \R^n$ at $t_0$ with the jet of
$\alpha(t-t_0)$ at $0$.

Define a polynomial map $P_*\colon J_0^k(1,n)\to J_0^k(1,1)$ by applying
$P$ to the $n$-tuple $(a_{1,0}+a_{1,1}t+\cdots+a_{1,k}t^k,
\ldots,a_{n,0}+a_{n,1}t+\cdots+a_{n,k}t^k)$, and then taking only the terms
up to degree~$k$. The preimage $P_*^{-1}(0)$ is the set of $\alpha$ for
which $P\circ \alpha\,(0)=0$ and the first $k$ derivatives of $P\circ \alpha$
at $0$ vanish, that is, the set of germs of paths that lie in $V$ up to
$k^{th}$-order.

\begin{lemma} If $P$ is nonconstant, then as a map from $\R^{(k+1)n}$ to 
$\R^{k+1}$, $P_*$ has maximal rank.
\label{lem:Jacobian}
\end{lemma}

\begin{proof}
We may select notation so that
$P(X,Y_1,\ldots,Y_m)=P_0(Y_1,\ldots,Y_m)+X^rP_r(Y_1,\ldots,Y_m)+\cdots+
X^sP_s(Y_1,\ldots,Y_m)$ with $P_r$ nonzero, and write elements of
$J_0^k(1,n)$ as $(a_0+a_1t+\cdots+a_kt^k,b_0+b_1t+\cdots)$.  The Jacobian
of $P_*$ is a $(k+1)\times ((k+1)n)$ matrix, and we will show that its
leftmost $(k+1)\times (k+1)$-block is lower triangular with nonzero (as
polynomials) diagonal entries.

Write $P_*(a_0+a_1t+\cdots+a_kt^k,b_0+b_1t+\cdots)$ as
$Q_0+Q_1t+\cdots+Q_kt^k$, where the $Q_i$ are polynomials on $\R^{(k+1)n}$.
Note that $Q_0=P(a_0,b_0,\ldots)$. So, the $(1,1)$-entry of the Jacobian is
$\dfrac{\partial Q_0}{\partial a_0}= \dfrac{\partial}{\partial a_0}
P(a_0,b_0,\ldots)$, while the $(1,i+1)$-entries, $1\leq i\leq k$, are
$\dfrac{\partial Q_0}{\partial a_i}=0$.

For $j\geq 1$, we have 
$Q_j=\dfrac{1}{j!}\dfrac{\partial^j P_*}{\partial t^j}\bigg\vert_{t=0}$, 
$\dfrac{\partial X}{\partial a_i}=t^i$, and
$\dfrac{\partial Y_\ell}{\partial a_i}=0$,
so $\dfrac{\partial Q_j}{\partial a_i}$ is
\[\dfrac{1}{j!}
\dfrac{\partial^j}{\partial t^j} \bigg( \dfrac{\partial P}{\partial
X}\dfrac{\partial X}{\partial a_i} +\sum_{\ell=1}^m\dfrac{\partial
P}{\partial Y_\ell}\dfrac{\partial Y_\ell}{\partial
a_i}\bigg)\bigg\vert_{t=0} =\dfrac{1}{j!}\dfrac{\partial^j}{\partial
t^j}\bigg(\dfrac{\partial P}{\partial X}t^i\bigg)\bigg\vert_{t=0}\] which
vanishes for $i>j$, and is
$\dfrac{\partial P}{\partial X}\bigg\vert_{t=0}=\dfrac{\partial}{\partial a_0}
P(a_0,b_0,\ldots)$ for $i=j$.
\end{proof}

For each $k$, put $Z_k=P_*^{-1}(0)$.  Lemma~\ref{lem:Jacobian} shows that
$Z_k$ is a variety of codimension $k+1$ in $J^k_0(1,n)$.  Observe that if
$\alpha\colon(\R,0)\to \R^n$ is a germ of a smooth map, and $0$ is a limit
point of $\alpha^{-1}(V)$, then all derivatives of $P\circ \alpha$ vanish
at $0$. That is, the $k$-jet of $\alpha$ at $t=0$ is contained in $Z_k$ for
every~$k$.

By Lemma~1.6 of \cite{Bruce}, there is a residual set of maps $G\in
C^\infty(I\times W,\R^n)$ such that the jet extensions $j^kG\colon I\times
W\to J^k(I,\R^n)$ defined by $j^kG(t,u)=j^kG_u(t)$ are transverse to $I\times
Z_k$.  For $k+1$ larger than the dimension of $I\times W$, this says that
no point of $G_u^{-1}(0)$ is a limit point, so each $G_u^{-1}(0)$ is
finite.
\end{proof}

\subsection{Finite singularity type}
For our later work, we will need some ideas from singularity theory.  Let
$f\colon (\R^m,0)\to (\R^p,0)$ be a germ of a smooth map. There is a
concept of \emph{finite singularity type} for $f$, whose definition is
readily available in the literature (for example, \cite[p.~117]{Bruce}).
The basic idea of the proof of theorem~\ref{thm:Bruce1} (given as
Theorem~1.1 in \cite{Bruce}) is to regard the submanifold $C$ locally as
the preimage of $0$ under a submersion $s$, then to perturb $f$ so that for
each $u$, the critical points of $s\circ f_u$ are of finite singularity
type. In fact, this is exactly the definition of what it means for $f_u$ to
be weakly transverse to $C$. In particular, when $C$ is a point, the
submersion can be taken to be the identity, so we have:
\begin{proposition} Let $f\colon M\to \R$ be smooth. If $f$ is weakly
transverse to a point $r\in \R$, then at each critical point in
$f^{-1}(r)$, the germ of $f$ has finite singularity type.
\label{prop:FST_weaktransverse}
\end{proposition}

Let $f$ and $g$ be germs of smooth maps from $(\R^m,a)$ to $(\R^p,f(a))$.
They are said to be \emph{$\mathcal{A}$-equivalent} if there exist a germ
$\varphi_1$ of a diffeomorphism of $(\R^m,a)$ and a germ $\varphi_2$ of a
diffeomorphism of $(\R^p,f(a))$ such that $g=\varphi_2\circ f\circ
\varphi_1$. If $\varphi_2$ can be taken to be the identity, then $f$ and
$g$ are called \emph{$\mathcal{R}$-equivalent} (for
\emph{right-equivalent}). There is also a notion of contact equivalence,
denoted by $\mathcal{K}$-equivalence, whose definition is readily
available, for example in \cite{Wall}. It is implied by
$\mathcal{A}$-equivalence.

We use $j^kf$ to denote the $k$-jet of $f$; for fixed coordinate systems at
points $a$ and $f(a)$ this is just the Taylor polynomial of $f$ of degree
$k$. For $\mathcal{G}$ one of $\mathcal{A}$, $\mathcal{K}$, or
$\mathcal{R}$, one says that $f$ is finitely $\mathcal{G}$-determined if
there exists a $k$ so that any germ $g$ with $j^kg=j^kf$ must be
$\mathcal{G}$-equivalent to $f$.  In particular, if $f$ is finitely
$\mathcal{G}$-determined, then for any fixed choice of coordinates at $a$
and $f(a)$, $f$ is $\mathcal{G}$-equivalent to a polynomial.

The elaborate theory of singularities of maps from $\R^m$ to $\R^p$
simplifies considerably when $p=1$.
\begin{lemma} Let $f$ be the germ of a map from $(\R^m,0)$ to $(\R,0)$,
with $0$ is a critical point of $f$. The following are equivalent.
\begin{enumerate}
\item[(i)] $f$ has finite singularity type.
\item[(ii)] $f$ is finitely $\mathcal{A}$-determined.
\item[(iii)] $f$ is finitely $\mathcal{R}$-determined.
\item[(iv)] $f$ is finitely $\mathcal{K}$-determined.
\end{enumerate}
\label{lem:FST}
\end{lemma}

\begin{proof}
In all dimensions, $f$ is finitely $\mathcal{K}$-determined if and only if
it is of finite singularity type (Corollary~III.6.9 of~\cite{GWLP}, or
alternatively the definition of finite singularity type on
\cite[p.~117]{Bruce} is exactly the condition given in Proposition~(3.6)(a)
of \cite{Mather} for $f$ to be finitely $\mathcal{K}$-determined).
Therefore (i) is equivalent to (iv).  Trivially (ii) implies (iii), and
(iii) implies (iv), and by Corollary~2.13 of \cite{Wall}, (iv)
implies~(ii).
\end{proof}

\subsection{Semialgebraic sets}
Recall (see for example Chapter~I.2 of \cite{GWLP}) that the class of
\emph{semialgebraic} subsets of $\R^m$ is defined to be the smallest
Boolean algebra of subsets of $\R^m$ that contains all sets of the form
$\{x\in \R^m\;\vert\;p(x)>0\}$ with $p$ a polynomial on $\R^m$. The
collection of semialgebraic subsets of $\R^m$ is closed under finite
unions, finite intersections, products, and complementation. The inverse
image of a semialgebraic set under a polynomial mapping is semialgebraic. A
nontrivial fact is the Tarski-Seidenberg Theorem (theorem~II.2(2.1) of
\cite{GWLP}), which says that a polynomial image of a semialgebraic set is
a semialgebraic set. Here is an easy lemma that we will need later.

\begin{lemma} Let $S$ be a semialgebraic subset of $\R^n$. If $S$ has empty
interior, then $S$ is contained in the zero set of a nontrivial polynomial
in $\R^n$.
\par
\label{lem:semialgebraic}
\end{lemma}

\begin{proof}
Since the union of the zero sets of two polynomials is the zero set of
their product, it suffices to consider a single semialgebraic set of the
form $(\cap_{i=1}^r \{x\;\vert\;p_i(x)\geq 0\})\cap (\cap_{j=1}^s
\{x\;\vert\;q_j(x) > 0\})$ where $p_i$ and $q_j$ are nontrivial
polynomials. We will show that if $S$ is of this form and has empty
interior, then $r\geq 1$ and $S$ is contained in the zero set of
$\prod_{i=1}^r p_i$. Suppose that $x\in S$ but all $p_i(x)>0$. Since all
$q_j(x)>0$ as well, there is an open neighborhood of $x$ on which all $p_i$
and all $q_j$ are positive. But then, $S$ has nonempty interior.
\end{proof}

\subsection{The codimension of a real-valued function}
\label{subsec:Sergeraert1}

It is, of course, fundamentally important that the Morse functions form an
open dense subset of $C^\infty(M,\R)$, the smooth maps from a closed
connected manifold $M$ of dimension $m$ to $\R$, with the
$C^\infty$-topology. But a great deal can also be said about the non-Morse
functions. There is a ``natural'' stratification of $C^\infty(M,\R)$ by
subsets $\mathcal{F}_i$, where stratification here means that the
$\mathcal{F}_i$ are disjoint subsets such that for every $n$ the union
$\cup_{i=0}^n \mathcal{F}_i$ is open. The functions in $\mathcal{F}_n$ are
those of ``codimension'' $n$, which we will define below. In particular,
$\mathcal{F}_0$ is exactly the open dense subset of Morse functions.

The union $\cup_{i=0}^\infty \mathcal{F}_i$ is not all of
$C^\infty(M,\R)$. However, the residual set
$C^\infty(M,\R)-\cup_{i=0}^\infty \mathcal{F}_i$ is of ``infinite
codimension,'' and any parameterized family of maps $F\colon M\times U\to
\R$ can be perturbed so that each $F_u$ is of finite codimension. In fact,
by applying theorem~\ref{thm:Bruce2} to the trivial bundle $1_{\R}\colon
\R\to \R$ and noting proposition~\ref{prop:FST_weaktransverse}, we may
perturb any parameterized family so that each $F_u$ is of finite
singularity type at each of its critical points. The definition of $f\in
C^\infty(M)$ being of finite codimension, given below, is exactly
equivalent to the algebraic condition given in (3.5) of \cite{Mather} for
$f$ to be finitely $\mathcal{A}$-determined at each of its critical points
(as noted in \cite{Mather}, this part of (3.5) was first due to Tougeron
\cite{Tougeron1}, \cite{Tougeron2}). By lemma~\ref{lem:FST}, this is
equivalent to $f$ having finite singularity type at each of its critical
points.  We summarize this as
\begin{proposition} A map $f\in C^\infty(M,\R)$ is of finite codimension 
if and only if it has finite singularity type at each of its critical
points.
\label{prop:FST_codimension}
\end{proposition}

We now recall material from section~7 of \cite{Sergeraert}. Denote the
smooth sections of a bundle $E$ over $M$ by $\Gamma(E)$. Until we reach
theorem~\ref{thm:Morse weak transversality}, we will denote
$C^\infty(M,\R)$ by $C(M)$. For a compact subset $K\subset \R$, define
$\Diff_K(\R)$ to be the diffeomorphisms of $\R$ supported on $K$.

Fix an element $f\in C(M)$ and a compact subset $K\subset \R$ for which
$f(M)$ lies in the interior of $K$. Define $\Phi\colon \Diff(M)\times
\Diff_K(\R)\to C(M)$ by $\Phi(\varphi_1,\varphi_2)=\varphi_2\circ
f\circ\varphi_1$. The differential of $\Phi$ at $(1_M,1_{\R})$ is defined
by $D(\xi_1,\xi_2)=f_*\xi_1+\xi_2\circ f$. Here, $\xi_1\in \Gamma(TM)$,
which is regarded as the tangent space at $1_M$ of $\Diff(M)$, $\xi_2\in
\Gamma_K(T\R)$, similarly identified with the tangent space at $1_{\R}$ of
$\Diff_K(\R)$, and $f_*\xi_1+\xi_2\circ f$ is regarded as an element of
$\Gamma(f^*T\R)$, which is identified with $C(M)$.  The \emph{codimension}
$\cdim(f)$ of $f$ is defined to be the real codimension of the image of $D$
in~$C(M)$. As will be seen shortly, the codimension of $f$ tells the real
codimension of the $\Diff(M)\times \Diff_K(\R)$-orbit of $f$ in~$C(M)$.

Suppose that $f$ has finite codimension $c$. In section~7.2 of
\cite{Sergeraert}, a method is given for computing $\cdim(f)$ using the
critical points of $f$.  Fix a critical point $a$ of $f$, with critical
value $f(a)=b$.  Consider $D_a\colon \Gamma_a(TM) \times C_b(\R)\to
C_a(M)$, where a subscript as in $\Gamma_a(TM)$ indicates the germs at $a$
of $\Gamma(TM)$, and so on. Notice that the codimension of the image
of $D_a$ is finite, indeed it is at most $c$.

Let $A$ denote the ideal $f_*\Gamma_a(TM)$ of $C_a(M)$. This can be
identified with the ideal in $C_a(M)$ generated by the partial derivatives
of $f$. An argument using Nakayama's lemma \cite[p.~645]{Sergeraert} shows
that $A$ has finite codimension in $C_a(M)$, and that some power of
$f(x)-f(a)$ lies in $A$. Define $\cdim(f,a)$ to be the dimension of
$C_a(M)/A$, and $\dim(f,a,b)$ to be the smallest $k$ such that
$(f(x)-f(a))^k\in A$.

Here is what these are measuring. The ideal $A$ tells what local
deformations of $f$ at $a$ can be achieved by precomposing $f$ with a
diffeomorphism of $M$ (near $1_M$), thus $\cdim(f,a)$ measures the
codimension of the $\Diff(M)$-orbit of the germ of $f$ at $a$.  The
additional local deformations of $f$ at $a$ that can be achieved by
postcomposing with a diffeomorphism of $\R$ (again, near $1_{\R}$) reduce
the codimension by $k$, basically because Taylor's theorem shows that the
germ at $a$ of any $\xi_2(f(x))$ can be written in terms of the powers
$(f(x)-f(a))^i$, $i<k$, plus a remainder of the form $K(x)(f(x)-f(a))^k$,
which is an element of the ideal $A$. Thus $\cdim(f,a)-\dim(f,a,b)$ is the
codimension of the image of $D_a$. For a noncritical point or a stable
critical point such as $f(x,y)=x^2-y^2$ at $(0,0)$, this local codimension
is $0$, but for unstable critical points it is positive.

Now, let $\dim(f,b)$ be the maximum of $\dim(f,a,b)$, taken over the
critical points $a$ such that $f(a)=b$ (put $\dim(f,b)=0$ if $b$ is not a
critical value). The codimension of $f$ is then $\sum_{a\in
M}\cdim(f,a)-\sum_{b\in\R}\dim(f,b)$.

Here is what is happening at each of the finitely many critical values $b$
of $f$. Let $a_1,\ldots\,$, $a_\ell$ be the critical points of $f$ with
$f(a_i)=b$. Let $f_i$ be the germ of $f-f(a_i)$ at $a_i$, and consider the
element $(f_1,\ldots,f_\ell)\in C_{a_1}(M)/A_1\oplus\cdots\oplus
C_{a_\ell}(M)/A_\ell$. The integer $\dim(f,b)$ is the smallest power of
$(f_1,\ldots,f_\ell)$ that is trivial in $C_{a_1}(M)/A_1\oplus\cdots\oplus
C_{a_\ell}(M)/A_\ell$.  The sum $\sum_i \cdim(f,a_i)$ counts how much
codimension of $f$ is produced by the inability to achieve local
deformations of $f$ near the $a_i$ by precomposing with local
diffeomorphisms at the $a_i$. This codimension is reduced by $\dim(f,b)$,
because the germs of the additional deformations that can be achieved by
postcomposition with diffeomorphisms of $\R$ near $b$ are the linear
combinations of $(1,\ldots,1)$, $(f_1,\ldots,f_\ell)$,
$(f_1^2,\ldots,f_\ell^2),\ldots\,$, $(f_1^{k-1},\ldots,f_\ell^{k-1})$. Thus
the contribution to the codimension from the critical points that map to
$b$ is $\sum_i\cdim(f,a_i)-\dim(f,b)$, and summing over all critical values
gives the codimension of $f$.

\subsection{The stratification of $C^\infty(M,\R)$ by codimension}
\label{subsec:Sergeraert2}

The functions whose codimension is finite and equal to $n$ form the stratum
$\mathcal{F}_n$. In particular, $\mathcal{F}_0$ are the Morse functions,
$\mathcal{F}_1$ are the functions either having all critical points stable
and exactly two with the same critical value, or having distinct critical
values and all critical points stable except one which is a birth-death
point. Moving to higher strata occurs either from more critical points
sharing a critical value, or from the appearance of more singularities of
positive but still finite local codimension.

We use the natural notations $\mathcal{F}_{\geq n}$ for $\cup_{i\geq
n}\mathcal{F}_i$, $\mathcal{F}_{> n}$ for $\cup_{i > n}\mathcal{F}_i$,
and so on. In particular, $\mathcal{F}_{\geq 0}$ is the set of all
elements of $C(M)$ of finite codimension, and $\mathcal{F}_{> 0}$ is the
set of all elements of finite codimension that are not Morse functions.

The main results of \cite{Sergeraert} (in particular, Theorem~8.1.1 and
Theorem~9.2.4) show that the Sergeraert stratification is locally trivial,
in the following sense.
\begin{theorem}[Sergeraert] Suppose that
$f\in \mathcal{F}_n$. Then there is a neighborhood $V$ of $f$ in
$C(M)$ of the form $U\times \R^n$, where
\begin{enumerate}
\item $U$ is a neighborhood of $1$ in $\Diff(M)\times \Diff_K(\R)$, and
\item
there is a stratification $\R^n=\cup_{i=0}^n F_i$, such that
$\mathcal{F}_i\cap V=U\times F_i$.
\end{enumerate}
\label{thm:Sergeraert}
\end{theorem}

The inner workings of this result are as follows. Select elements
$f_1,\ldots,\,$ $f_n\in C(M)$ that represent a basis for the quotient of
$C(M)$ by the image of the differential $D$ of $\Phi$ at $(1_M,1_{\R})$.  For
$x=(x_1,\ldots,x_n)\in \R^n$, the function $g_x=f+\sum_{i=1}^n x_i f_i$ is
an element of $C(M)$. If the $x_i$ are chosen in a sufficiently small ball
around $0$, which is again identified with $\R^n$, then these $g_x$ form a
copy $E$ of $\R^n$ ``transverse'' to the image of $\Phi$. Then, $F_i$ is
defined to be the intersection $E\cap \mathcal{F}_i$. A number of subtle
results on this local structure and its relation to the action of
$\Diff(M)\times \Diff_K(\R)$ are obtained in \cite{Sergeraert}, but we will
only need the local structure we have described here.

We remark that $F_n$ is not necessarily just $\{0\}\in \R^n$, that is, the
orbit of $f$ under $\Diff(M)\times \Diff_K(\R)$ might not fill up the
stratum $\mathcal{F}_n$ near $f$. This result, due to H. Hendriks
\cite{Hendriks}, has been interpreted as saying that the Sergeraert
stratification of $C(M)$ is \emph{not} locally trivial (a source of some
confusion), or that it is ``pathological'' (which we find far too
pejorative). 

Denoting $\cup_{i\geq 1}F_i$ by $F_{\geq 1}$, we have the following key
technical result.

\begin{proposition}
For some coordinates on $E$ as $\R^n$, there are a neighborhood $L$ of
$0$ in $\R^n$ and a nonzero polynomial $p$ on $\R^n$ such that $p(L\cap
F_{\geq 1})=0$.
\label{prop:Sergeraert}
\end{proposition}
\begin{proof}
We begin with a rough outline of the proof. By lemma~\ref{lem:FST}, we may
choose local coordinates at the critical points of $f$ for which
$f$ is polynomial near each
critical point. We will select the $f_i$ in the construction of the
transverse slice $E=\R^n$ to be polynomial on these neighborhoods.  Now
$F_{\geq 1}$ consists exactly of the choices of parameters $x_i$ for which $f+\sum
x_if_i$ is not a Morse function, since they are the intersection of $E$
with $\mathcal{F}_{\geq 1}$. We will show that they form a semialgebraic
set. But $F_{\geq1}$ has no interior, since otherwise (using
theorem~\ref{thm:Sergeraert}) the subset of Morse functions $\mathcal{F}_0$
would not be dense in $C(M)$. So lemma~\ref{lem:semialgebraic} show that
$F_{\geq 1}$ lies in the zero set of some nontrivial polynomial.

Now for the details. Recall that $m$ denotes the dimension of $M$. Consider
a single critical value $b$, and let $a_1,\ldots\,$, $a_\ell$ be the
critical points with $f(a_i)=b$. Fix coordinate neighborhoods $U_i$ of the
$a_i$ with disjoint closures, so that $a_i$ is the origin $0$ in $U_i$.  By
lemma~\ref{lem:FST}, $f$ is finitely $\mathcal{R}$-determined near each
critical point, so on each $U_i$ there is a germ $\varphi_i$ of a
diffeomorphism at $0$ so that $f\circ \varphi_i$ is the germ of a
polynomial. That is, by reducing the size of the $U_i$ and changing the
local coordinates, we may assume that on each $U_i$, $f$ is a polynomial
$p_i$. As explained in subsection~\ref{subsec:Sergeraert1}, the
contribution to the codimension of $f$ from the $a_i$ is the dimension of
the quotient
\[ Q_b=\big(\oplus_{i=1}^\ell C_{a_i}(U_i)/A_i\big)/B\]
where $B$ is the vector subspace spanned by $\{1,
(p_1(x)-b,\ldots,p_\ell(x)-b),\ldots,\allowbreak
((p_1(x)-b)^{k-1},\ldots,(p_\ell(x)-b)^{k-1})\}$. Choose $q_{i,j}$, $1\leq
j\leq n_i$, where $q_{i,j}$ is a polynomial on $U_i$, so that the germs of
the $q_{i,j}$ form a basis for $Q_b$. Fix vector spaces $\Lambda_i\cong
\R^{n_i}=\{(x_{i,1},\ldots,x_{i,n_i})\}$; these will eventually be some of
the coordinates on~$E$.

In each $U_i$, select round open balls $V_i$ and $W_i$ centered at $0$ so
that $W_i\subset\overline{W_i} \subset V_i\subset \overline{V_i}\subset
U_i$.  We select them small enough so that the closures in $\R$ of their
images under $f$ do not contain any critical value except for $b$. Fix a
smooth function $\mu\colon M\to [0,1]$ which is $1$ on $\cup
\overline{W_i}$ and is $0$ on $M-\cup V_i$, and put $f_{i,j}=\mu\cdot
q_{i,j}$, a smooth function on all of $M$. Now choose a product $L=\prod_i
L_i$, where each $L_i$ is a round open ball centered at $0$ in $\Lambda_i$,
small enough so that if each $(x_{i,1},\ldots, x_{i,n_i})\in L_i$, then
each critical point of $f+\sum x_{i,j} f_{i,j}$ either lies in $\cup W_i$,
or is one of the original critical points of $f$ lying outside of~$\cup
U_i$.

We repeat this process for each of the finitely many critical values of
$f$, choosing additional $W_i$ and $L_i$ so small that all critical points
of $f+\sum x_{i,j} f_{i,j}$ lie in $\cup W_i$. That is, these perturbations
of $f$ are so small that each of the original critical points of $f$ breaks
up into critical points that lie very near the original one and far from
the others.

The sum of all $n_i$ is now $n$. We again use $\ell$ for the number of
$U_i$, and write $\Lambda$ and $L$ for the direct sum of all the
$\Lambda_i$ and the product of all the $L_i$ respectively. For $x\in
L$, write $g_x=f+\sum x_{i,j}f_{i,j}$. It remains to show that the set of
$x$ for which $g_x$ is not a Morse function--- that is, has a critical
point with zero Hessian or has two critical points with the same value---
is contained in a union of finitely many semialgebraic sets.

Denote elements of $W_i$ by $\overline{u_i} = (u_{i,1},\ldots,u_{i,m})$,
and similarly for elements $\overline{x_i}$ of $L_i$.  Define $G_i\colon
W_i\times L_i\to \R$ by $G_i(\overline{u_i},\overline{x_i}) =
p_i(\overline{u_i})+\sum_{j=1}^{n_i}x_{i,j}q_{i,j}(\overline{u_i})$. Note
that for $x=(\overline{x_1},\ldots,\overline{x_\ell})$,
$(G_i)_{\overline{x_i}}$ is exactly the restriction of $g_x$ to $W_i$.

We introduce one more notation that will be convenient.  For $X\subseteq
L_i$ define $E(X)$ to be the set of all
$(\overline{x_1},\ldots,\overline{x_\ell})$ in $L$ such that
$\overline{x_i}\in X$. When $X$ is a semialgebraic subset of $L_i$, $E(X)$
is a semialgebraic subset of $L$. Similarly, if $X\times Y\subseteq
L_i\times L_j$, we use $E(X\times Y)$ to denote its extension to a subset
of $L$, that is, $E(X)\cap E(Y)$.

For each $i$, let $S_i$ be the set of all $(\overline{u_i},\overline{x_i})$
in $W_i\times L_i$ such that that $\partial G_i/\partial u_{i,j} $ for
$1\leq j\leq n_i$ all vanish at $(\overline{u_i},\overline{x_i})$, that is,
the pairs such that $\overline{u_i}$ is a critical point of
$(G_i)_{\overline{x_i}}$.  Since $S_i$ is the intersection of an algebraic
set in $\R^m\times \R^{n_i}$ with $W_i\times L_i$, and the latter are round
open balls, $S_i$ is semialgebraic. Let $H_i$ be the set of all
$(\overline{u_i},\overline{x_i})$ in $W_i\times L_i$ such that the Hessian
of $(G_i)_{\overline{x_i}}$ vanishes at $\overline{u_i}$, again a
semialgebraic set. The intersection $H_i\cap S_i$ is the set of all
$(\overline{u_i},\overline{x_i})$ such that $(G_i)_{\overline{x_i}}$ has an
unstable critical point at $\overline{u_i}$. By the Tarski-Seidenberg
Theorem, its projection to $L_i$ is a semialgebraic set, which we will
denote by $A_i$. The union of the $E(A_i)$, $1\leq i\leq \ell$, is
precisely the set of $x$ in $L$ such that $g_x$ has an unstable critical
point.

Now consider $G_i\times G_i\colon S_i\times S_i-\Delta_i\to \R^2$, where
$\Delta_i$ is the diagonal in $S_i\times S_i$.  Let
$\widetilde{B_i}=(G_i\times G_i)^{-1}(\Delta_2)$, where $\Delta_2$ is the
diagonal of $\R^2$. Now, let $\Delta_i'$ be the set of all
$((\overline{u_i},\overline{x_i}), (\overline{u_i}',\overline{x_i}'))$ in
$W_i\times L_i\times W_i\times L_i$ such that
$\overline{x_i}=\overline{x_i}'$. Then the projection of
$\widetilde{B_i}\cap \Delta_i'$ to its first two coordinates is the set of
all $(\overline{u_i},\overline{x_i})$ in $W_i\times L_i$ such that
$\overline{u_i}$ is a critical point of $(G_i)_{\overline{x_i}}$ and
$(G_i)_{\overline{x_i}}$ has another critical point with the same value.
The projection to the second coordinate alone is the set $B_i$ of
$\overline{x_i}$ for which $(G_i)_{\overline{x_i}}$ has two critical points
with the same value.

Finally, for $i\neq j$, consider $G_i\times G_j\colon S_i\times S_j\to
\R^2$ and let $\widetilde{B_{i,j}}$ be the preimage of $\Delta_2$. Let
$B_{i,j}$ be the projection of $\widetilde{B_{i,j}}$ to a subset of
$L_i\times L_j$.  The union of the $E(B_i)$ and the $E(B_{i,j})$ is
precisely the set of all $x$ such that $g_x$ has two critical points with
the same value. Since these are semialgebraic sets, the proof is complete.
\end{proof}

Here is the main result of this subsection.
\begin{theorem}
Let $M$ and $W$ be compact smooth manifolds. Then for a residual set of
smooth maps $F$ from $I\times W$ to $C^\infty(M,\R)$, the following hold.
\begin{enumerate}
\item[(i)] $F(I\times W)\subset \mathcal{F}_{\geq 0}$.
\item[(ii)] Each $F_u^{-1}(\mathcal{F}_{>0})$ is finite.
\end{enumerate}
\label{thm:Morse weak transversality}
\end{theorem}

\begin{proof}
Start with a smooth map $G\colon I\times W\to C^\infty(M,\R)$.  Regarding
it as a parameterized family of maps $M\times (I\times W)\to \R$, we apply
theorem~\ref{thm:Bruce2} to perturb $G$ so that each $G_u$ is weakly
transverse to the points of $\R$. By
proposition~\ref{prop:FST_codimension}, this implies that $G(I\times
W)\subset \mathcal{F}_{\geq 0}$. Since $I\times W$ is compact, $G(I\times
W)\subset \mathcal{F}_{\leq n}$ for some~$n$.

For each $f\in \mathcal{F}_{>0}$, choose a neighborhood $V_f=U_f\times \R^n$ as in
theorem~\ref{thm:Sergeraert}. Using proposition~\ref{prop:Sergeraert},
select a neighborhood $L_f$ of $0$ in $\R^n$ and a nonzero polynomial
$p_f\colon L_f\to \R$ such that $p_f(L\cap F_{i\geq 1})=0$.

Now, partition $I$ into subintervals and triangulate $W$ so that for each
subinterval $J$ and each simplex $\Delta$ of maximal dimension in the
triangulation, $G(J\times \Delta)$ lies either in $\mathcal{F}_0$ or in
some $U_f\times L_f$. Fix a particular $J\times \Delta$.  If $G(J\times
\Delta)$ lies in $\mathcal{F}_0$, do nothing. If not, choose $f$ so that
$G(J\times \Delta)$ lies in $U_f\times L_f$.  Let $\pi\colon U_f\times
L_f\to L_f$ be the projection, so that $p_f\circ \pi(U_f\times
F_{\geq 1})=0$. By proposition~\ref{prop:polynomial_Bruce}, we
may perturb $G\vert_{J\times \Delta}$ (changing only its $L_f$-coordinate
in $U_f\times L_f$) so that for each $u\in \Delta$,
$G_u\vert_J^{-1}(\mathcal{F}_{i\geq 1})$ is finite, and any map
sufficiently close to $G\vert_{J\times \Delta}$ on $J\times \Delta$ will
have this same property. As usual, of course, this is extended to a
perturbation of~$G$.

%This perturbation of $G$ can be tapered off near $J\times \Delta$, changing
%$G$ only on the preimage of $L_f\times U_f$.

This process can be repeated sequentially on the remaining $J\times
\Delta$. The perturbations must be so small that the property of having
each $G_u\vert_J^{-1}(\mathcal{F}_{i\geq 1})$ finite is not lost on
previously considered sets. When all $J\times \Delta$ have been considered,
each $G_u^{-1}(\mathcal{F}_{i\geq 1})$ is finite.
\end{proof}

\subsection{Border label control}
\label{sec:borderlabel}
We now return to the case when $M$ is a closed $3$-manifold, as in the
introduction of section~\ref{sec:generalposition}. In this subsection, we
will obtain a deformation of $f\colon M\times W\to M$ for which some
$\epsilon$ gives border label control.

%We recall from subsection~\ref{subsec:PEP} that any sufficiently close
%approximation of $f$ can be obtained by a deformation of $f$ through
%parameterized families of diffeomorphisms.  Our procedure will require many
%deformations of the initial $f$, each of which must be small enough to
%maintain all previously achieved properties. We use the terms ``perturb''
%and ``perturbation'' to signal this.  The new maps will be called $f$
%again, to avoid a proliferation of notation.

We begin by ensuring that no $f_u$ carries a component of the singular set
$T$ of $\tau$ into $T$. Consider two circles $C_1$ and $C_2$ imbedded in
$M$. By theorem~\ref{thm:Bruce1}, applied with $A=C_1\times W$, $B=M$, and
$C=C_2$, we may perturb $f\vert_{C_1\times W}$ so that for each $u\in W$,
$f_u\vert_{C_1}$ meets $C_2$ in only finitely many points.

Recall that $T$ consists of smooth circles and arcs in $M$. Each arc is
part of some smoothly imbedded circle, so $T$ is contained in a union
$\cup_{i=1}^n C_i$ of imbedded circles in $M$. By a sequence of
perturbations as above, we may assume that at each $u$, each $f_u(C_i)$
meets each $C_j$ in a finite set (including when $i=j$), so that $f_u(T)$
meets $T$ in a finite set.

The next potential problem is that at some $u$, $f_u(T_0)$ or $f_u(T_1)$
might be contained in a single level $P_t$. Recall that the notation
$R(s,t)$, introduced in section~\ref{sec:RSgraphic}, means
$\tau^{-1}([s,t])$.  For some $\delta>0$, every $f_u(T_0)$ meets
$R(3\delta,1-3\delta)$, since otherwise the compactness of $W$ would lead
to a parameter $u$ for which $f_u(T_0)\subset T$.  Let $\phi\colon
R(\delta,1-\delta)\to [\delta,1-\delta]$ be the restriction of the map
$\pi(\tau(x,t))=t$. This $\phi$ makes $R(\delta,1-\delta)$ a bundle with
fibers that are level tori. As before, let $C_1$ be one of the circles
whose union contains $T$. Only the most superficial changes are needed to
the proof of theorem~\ref{thm:Bruce2} given in \cite{Bruce} so that it
applies when $\phi$ is a bundle map defined on a codimension-zero
submanifold of $B$ rather than on all of $B$; the only difference is that
the subsets of jets which are to be avoided are defined only at points of
the subspace rather than at every point of $B$. Using this slight
generalization of theorem~\ref{thm:Bruce2} (and as usual, the Parameterized
Extension Principle), we perturb $f$ so that each $f_u\vert_{C_1}$ is
weakly transverse to each $P_t$ with $\delta\leq t\leq 1-\delta$. Since
$C_1$ is $1$-dimensional, weakly transverse implies that $f_u(C_1)$ meets
each such $P_t$ in only finitely many points.  Repeating for the other
$C_i$, we may assume that each $f_u(T_0)$ meets the $P_t$ with $\delta\leq
t \leq 1-\delta$ in only finitely many points. We also choose the
perturbations small enough so that each $f_u(T_0)$ still meets
$R(2\delta,1-2\delta)$. So $f_u^{-1}(P_t)\cap T_0$ is nonempty and finite
at least some~$t$.  In particular, $\pi(f_u(T_0))$ contains an open set, so
by Sard's Theorem applied to $\pi\circ f_u\vert_{T_0}$, for each $u$, there
is an $r$ so that $f_u(T_0)$ meets $P_r$ transversely in a nonempty set (we
select $r$ so that $P_r$ does not contain the image of a vertex of
$T_0$). For a small enough $\epsilon$, a component of $X_s\cap P_r$ will be
a compressing disk of $X_s$ whenever $s\leq 2\epsilon$, and by compactness
of $W$, there is an $\epsilon$ such that for every $u$, there is a level
$P_r$ such that some component of $X_s\cap P_r$ contains a compressing disk
of $X_s$ whenever $s\leq 2\epsilon$.

Applying the same procedure to $T_1$, we may assume that for every $u$,
there there is a level $P_r$ such that some component of $Y_s\cap P_r$ is a
compressing disk of $Y_s$ whenever $s\geq 1-2\epsilon$.

Let $h\colon M\times W\to M$ be defined by $h(x,u)=f_u^{-1}(x)$. Applying
the previous procedure to $h$, making sure that all perturbations are small
enough to preserve the conditions developed for $f$, and perhaps making
$\epsilon$ smaller, we may assume that for each $u$, there is a level $Q_r$
such that $V_t\cap Q_r$ is a compressing disk of $V_t$ whenever $t\leq
2\epsilon$, and a similar $Q_r$ for $W_t$ with $t\geq 1-2\epsilon$. Thus
the number $\epsilon$ gives border label control for $f$. Since border
label control holds, with the same $\epsilon$, for any map sufficiently
close to $f$, we may assume it is preserved by all future perturbations.
\longpage 

\subsection{Building the graphics}
It remains to deform $f$ to satisfy conditions (GP1), (GP2), and (GP3).  As
before, let $i\colon I\to \R$ be the inclusion, and consider the smooth map
$i\circ \pi \circ f\circ (\tau\times 1_W)\colon P\times I\times W\to \R$.
Regard this as a family of maps from $I$ to $C^\infty(P,\R)$, parameterized
by $W$. Apply theorem~\ref{thm:Morse weak transversality} to obtain a
family $k\colon P\times I\times W\to \R$. For each $I\times \{u\}$, there
will be only finitely many values of $s$ in $I$ for which the restriction
$k_{(s,u)}$ of $k$ to $P\times \{s\}\times \{u\}$ is not a Morse
function. At these levels, the projection from $Q_s$ into the transverse
direction to $P_t$ is an element of some $\mathcal{F}_n$, so each tangency
of $Q_s$ with $P_t$ looks like the graph of a critical point of finite
multiplicity. This will ultimately ensure that condition (GP1) is attained
when we complete our deformations of~$f$.

We will use $k$ to obtain a deformation of the original $f$, by moving
image points vertically with respect to the levels of the range. This would
not make sense where the values of $k$ fall outside $(0,1)$, so the motion
will be tapered off so as not to change $f$ at points that map near $T$. It
also would not be well-defined at points of $T\times W$, so we taper off
the deformation so as not to change $f$ near $T\times W$. The fact that $f$
is unchanged near $T\times W$ and near points that map to $T$ will not
matter, since border label control will allow us to ignore these regions in
our later work.

Regard $P\times I\times W$ as a subspace of $P\times \R\times W$. For
each $(x,r,u)\in P\times I\times W$, let $w_{(x,r,u)}'$ be
$k(x,r,u)-i\circ \pi \circ f_u\circ \tau(x,r)$, regarded as a tangent vector to
$\R$ at $i\circ \pi\circ f_u\circ \tau(x,r)$.

We will taper off the $w_{(x,r,u)}'$ so that for each fixed $u$ they will
produce a vector field on $M$. Fix a number $\epsilon$ that gives border
label control for $f$, and a smooth function $\mu\colon \R\to I$ which
carries $(-\infty,\epsilon/4]\cup [1-\epsilon/4,\infty)$ to $0$ and carries
$[\epsilon/2,1-\epsilon/2]$ to $1$. Define $w_{(x,r,u)}$ to be
$\mu(r)\,\mu(i\circ \pi \circ f_u\circ \tau(x,r))\,w_{(x,r,u)}'$. These
vectors vanish whenever $r\notin [\epsilon/4,1-\epsilon/4]$ or $i\circ
\pi\circ f_u\circ \tau(x,r,u)\notin [\epsilon/4,1-\epsilon/4]$, that is,
whenever $\tau(x,r)$ or $f_u\circ \tau(x,r)$ is close to $T$.  Using the
map $i\circ \pi\colon M\to \R$, we pull these back to vectors in $M$ that
are perpendicular to $P_t$; this makes sense near $T$ since the
$w_{(x,r,u)}'$ are zero at these points). For each $u$, we obtain at each
point $f_u\circ \tau(x,r)\in M$ a vector $v_{(x,r,u)}$ that points in the
$I$-direction (i.~e.~is perpendicular to $P_t$) and maps to $w_{(x,r,u)}$
under~$(i\circ \pi)_*$.

If $k$ was a sufficiently small perturbation, the $v_{(x,r,u)}$ define a
smooth map $j_u\colon M\to M$ by $j_u(\tau(x,r))=\Exp(v_{(x,r,u)})$. Put
$g_u=j_u\circ f_u$.  Since $\mu(r)=1$ for $\epsilon/2\leq r \leq
1-\epsilon/2$, we have $i\circ \pi \circ g_u\circ \tau(x,r)=k(x,r,u)$
whenever both $\epsilon/2\leq r\leq 1-\epsilon/2$ and $\epsilon/2\leq
i\circ \pi \circ f_u\circ \tau(x,r)\leq 1-\epsilon/2$. The latter condition
says that $f_u\circ \tau(x,r)$ is in $P_s$ for some $\epsilon/2\leq s \leq
1-\epsilon/2$. Assuming that $k$ was close enough to $i\circ \pi \circ
f\circ (\tau\times 1_W)$ so that each $\pi\circ g_u\circ \tau(x,r)$ is
within $\epsilon/4$ of $\pi\circ f_u\circ \tau(x,r)$, the equality $i\circ
\pi \circ g_u\circ \tau(x,r)=k(x,r,u)$ holds whenever $\tau(x,r)$ is in a
$P_s$ and $g_u\circ \tau(x,r)$ is in a $P_t$ with $\epsilon\leq s,t\leq
1-\epsilon$.

Carrying out this construction for a sequence of $k$ that converge to
$i\circ \pi \circ f\circ (\tau\times 1_W)$, we obtain vector fields
$v_{(x,r,u)}$ that converge to the zero vector field. For those
sufficiently close to zero, $g$ will be a perturbation of $f$. Choosing $g$
sufficiently close to $f$, we may ensure that $\epsilon$ still gives border
label control for~$g$.

We will now analyze the graphic of $g_u$ on $I^2_\epsilon$. For $s,t\in
[\epsilon,1-\epsilon]$, $\pi\circ g_u(x)$ equals $k_{(s,u)}(x)$ whenever
$x\in P_s$ and $g_u(x)\in P_t$. Therefore the tangencies of $g_u(P_s)$ with
$P_t$ are locally just the graphs of a critical point of $k_{(s,u)}\colon
P\to \R$, so $g$ has property~(GP1).

\begin{figure}
\includegraphics[width=35 ex]{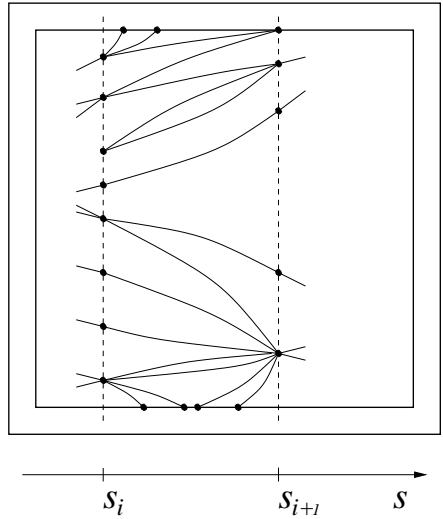}
\caption{A portion of the graphic of $g_u$.}
\label{fig:graphic}
\end{figure}

Let $s_1,\ldots\,$, $s_n$, be the values of $s$ in $[\epsilon,1-\epsilon]$
for which $k_{(s_i,u)}\colon P\to \R$ is not a Morse function.  Each
$k_{(s_i,u)}$ is still a function of finite codimension, so has finitely
many critical points. Those with critical values in $[\epsilon,1-\epsilon]$
produce the points of the graphic of $g_u$ that lie in the vertical line
$s=s_i$, as suggested in figure~\ref{fig:graphic}.  We declare the
$(s_i,t)$ at which $k_{(s_i,u)}$ has a critical point at $t$ to be vertices
of $\Gamma_u$.

When $s$ is not one of the $s_i$, $k_{(s,u)}$ is a Morse function, so any
tangency of $g_u(P_s)$ with $P_t$ is stable, and there is at most one such
tangency. Since these tangencies are stable, all nearby tangencies are
equivalent to them and hence also stable, so in the graphic for $g_u$ in
$I^2_\epsilon$, the pairs $(s,t)$ corresponding to levels with a single
stable tangency form ascending and descending arcs as suggested in
figure~\ref{fig:graphic}. These arcs may enter or leave $I^2_\epsilon$, or
may end at a point corresponding to one of the finitely many points of the
graphic with $s$-coordinate equal to one of the $s_i$.  We declare the
intersection points of these arcs with $\partial I_\epsilon$ to be vertices
of $\Gamma_u$.  The conditions (GP2) and (GP3) have been achieved,
completing the proof of theorem~\ref{thm:generalposition}.

\newpage
\section[Finding good regions]
{Finding good regions}
\label{sec:goodregions}

In this section, we adapt the arguments of
section~\ref{sec:Rubinstein-Scharlemann} to general position families.  The
graphics associated to the $f_u$ of a general position family $f\colon
M\times W\to M$ satisfy property (RS\ref{item:noAB}) of
section~\ref{sec:Rubinstein-Scharlemann} (provided that the Heegaard
splittings associated to the sweepout are strongly irreducible) and
property (RS\ref{item:edge labeling}) (since the open edges of the $\Gamma$
correspond to pairs of levels that have a single stable tangency, see the
remark after the definition of (RS\ref{item:edge labeling}) in
section~\ref{sec:Rubinstein-Scharlemann}), but not property
(RS\ref{item:2-cell labeling}). Indeed, property (RS\ref{item:2-cell
labeling}) does not even make sense, since the vertices of $\Gamma$ can
have high valence. Property (RS\ref{item:noAB}) is what allows the map
from the $0$-cells of $K$ to the $0$-simplices of the Diagram to be
defined. Property (RS\ref{item:edge labeling}) (plus conditions on regions
near $\partial K$, which we will still have due to border label control)
allows it to be extended to a cellular map from the $1$-skeleton of $K$ to
the $1$-skeleton of the Diagram. What ensures that it still extends to the
$2$-cells is a topological fact about pairs of levels whose intersection
contains a common spine, lemma~\ref{lem:common spine}. Because it involves
surfaces that do not meet transversely, its proof is complicated and
somewhat delicate. Since the proof does not introduce any ideas needed
elsewhere, the reader may wish to skip it on a first reading, and go
directly from the statement of lemma~\ref{lem:common spine} to the last
four paragraphs of the section.

We specialize to the case of a parameterized family $f\colon L\times W\to
L$, where $L$ is a lens space and $W$ is a compact manifold.  We retain the
notations $P_t$, $Q_s$, $V_t$, $W_t$, $X_s$, and $Y_s$ of
section~\ref{sec:generalposition}. As usual, only $P_t$, $V_t$, $W_t$, and
a number $\epsilon$ which gives border label control for $f$ are
independent of the parameter~$u$. As was mentioned above, properties
(RS\ref{item:noAB}) and (RS\ref{item:edge labeling}) still hold for the
labels of the regions of the graphic of each~$f_u$.

\begin{theorem}
Suppose that $f\colon L\times W\to L$ is in general position with respect
to~$\tau$. Then for each $u$, there exists $(s,t)$ such that $Q_s$ meets
$P_t$ in good position.\par
\label{thm:finding good levels}
\end{theorem}
\noindent The proof of theorem~\ref{thm:finding good levels} will
constitute the remainder of this section.

We begin by examining the labels of parameters near the boundary of
$I^2_\epsilon$; this will ultimately ensure that the boundary of
$I^2_\epsilon$ maps to the Diagram in an esssential way. Fix a parameter
$u$, and suppose that $(s,t)$ is a point in the interior of $I^2_\epsilon$
such that $Q_s$ meets $P_t$ transversely.  The next lemma is immediate from
the definition of border label control and the labeling rules for regions.
It does not require that we be working with lens spaces, so we state it as
a lemma with weaker hypotheses.
\begin{lemma}
Suppose that $f\colon M\times W\to M$ is in general position with respect
to $\tau$. Assume that $M\neq S^3$ and that the Heegaard splittings
associated to $\tau$ are strongly irreducible. Suppose that $\epsilon$
gives border label control for $f$.
\begin{enumerate}
\item
If $t\leq \epsilon$, then the label of $(s,t)$ contains $A$.
\item
If $t\geq 1-\epsilon$, then the label of $(s,t)$ contains $B$.
\item
If $s\leq \epsilon$, then the label of $(s,t)$ contains $X$.
\item
If $s\geq 1-\epsilon$, then the label of $(s,t)$ contains $Y$.
\end{enumerate}
\label{lem:borderlabelconstraints}
\end{lemma}

We now prove a key geometric lemma that is particular to lens spaces.
\begin{lemma}
Let $f\colon L \times W\to L$ be a parameterized family of diffeomorphisms
in general position, and let $(s,t)\in I^2_\epsilon$. If $Q_s\cap P_t$
contains a spine of $P_t$, then either $V_t$ or $W_t$ contains a core
circle which is disjoint from $Q_s$.
\label{lem:common spine}
\end{lemma}

\begin{proof}
We will move $Q_s$ by a sequence of isotopies. All isotopies will have the
property that if $V_t-Q_s$ (or $W_t-Q_s$) did not contain a core circle of
$V_t$ (or $W_t$) before the isotopy, then the same is true after the
isotopy.  We say this succinctly with the phrase that the isotopy
\emph{does not create core circles.} Typically some of the isotopies will
not be smooth, so we work in the PL category. At the end of an initial
``flattening'' isotopy, $Q_s$ will intersect $P_t$ nontransversely in a
$2$-dimensional simplicial complex $X$ in $P_t$ whose frontier consists of
points where $Q_s$ is PL imbedded but not smoothly imbedded.  A sequence of
simplifications called tunnel moves and bigon moves, plus isotopies that
push disks across balls, will make $Q_s\cap P_t$ a single component $X_0$,
which will then undergo a few additional improvements.  After this has been
completed, an Euler characteristic calculation will show that a core circle
disjoint from the repositioned $Q_s$ exists in either $V_t$ or $W_t$, and
consequently one existed for the original~$Q_s$.

Since $f$ is in general position, $Q_s\cap P_t$ is a $1$-complex satisfying
the property (GP1) of section~\ref{sec:generalposition}. Each isolated
vertex of $Q_s\cap P_t$ is an isolated tangency of $Q_s\cap P_t$, so we can
move $Q_s$ by a small isotopy near the vertex to eliminate it from the
intersection. After this step, $Q_s\cap P_t$ is a graph $\Gamma$ which
contains a spine of $Q_s\cap P_t$, such that each vertex of $\Gamma$ has
positive valence.

By property (GP1), each vertex $x$ of $\Gamma$ is a point where $Q_s$ is
tangent to $P_t$, and the edges of $\Gamma$ that emanate from $x$ are arcs
where $Q_s$ intersects $P_t$ transversely. Along each arc, $Q_s$ crosses
from $V_t$ into $W_t$ or vice versa, so there is an even number of these
arcs. Near $x$, the tangent planes of $Q_s$ are nearly parallel to those of
$P_t$, and there is an isotopy that moves a small disk neighborhood of $x$
in $Q_s$ until it coincides with a small disk neighborhood of $x$ in
$P_t$. Perform such isotopies near each vertex of~$\Gamma$. This enlarges
$\Gamma$ in $Q_s\cap P_t$ to the union of $\Gamma$ with a union $E$ of
disks, each disk containing one of the original vertices.

The closure of the portion of $\Gamma$ that is not in $E$ now consists of a
collection of arcs and circles where $Q_s$ intersects $P_t$ transversely,
except at the endpoints of the arcs, which lie in $E$. Consider one of
these arcs, $\alpha$. At points of $\alpha$ near $E$, the tangent planes to
$Q_s$ are nearly parallel to those of $P_t$, and starting from each end
there is an isotopy that moves a small regular neighborhood of a portion of
$\alpha$ in $Q_s$ onto a small regular neighborhood of the same portion of
$\alpha$ in $P_t$. This flattening process can be continued along
$\alpha$. If it is started from both ends of $\alpha$, it may be possible
to flatten all of a regular neighborhood of $\alpha$ in $Q_s$ onto one in
$P_t$. This occurs when the vectors in a field of normal vectors to
$\alpha$ in $Q_s$ are being moved to normal vectors on the same side of
$\alpha$ in $P_t$. If they are being moved to opposite sides, then we
introduce a point where the configuration is as in
figure~\ref{fig:flatten}, in which $P_t$ appears as the $xy$-plane,
$\alpha$ appears as the points in $P_t$ with $x=-y$, and $Q_s$ appears as
the four shaded half- or quarter-planes. These points will be called
\emph{crossover} points. Perform such isotopies in disjoint neighborhoods
of all the arcs of $\Gamma-E$. For the components of $\Gamma$ that are
circles of transverse intersection points, we flatten $Q_s$ near each
circle to enlarge it to an annulus.
\begin{figure}
\includegraphics[width=30 ex]{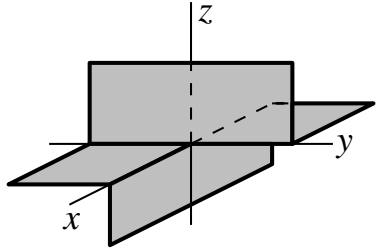}
\caption{The crossover configuration.}
\label{fig:flatten}
\end{figure}

At the end of this initial process, $\Gamma$ has been been enlarged to a
$2$-complex $X$ in $Q_s\cap P_t$ that is a regular neighborhood of
$\Gamma$, except at the crossover points where $\Gamma$ and $X$ look
locally like the antidiagonal $x=-y$ of the $xy$-plane and the set of
points with $xy\leq 0$. We will refer to $X$ as a \emph{pinched regular
neighborhood} of~$\Gamma$.

Since $\Gamma$ originally contained a spine of $P_t$, $X$ contains two
circles that meet transversely in one point that lies in the interior (in
$P_t$) of $X$. Therefore $X$ contains a common spine of $P_t$ and~$Q_s$.
Let $X_0$ be the component of $X$ that contains a common spine of $Q_s$ and
$P_t$. All components of $P_t-X_0$ and $Q_s-X_0$ are open disks. Let
$X_1=X-X_0$, and for each $i$, denote $\Gamma\cap X_i$ by~$\Gamma_i$.

The next step will be to move $Q_s$ by isotopy to remove $X_1$ from
$Q_s\cap P_t$. These isotopies will be fixed near $X_0$.  Some of them will
have the effect of joining two components of $V_t-Q_s$ (or of $W_t-Q_s$)
into a single component of $V_t-Q_s$ (or of $W_t-Q_s$) for the repositioned
$Q_s$, so we must be very careful not to create core circles.

The frontier of $X_1$ in $P_t$ is a graph $\Fr(X_1)$ for which each vertex
is a crossover point, and has valence $4$ (as usual, our ``graphs'' can
have open edges that are circles). Its edges are of two types: \emph{up}
edges, for which the component of $\overline{Q_s-X}$ that contains the edge
lies in $W_t$, and \emph{down} edges, for which it lies in $V_t$. At each
disk of $E$, the up and down edges alternate as one moves around $\partial
E$ (see figure~\ref{fig:updown}). For each of the arcs of $\Gamma_1-E$, the
flattening process creates an up edge on one side and a down edge on the
other, but there is a fundamental difference in the way that the up and
down edges appear in $Q_s$ and $P_t$. As shown in figure~\ref{fig:updown},
up edges (the solid ones) and down edges (the dotted ones) alternate as one
moves around a crossover point, while in $Q_s$ they occur in adjacent
pairs. This is immediate upon examination of figure~\ref{fig:flatten}.

\begin{figure}
\includegraphics[width=0.95\textwidth]{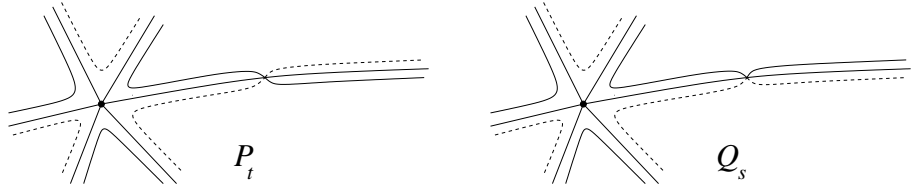}
\caption{Up and down edges of $X$ as they appear in $P_t$ and $Q_s$.}
\label{fig:updown}
\end{figure}

For our inductive procedure, we start with a pinched regular neighborhood
$X_1\subset Q_s\cap P_t$ of a graph $\Gamma_1$ in $Q_s\cap P_t$, all of
whose vertices have positive even valence.  Moreover, the edges of the
frontier of $X_1$ are up or down according to whether the portion of
$\overline{Q_s-X}$ that contains them lies in $W_t$ or $V_t$. We call this
an \emph{inductive configuration.}

To ensure that our isotopy process will terminate, we use the complexity
$-\chi(\Gamma_1)-\chi(\Fr(X_1))+N$, where $N$ is the number of components
of $\Gamma_1$.  Since all vertices of $\Gamma_1$ and $\Fr(X_1)$ have
valence at least~$2$, each of their components has nonpositive Euler
characteristic, so the complexity is a non-negative integer. The remaining
isotopies will reduce this complexity, so our procedure must terminate.

We may assume that the complexity is nonzero, since if $N=0$ then $X_1$ is
empty. Consider $X_1$ as a subset of the union of open disks
$Q_s-X_0$. Since $X_1$ is a regular neighborhood of a graph with vertices
of valence at least~$2$, it separates these disks, and we can find a closed
disk $D$ in $Q_s$ with $\partial D\subset X_1$ and $D\cap X=\partial D$.
It lies either in $V_t$ or $W_t$. Assume it is in $W_t$ (the case of $V_t$
is similar), in which case all of its edges are up edges. Since $\partial
D\subset P_t-X_0$, $\partial D$ bounds a disk $D'$ in $P_t-X_0$.  Since the
interior of $D$ is disjoint from $P_t$, $D\cup D'$ bounds a $3$-ball
$\Sigma$ in~$L$. Of course, $D'$ may contain portions of the component of
$X_1$ that contains $\partial D'$, or other components of~$X_1$. Let $X_1'$
be the component of $X_1$ that contains $\partial D'$; it is a pinched
regular neighborhood of a component $\Gamma_1'$ of~$\Gamma$.

Suppose that $X_1'$ contains some vertices of $\Gamma_1$ of valence more
than $2$. We will perform an isotopy of $Q_s$ that we call a \emph{tunnel
move,} illustrated in figure~\ref{fig:tunnel}, that reduces the complexity
of the inductive configuration. Near the vertex, select an arc in $X_1'$
that connects the edge of $\Fr(X_1')$ in $D'$ with another up edge of
$\Fr(X_1')$ that lies near the vertex (this arc may lie in $D'$, in a
portion of $X_1$ contained in $D'$). An isotopy of $Q_s$ is performed near
this arc, that pulls an open regular neighborhood of the arc in $X_1'$ into
$W_t$. This does not change the interior of $V_t-Q_s$ (it just adds the
regular neighborhood of the arc to $V_t-Q_s$), but in $W_t$ it creates a
tunnel that joins two different components of $W_t-Q_s$.  One of these
components was in $\Sigma$, so the isotopy cannot create core
circles. After the tunnel move, we have a new inductive configuration. The
Euler characteristic of $\Gamma_1$ has been increased by the addition of
one vertex, while $\chi(\Fr(X_1))$ and $N$ are unchanged, so the new
inductive configuration is of lower complexity. The procedure continues by
finding a new $D$ and $D'$ and repeating the process.
\begin{figure}
\includegraphics[width=60 ex]{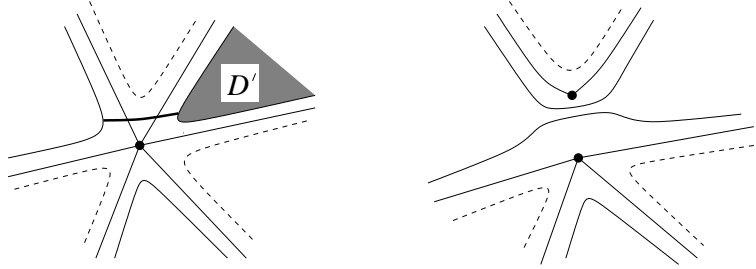}
\caption{A portion of $P_t$ showing a tunnel arc in $X_1$, and the new
$\Gamma_1$ and $X_1$ after the tunnel move.}
\label{fig:tunnel}
\end{figure}

When a $D$ has been found for which no tunnel moves are possible, all
vertices of $\Gamma_1'$ (if any) have valence $2$. Suppose that $X_1'$
contains crossover points. It must contain an even number of them, since up
and down edges alternate in $P_t$ around vertices of $\Gamma_1'$. Some
portion of $X_1'$ is a disk $B$ whose frontier consists of two crossover
points and two edges of $\Fr(X_1)$, each connecting the two crossover
points. There is an isotopy of $Q_s$, supported in a neighborhood of $B$,
that repositions $Q_s$ and replaces a neighborhood of $B$ in $X$ with a
rectangle containing no crossover points. Figure~\ref{fig:bigon}
illustrates this isotopy. It cannot create core circles, indeed the
interiors of $V_t-Q_s$ and $W_t-Q_s$ are unchanged during the isotopy. We
call such an isotopy a \emph{bigon move.}
\begin{figure}
\includegraphics[width=65 ex]{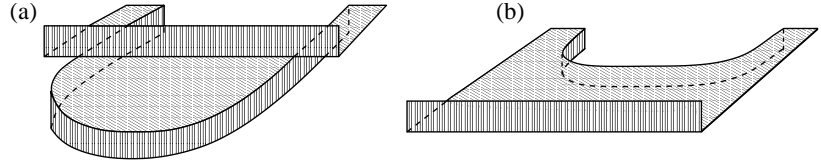}
\caption{Elimination of a bigon of $Q_s\cap P_t$ by isotopy.}
\label{fig:bigon}
\end{figure}

Since bigon moves increase the Euler characteristic of $\Fr(X_1)$, without
changing $\Gamma_1$ or $N$, they reduce complexity. So we eventually arrive
at the case when $X_1'$ is an annulus. Assume for now that the interior of
$D'$ is disjoint from $X_1$.  There is an isotopy of $Q_s$ that pushes $D$
across $\Sigma$, until it coincides with $D'$.  This cannot create core
circles, since its effect on the homeomorphism type of $W_t-Q_s$ is simply
to remove the component $\Sigma-Q_s$. Perform a small isotopy that pulls
$D'$ off into the interior of $W_t$, again creating no new core circles. An
annulus component of $X_1$ has been eliminated, reducing the complexity. If
the interior of $D'$ meets $X_1$, then $D'\cap X_1=X_1'$, and a similar
isotopy eliminates~$X_1'$.

Suppose now that the interior of $D'$ contains components of $X_1$ other
than perhaps $X_1'$.  Let $X_1''$ be their union. It is a pinched regular
neighborhood of a union $\Gamma_1''$ of components of $\Gamma_1$. If
$\Gamma_1''$ has vertices of valence more than $2$, then tunnel moves can
be performed. These cannot create new core circles, since they do not
change the interior of $V_t-Q_s$, and in $W_t-Q_s$ they only connect
regions that are contained in $\Sigma$. If no tunnel move is possible, but
there are crossover points, then a bigon move may be performed. So we may
assume that every component of $X_1''$ is an annulus.

Let $S$ be a boundary circle of $X_1''$ innermost on $D'$, bounding a disk
$D''$ in $D'$ whose interior is disjoint from $X$. Let $E''$ be the disk in
$Q_s$ bounded by $S$, so that $D''\cup E''$ bounds a $3$-ball $\Sigma''$
in~$L$.

We claim that if $(V_t -Q_s)\cup E''$ contains a core circle of $V_t$, then
$V_t-Q_s$ contained a core circle of $V_t$ (and analogously for $W_t$).
The closures of the components of $E''- X_1$ are planar surfaces, each
lying either in $V_t$ or $W_t$. Let $F$ be one of these, lying (say) in
$V_t$. Its boundary circles lie in $P_t-X_0$, so bound disks in $P_t$. A
regular neighborhood in $V_t$ of the union of $F$ and these disks is a
punctured $3$-cell $Z(F)$ meeting $P_t$ in disks.  Suppose that $C$ is a
core circle in $V_t$ that is disjoint from $Q_s-F$.  We may assume that $C$
meets $\partial Z(F)$ transversely, so cuts through $Z(F)$ is a collection
of arcs. Since $Z(F)$ is a punctured $3$-cell, there is an isotopy of $C$
that pushes the arcs to the frontier of $Z(F)$ and across it, removing the
intersections of $C$ with $F$ without creating new intersections (since the
arcs need only be pushed slightly outside of $Z(F)$). Performing such
isotopies for all components of $E''-X_1$ in $V_t$ produces a core circle
disjoint from $E''$, proving the claim.

By virtue of the claim, an isotopy that pushes $E''$ across $\Sigma''$
until it coincides with $D''$ does not create core circles. Then, a slight
additional isotopy pulls $D''$ and the component of $X_1$ that contained
$\partial D''$ off of $P_t$, reducing the complexity.

Since we can always reduce a nonzero complexity by one of these isotopies,
we may assume that $Q_s\cap P_t=X_0$. The frontier $\Fr(X_0)$ in $P_t$ is
the union of a graph $\Gamma_2$, each of whose components has vertices of
valence~$4$ corresponding to crossover points, and a graph $\Gamma_3$ whose
components are circles.

A component of $\Gamma_3$ must bound both a disk $D_Q$ in
$\overline{Q_s-X_0}$ and a disk $D_P$ in $\overline{P_t-X_0}$. Since
$Q_s\cap P_t=X_0$, the interiors of $D_P$ and $D_Q$ are disjoint, and $D_Q$
lies either in $V_t$ or in $W_t$. So we may push $D_Q$ across the $3$-ball
bounded by $D_Q\cup D_P$ and onto~$D_P$, without creating core circles.
Repeating this procedure to eliminate the other components of $\Gamma_3$,
we achieve that the frontier of $Q_s\cap P_t$ equals the graph $\Gamma_2$.

Figure~\ref{fig:meridian disks} shows a possible intersection of $Q_s$ with
$P_t$ at this stage. The shaded region is $Q_s\cap P_t$; the portion of it
that lies between $C_1$ and $C_2$ is a single octagon that passes around
the back of the torus. The closure of $Q_s-(Q_s\cap P_t)$ consists of two
meridian disks in $V_t$, bounded by the circles $C_1$ and $C_2$, and two
boundary-parallel disks in $W_t$, bounded by the circles $C_3$ and $C_4$.
\begin{figure}
\includegraphics[width=42 ex]{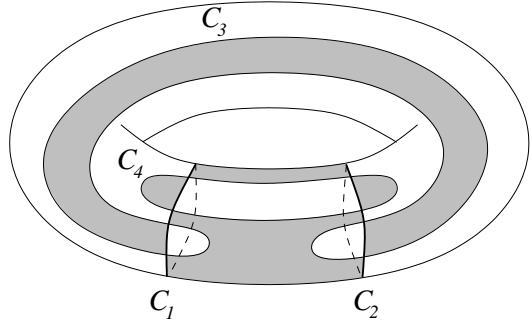}
\caption{A flattened torus containing two meridian disks.}
\label{fig:meridian disks}
\end{figure}

Suppose that there are now $2k_1$ meridian disks of $Q_s$ in $V_t$ and
$2k_2$ in $W_t$ (their numbers must be even since $Q_s$ is zero in
$H_2(L)$), and a total of $k_0$ boundary-parallel disks in $Q_s\cap V_t$
and $Q_s\cap W_t$. Since $\chi(Q_s)=0$, we have $\chi(Q_s\cap
P_t)=-k_0-2k_1-2k_2$. To prove the lemma, we must show that either $k_1$ or
$k_2$ is $0$.

Let $V$ be the number of vertices of $\Gamma_2$. Since all of its vertices
have valence $4$, $\Gamma_2$ has $2V$ edges. The remainder of $Q_s\cap P_t$
consists of $2$-dimensional faces. Each of these faces has boundary
consisting of an even number of edges, since up and down edges alternate
around a face.  If some of the faces are bigons, such as two of the faces
in figure~\ref{fig:meridian disks}, they may be eliminated by bigon
moves. These may create additional components of the frontier of $X_0$ that
are circles, indeed this happens in the example of figure~\ref{fig:meridian
disks}. These are eliminated as before by moving disks of $Q_s$ onto disks
in $P_t$. After all bigons have been eliminated, each face has at least
four edges, so there are at most $V/2$ faces. So we have $\chi(Q_s\cap
P_t)\leq V-2V+V/2=-V/2$.

Each boundary-parallel disk in $Q_s\cap V_t$ or $Q_s\cap W_t$ contributes
at least two vertices to the graph, since at each crossover point, $X_0$
crosses over to the other side in $P_t$ of the boundary of the disk. This
gives at least $2k_0$ vertices. The meridian disks on the two sides
contribute at least $2k_1\cdot 2k_2\cdot m$ additional vertices, where
$L=L(m,q)$, since the meridians of $V_t$ and $W_t$ have algebraic
intersection $\pm m$ in $P_t$. Thus $V\geq 2k_0+4k_1k_2m$. We calculate:
\[-k_0-2k_1-2k_2=\chi(Q_s\cap P_t)\leq -V/2\leq -k_0 -2k_1k_2m\ .\]
Since $m>2$, this can hold only when either $k_1$ or $k_2$ is $0$.
\end{proof}

Lemma~\ref{lem:common spine} fails (at the last sentence of the proof) for
the case of $L(2,1)$. Indeed, there is a flattened Heegaard torus in
$L(2,1)$ which meets $P_{1/2}$ in four squares and has two meridian disks
on each side. In a sketch somewhat like that of figure~\ref{fig:meridian
disks}, the boundaries of these disks are two meridian circles and two
$(2,1)$-loops intersecting in a total of $8$ points, and cutting the torus
into $8$ squares. There are two choices of four of these squares to
form~$Q_s\cap P_t$.

Now, we will complete the proof of theorem~\ref{thm:finding good
levels}. As in section~\ref{sec:Rubinstein-Scharlemann}, assume for
contradiction that all regions are labeled, and triangulate $I^2_\epsilon$.
The map on the $1$-skeleton is defined exactly as in
section~\ref{sec:Rubinstein-Scharlemann}, using
lemma~\ref{lem:borderlabelconstraints} and the fact that the labels satisfy
property~(RS\ref{item:edge labeling}). Using
lemma~\ref{lem:borderlabelconstraints}, each $1$-cell maps either to a
$0$-simplex or a $1$-simplex of the Diagram, and exactly as before the
boundary circle of $K$ maps to the Diagram in an essential way. The
contradiction will be achieved once we show that the map extends over the
$2$-cells.

There is no change from before when the $2$-cell meets $\partial K$ or lies
in the interior of $K$ but does not contain a vertex of $\Gamma$, so we fix
a $2$-cell in the interior of $K$ that is dual to a vertex $v_0$ of
$\Gamma$, located at a point~$(s_0,t_0)$.

Suppose first that $Q_{s_0}\cap P_{t_0}$ contains a spine of $P_{t_0}$. By
lemma~\ref{lem:common spine}, either $V_{t_0}$ or $W_{t_0}$ has a core
circle $C$ which is disjoint from $Q_{s_0}$; we assume it lies in
$V_{t_0}$, with the case when it lies in $W_{t_0}$ being similar.  The
letter $A$ cannot appear in the label of any region whose closure contains
$v_0$, since $C$ is a core circle for all $P_t$ with $t$ near $t_0$, and
$Q_s$ is disjoint from $C$ for all $s$ near $s_0$.  By lemma~\ref{lem:no
lone small letters}, any letter $a$ that appears in the label of one of the
regions whose closure contains $v_0$ must appear in a combination of either
$ax$ or $ay$, so none of these regions has label $\sA$. Since each $1$-cell
maps to a $0$- or $1$-simplex of the Diagram, the map defined on the
$1$-cells of $K$ maps the boundary of the $2$-cell dual to $v_0$ into the
complement of the vertex $\sA$ of the Diagram. Since this complement is
contractible, the map can be extended over the $2$-cell.

Suppose now that $Q_{s_0}\cap P_{t_0}$ does not contain a spine of
$P_{t_0}$. Then there is a loop $C_{(s_0,t_0)}$ essential in $P_{t_0}$ and
disjoint from $Q_{s_0}$. For every $(s,t)$ near $(s_0,t_0)$, there is a
loop $C_{(s,t)}$ essential in $P_t$ and disjoint from $Q_s$, with the
property that $C_{(s,t)}$ is a meridian of $V_t$ (respectively $W_t$) if and
only if $C_{(s_0,t_0)}$ is a meridian of $V_{t_0}$ (respectively $W_{t_0}$).
In particular, any intersection circle of $Q_s$ and $P_t$ which bounds a
disk in $Q_s$ which is precompressing for $P_t$ in $V_t$ or in $W_t$ must
be disjoint from $C_{(s,t)}$. Since the meridian disks of $V_t$ and $W_t$
have nonzero algebraic intersection, the meridians for $V_t$ and $W_t$
cannot both be disjoint from $C_{(s,t)}$. So for all $(s,t)$ in this
neighborhood of $(s_0,t_0)$, either all disks in $Q_s$ that are
precompressions for $P_t$ are precompressions in $V_t$, or all are
precompressions in $W_t$. In the first case, the letter $B$ does not appear
in the label of any of the regions whose closure contain $v_0$, while in
the second case, the letter $A$ does not. In either case, the extension to
the $2$-cell can now be obtained just as in the previous paragraph. This
completes the proof of theorem~\ref{thm:finding good levels}.

\newpage
\section[From good to very good]
{From good to very good}
\label{sec:from good to very good}

By virtue of theorem~\ref{thm:finding good levels}, we may perturb a
parameterized family of diffeomorphisms of $M$ so that at each parameter
$u$, some level $P_t$ and some image level $f_u(P_s)$ meet in good
position. In this section, we use the methodology of A. Hatcher \cite{H,
Hold} (see \cite{Hnew} for a more detailed version of \cite{Hold}, see also
N. Ivanov \cite{I3}) to change the family so that we may assume that $P_t$
and $f_u(P_s)$ meet in very good position. In fact, we will achieve a
rather stronger condition on discal intersections.

Following our usual notation, we fix a sweepout $\tau\colon P\times
[0,1]\to M$ of a closed orientable $3$-manifold $M$, and give $P_t$, $V_t$,
and $W_t$ their usual meanings. Given a parameterized family of
diffeomorphisms $f\colon M\times W\to M$, we give $f_u$, $Q_s$, $X_s$, and
$Y_s$ their usual parameter-dependent meanings. From now on, we refer to
the $P_t$ as \emph{levels} and the $Q_s$ as \emph{image levels.}

Throughout this section, we assume that for each $u\in W$, there is a pair
$(s,t)$ such that $Q_s$ and $P_t$ are in good position. Before stating the
main result, we will need to make some preliminary selections.

By transversality, being in good position is an open condition, so there
exist a finite covering of $W$ by open sets $U_i$, $1\leq i\leq n$, and
pairs $(s_i,t_i)$, so that for each $u\in U_i$, $Q_{s_i}$ and $P_{t_i}$
meet in good position. Note that we may reselect the $U_i$ to be open
$d$-balls whose closures are $d$-balls, by shrinking them and covering the
shrunken ones with finitely many such $d$-balls contained in the
original~$\overline{U(t_i)}$. Moreover, by shrinking of the open cover, we
can and always will assume that all transversality and good-position
conditions that hold at parameters in $U_i$ actually hold
on~$\overline{U_i}$.

We want to select the sets and parameters so that at parameters in $U_i$,
$Q_{s_i}$ is transverse to $P_{t_j}$ for all $t_j$. First note that for any
$s$ sufficiently close to $s_i$, $Q_s$ is transverse to $P_{t_i}$ at all
parameters of $U_i$ (here we are already using our condition that the
transversality for the $Q_{s_i}$ holds for all parameters in
$\overline{U_i}$). On $U_1$, $Q_{s_1}$ is already transverse to
$P_{t_1}$. Sard's Theorem ensures that at each $u\in U_2$, there is a value
$s$ arbitrarily close to $s_2$ such that $Q_s$ is transverse to $P_{t_1}$
at all parameters in a neighborhood of $u$. Replace $U_2$ by finitely many
open sets (with associated $s$-values), for which on each of these sets the
associated $Q_s$ are transverse to $P_{t_1}$. The new $s$ are selected
close enough to $s_2$ so that these $Q_s$ still meet $P_{t_2}$ in good
position. Repeat this process for $U_3$, that is, replace $U_3$ by a
collection of sets and associated values of $s$ for which the associated
$Q_s$ are transverse to $P_{t_1}$ and still meet $P_{t_3}$ in good
position.  Proceeding through the remaining original $U_i$, we have a new
collection, with many more sets $U_i$, but only the same $t_i$ values that
we started with, and at each parameter in one of the new $U_i$, $Q_{s_i}$
is transverse to $P_{t_1}$. Now proceed to $P_{t_2}$. For the $U_i$ whose
associated $t$-value is not $t_2$, we perform a similar process, and we
also select the new $s$-values so close to $s_i$ that the new $Q_s$ are
still transverse to $P_{t_1}$ and still meet their associated $P_{t_i}$ in
good position. After finitely many repetitions, all $Q_{s_i}$ are
transverse to each~$P_{t_j}$.

Since transversality is an open condition, we are free to replace $s_i$ by
a very nearby value, while still retaining the good position of $Q_{s_i}$
and $P_{t_i}$ and the transverse intersection of $Q_{s_i}$ with all
$P_{t_j}$, for all parameters in $U_i$, and similarly we may reselect any
$t_j$.  So (with the argument in the previous paragraph now completed) we
can and always will assume that all $s_i$ are distinct, and all $t_i$ are
distinct. We may then use $U(t_i)$ to denote the open set in $W$ associated
to the pair~$(s_i,t_i)$.
\shortpage

We can now state the main result of this section. With notation as above:
\begin{theorem}
Let $f\colon W\to \diff(M)$ be a parameterized family, such that for each
$u$ there exists $(s,t)$ such that $Q_s$ and $P_t$ meet in good
position. Then $f$ may be changed by homotopy so that there exists a
covering $\set{U(t_i)}$ as above, with the property that for all $u\in
U(t_i)$, $Q_{s_i}$ and $P_{t_i}$ meet in very good position, and $Q_{s_i}$
has no discal intersection with any $P_{t_j}$. If these conditions already
hold for all parameters in some closed subset $W_0$ of $W$, then the
deformation of $f$ may be taken to be constant on some neighborhood of
$W_0$.\par
\label{thm:from good to very good}
\end{theorem}

Before starting the proof, we introduce a simplifying convention. Although
strictly speaking, $Q_{s_i}$ is meaningful at every parameter, as is every
$Q_s$, throughout the remainder of this section we speak of $Q_{s_i}$ only
for parameters in $\overline{U(t_i)}$. That is, unless explicitly stated
otherwise, an assertion made about $Q_{s_i}$ means that the assertion holds
at parameters in $\overline{U(t_i)}$, but not necessarily at other
parameters. Also, to refer to $Q_{s_i}$ at a single parameter $u$, we use
the notation $Q_{s_i}(u)$. By our convention, $Q_{s_i}(u)$ is meaningful
only when $u$ is a value in~$\overline{U(t_i)}$.

Now, to preview some of the complications that appear in the proof of
theorem~\ref{thm:from good to very good}, consider the problem of removing,
just for a single parameter $u\in U(t_i)$, a discal component $c$ of the
intersection of $Q_{s_i}(u)$ with some $P_{t_j}$. Suppose that the disk
$D'$ in $Q_{s_i}(u)$ bounded by $c$ is innermost among all disks in
$Q_{s_i}(u)$ bounded by discal intersections of $Q_{s_i}(u)$ with the
$P_{t_k}$. Note that $D'$ can contain a nondiscal intersection of
$Q_{s_i}(u)$ with a $P_{t_k}$; such an intersection will be a meridian of
either $V_{t_k}$ or $W_{t_k}$ (although $k$ cannot equal $i$, since
$Q_{s_i}(u)$ and $P_{t_i}$ meet in good position). Let $D$ be the disk in
$P_{t_j}$ bounded by $c$, so that $D\cup D'$ is the boundary of a $3$-ball
$E$. There is an isotopy of $f_u$ that moves $D'$ across $E$ to $D$, and on
across $D$, eliminating $c$ and possibly other intersections of the
$Q_{s_\ell}(u)$ with the $P_{t_k}$.  We will refer to this as a \emph{basic
isotopy.}

It is possible for a basic isotopy to remove a biessential component of
some $Q_{s_k}(u)\cap P_{t_k}$. Examples are a bit complicated to describe,
but involve ideas similar to the construction in figure~\ref{fig:bad
torus}. Fortunately, the following lemma ensures that good position is not
lost.
\begin{lemma} 
After a basic isotopy as described above, each $Q_{s_k}(u)\cap P_{t_k}$ still
has a biessential component.
\label{lem:no biessential elimination}
\end{lemma}

\begin{proof}
Throughout the proof of the lemma, $Q_s$ is understood to mean~$Q_s(u)$.

Suppose that a biessential component of some $Q_{s_k}\cap P_{t_k}$ is
contained in the ball $E$, and hence is removed by the isotopy. Since a
spine of $Q_{s_k}$ cannot be contained in a $3$-ball, there must be a
circle of intersection of $Q_{s_k}$ with $D$ that is essential in
$Q_{s_k}$. This implies that $k\neq j$. Now $D'$ must have nonempty
intersection with $P_{t_k}$, since otherwise $P_{t_k}$ would be contained
in $E$. An intersection circle innermost on $D'$ cannot be inessential in
$P_{t_k}$, since $c$ was an innermost discal intersection on $Q_{s_i}$, so
$D'$ contains a meridian disk $D_0'$ for either $V_{t_k}$ or $W_{t_k}$.
Choose notation so that $D$ is contained in~$V_{t_k}$ (that is, $t_j<t_k$).

Suppose first that $D_0'\subset V_{t_k}$. The basic isotopy pushing $D'$
across $E$ moves $Q_{s_k}\cap E$ into a small neighborhood of $D$, so that
it is contained in $V_{t_k}$. If there is no longer any biessential
intersection of $Q_{s_k}$ with $P_{t_k}$, then the complement in $V_{t_k}$
of the original $D_0'$ contains a spine of $Q_{s_k}$ (since the original
intersection of $Q_{s_k}$ with $D$ contained a loop essential in $Q_{s_k}$,
the spine of $Q_{s_k}$ is now on the $V_{t_k}$-side of $P_{t_k}$). This is
a contradiction, since $Q_{s_k}$ is a Heegaard torus.

Suppose now that $D_0'\subset W_{t_k}$. Since the biessential circles of
$Q_{s_k}\cap P_{t_k}$ are disjoint from $D_0'$, they are meridians for
$W_{t_k}$ and hence are essential in $V_{t_k}$. Now, let $A$ be innermost
among the annuli on $Q_{s_k}$ bounded by a biessential component $C$ of
$Q_{s_k}\cap P_{t_k}$ and a circle of $Q_{s_k}\cap D$. Since $Q_{t_k}$ and
$P_{t_k}$ meet in good position, the intersection of the interior of $A$
with $P_{t_k}$ is discal. This implies that $C$ is contractible in
$V_{t_k}$, a contradiction.
\end{proof}

\begin{proof}[Proof of theorem~\ref{thm:from good to very good}]
We will adapt the approach of Hatcher \cite{H}. The principal difference
for us is that in \cite{H}, there is only a single domain level, whereas we
have the different $Q_{s_i}$ on the sets~$U(t_i)$.

The first step is to construct a family $h_{u,t}$, $0\leq t \leq 1$ of
isotopies of the $f_u=h_{u,0}$, which eliminates the discal intersections
of every $Q_{s_i}(u)$ with every $P_{t_j}$. Let $\mathcal{C}$ be the set of
all discal intersection curves of $Q_{s_i} \cap P_{t_j}$.  Since $Q_{s_i}$
is transverse to $P_{t_j}$ at all $u\in \overline{U(t_i)}$, the curves in
$\mathcal{C}$ fall into finitely many families which vary by isotopy as the
parameter moves over (the connected set) $\overline{U(t_i)}$. Thus we may
regard $\mathcal{C}$ as a disjoint union containing finitely many copies of
each $\overline{U(t_i)}$. It projects to $W$, with the preimage of $u$
consising of the discal intersection curves $\mathcal{C}_u$ of the
$Q_{s_i}(u)$ and $P_{t_j}$ for which $u\in U(t_i)$. By assumption, no
element of $\mathcal{C}$ projects to any parameter $u\in W_0$.

Each $c \in \mathcal{C}_u$ bounds unique disks $D_c \subset P_{t_j}$ and
$D'_c \subset Q_{s_i}(u)$ for some $i$ and $j$. The inclusion relations among
the $D_c$ define a partial ordering $<_P$ on $\mathcal{C}_u$, by the rule
that $c_1 <_P c_2$ when $D_{c_1} \subset D_{c_2}$. Similarly, $c_1 <_Q c_2$
when $D'_{c_1} \subset D'_{c_2}$.

If $c$ is minimal for $<_Q$, then $D'_c \cup D_c$ is an imbedded $2$-sphere
in $M$ which bounds a $3$-ball $E_c$.  By lemma~\ref{lem:no biessential
elimination}, the basic isotopy that pushes $D'_c$ across $E_c$ to $D_c$
and on to the other side of $D_c$ retains the property that every
$Q_{s_k}(u)\cap P_{s_k}$ has a biessential intersection. This ensures that
when all discal intersections have been eliminated, each $Q_{s_k}(u)\cap
P_{t_k}$ will still have an intersection, so they will be in very good
position.

Shrink the open cover $\set{U(t_i)}$ to an open cover $\set{U(t_i)'}$ for
which each $\overline{U(t_i)'}\subset U(t_i)$. To construct the $h_{u,t}$,
Hatcher introduced an auxiliary function $\Psi\colon \mathcal{C}\to (0,2)$
that gives the order in which the elements of $\mathcal{C}$ are to be
eliminated, and allows the basic isotopies to be tapered off as one nears
the frontier of $U(t_i)$. Denoting by $\psi_u$ the restriction of $\Psi$ to
$\mathcal{C}_u$, we will select $\Psi$ so that the following conditions are
satisfied:
\begin{enumerate}
\item $\psi_u (c) < \psi_u (c')$ whenever $c <_Q c'$
\item $\psi_u (c) < 1 $ if $c \subset Q_{s_i}(u)$ and $u \in \overline{U(t_i)'}$
\item $\psi_u (c) > 1$ if $c \subset Q_{s_i}(u)$ and $u \in
\overline{U(t_i)}-U(t_i)$.
\end{enumerate}
One way to construct such a $\Psi$ is to choose a Riemannian metric on
$\tau(P\times (0,1))$ for which each $P_t$ has area $1$, and define
$\Psi_0(c)$ to be the area of $f_u^{-1}(D_c')$ in $P_{s_i}$. Then,
choose continuous functions $\alpha_{t_i}$ which are $0$ on
$\overline{U(t_i)'}$ and $1$ on $W-U(t_i)$, and define
$\Psi(c)=\Psi_0(c)+\alpha_{t_i}(u)$ for $c\subset Q_{s_i}(u)$.

Roughly speaking, the idea of Hatcher's construction is to have $h_{u,t}$
perform the basic isotopy that eliminates $c$ during a small time interval
$I_u(c)$ which starts at the number $\psi_u(c)$. In order to retain control
of this process, preliminary steps must be taken to ensure that basic
isotopies that move points in intersecting $3$-balls $E_c$ do not occur at
the same time.

Define $G_0$ to be the subset of $W\times [0,2]$ consisting of all
$(u,\psi_u(c))$ with $c\in \mathcal{C}_u$. For a fixed isotopic family of
$c\in \mathcal{C}$ with $c \subset Q_{s_i}$, the points $(u,\psi_u(c))$
form a $d$-dimensional sheet $i(c)$ lying over $\overline{U(t_i)}$, where
$d$ is the dimension of $W$. If $i(c_1)$ meets $i(c_2)$, then by the first
property of $\Psi$, $c_1$ and $c_2$ cannot be $<_Q$-related.

Thicken each $i(c)$ to a plate $I(c)$ intersecting each $\set{u} \times
[0,2]$ in an interval $I_u(c)=[\psi_u (c), \psi_u (c) + \epsilon]$, for
some small positive $\epsilon$. This interval will contain the $t$-support of
the portion of $h_{u,t}$ that eliminates $c$, assuming that all other loops
in $\mathcal{C}_u$ with smaller $\psi_u$-value have already been
eliminated. By condition~(1), $c$ will be $<_Q$-minimal at the times $t\in
I_u(c)$. Since $\mathcal{C}_u$ is empty for $u\in W_0$, the $h_{u,t}$ will
be constant for all $u\in W_0$.

Choose the $\epsilon$ small enough so that $I(c_1) \cap I(c_2)$ is nonempty
only near the intersections of $i(c_1)$ and $i(c_2)$. This ensures that if
basic isotopies eliminating $c_1$ and $c_2$ occur on overlapping time
intervals, then $c_1$ and $c_2$ are $<_Q$-unrelated. Also, choose
$\epsilon$ small enough so that $I_u(c)\subset [0,1]$ whenever $u\in
U(t_i)'$.

It may happen that for some $c_1,c_2\in \mathcal {C}_u$ with $\psi_u (c_1) <
\psi_u (c_2)$, we have $c_2 <_P c_1$. In this case the isotopy which
eliminates $c_1$ will also eliminate $c_2$. So reduce $G_0$ by deleting all
points $(u,\psi_u (c_2))$ such that $\psi_u (c_1) < \psi_u (c_2)$ for some
$c_1$ with $c_2 <_P c_1$. Make a corresponding reduction of $I(c_2)$ by
deleting points $t \in I_u(c_2)$ such that $t > \psi_u (c_1)$ for some
$c_1$ with $c_2 <_P c_1$.

At values of $t$ where the interiors of $I(c_1)$ and $I(c_2)$ still
overlap, $c_1$ and $c_2$ are $<_Q$-unrelated, and the reduction just made
ensures that they are not $<_P$-related. In Hatcher's context, all
intersections are discal, so the combined effect of these is to eliminate
the possibility of simultaneous isotopies on intersecting $3$-balls
$E_{c_1}$ and $E_{c_2}$. In our context, however, $E_{c_1}$ and $E_{c_2}$
can intersect on overlaps of $I(c_1)$ and $I(c_2)$ even when $c_1$ and
$c_2$ are neither $<_P$-related nor $<_Q$-related. Figure~\ref{fig:nested}
shows a simple example. The intersections of $P_{t_1}$ with $Q_{s_2}$, are
not discal, nor are the intersections of $P_{t_2}$ with $Q_{s_1}$, but
$Q_{s_2}$ has a discal intersection with $P_{t_2}$ inside $E(c_1)$.
When this happens, however, $E_{c_1}$ and $E_{c_2}$
must be either disjoint or nested:\par
\begin{figure}
\includegraphics[width=70 ex]{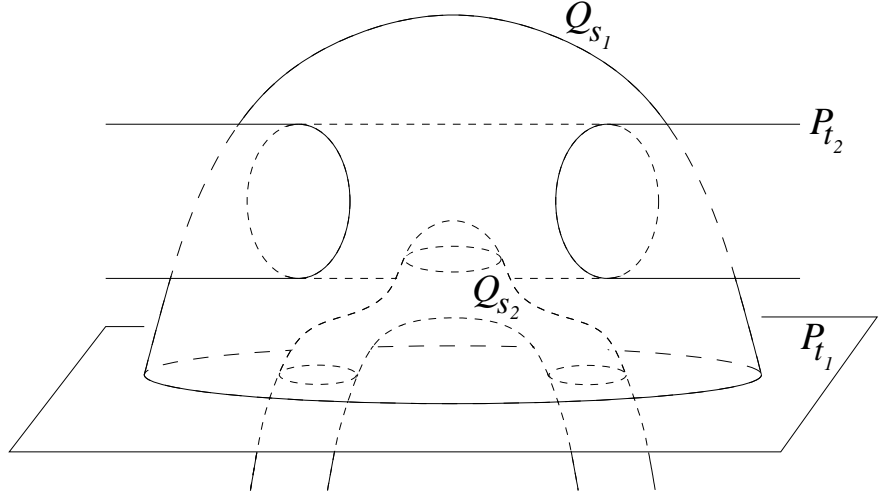}
\caption{Nested ball regions for basic isotopies.}
\label{fig:nested}
\end{figure}

\begin{lemma}
Suppose that $c_1$ and $c_2$ are $<_Q$-minimal discal intersections, and
are neither $<_P$-related nor $<_Q$-related. Then $\partial E_{c_1}$ and
$\partial E_{c_2}$ are disjoint.
\label{lem:nested E}
\end{lemma}

\begin{proof}
Since $c_1$ and $c_2$ are not $<_Q$-related, $D'(c_1)$ and $D'(c_2)$ are
disjoint, and since they are not $<_P$-related, $D(c_1)$ and $D(c_2)$ are
disjoint. An intersection circle of $D(c_1)$ and $D'(c_2)$ would be smaller
than $c_2$ in the $<_Q$-ordering, and similarly an intersection circle of
$D'(c_1)$ and $D(c_2)$ would be smaller than $c_1$ in the $<_Q$-ordering.
\end{proof}

When $E_{c_1}$ and $E_{c_2}$ are nested, say, $E_{c_2}$ lies in $E_{c_1}$,
a basic isotopy that removes $c_1$ will also remove $c_2$. So we make the
further reduction in $G_0$ of deleting all $(u,\psi_u (c_2))$ for which
there is a $c_1$ such that $i(c_1)$ meets $i(c_2)$, $\psi_u (c_1) < \psi_u
(c_2)$, and $E_{c_2}\subset E_{c_1}$. Also, reduce $I(c_2)$ by removing
any $t$ in $I_u(c_2)$ with $t>\psi_u(c_1)$.

For fixed $u \in W$, the basic isotopies are combined by proceeding upward
in $W\times [0,2]$ from $t = 0 $ to $t = 1$, performing each basic isotopy
involving $c$ on the interval $I_u(c)$. Condition~(3) on the $\psi_u$
ensures that the basic isotopies involving $c\subset Q_{s_i}(u)$ taper off
at parameters near the frontier of $U(t_i)$.  On a reduced interval
$I_u(c)$, which is an initial segment of $[\psi_u(c), \psi_u(c) +
\epsilon]$, peform only the corresponding initial portion of the basic
isotopy.  On the overlaps of the $I(c)$, perform the corresponding basic
isotopies concurrently; the reductions of the $I(c)$ have ensured that
these basic isotopies will have disjoint supports.  Since $\epsilon$ was
chosen small enough so that $I_u(c)\subset [0,1]$ whenever $u\in U(t_i)'$,
the basic isotopies involving $Q_{s_i}$ will be completed at all $u$ in
$U(t_i)'$. Since $\mathcal{C}_u$ is empty for $u\in W_0$, no isotopies take
place at parameters in $W_0$.

The remaining concern is that the basic isotopies eliminating $c\subset
Q_{s_i}(u)$ must be selected so that they fit together continuously in the
parameter $u$ on $U(t_i)$. This can be achieved using the method in the
last paragraph on p.~345 of \cite{H} (which applies in the smooth category
by virtue of \cite{H1}, see also the more detailed version in~\cite{Hnew}).
\end{proof}

\newpage
\section[Setting up the last step]
{Setting up the last step}
\label{sec:lemmas}

In this section, we present some technical lemmas that will be needed for
the final stage of the proof.

The first two lemmas give certain uniqueness properties for the fiber of
the Hopf fibration on $L$. Both are false for $\R\P^3$, so require our
convention that $L=L(m,q)$ with $m>2$, and as usual we select $q$ so that
$1\leq q\leq m/2$. From now on, we endow $L$ with the Hopf fibering and
assume that our sweepout of $L$ is selected so that each $P_t$ is a union
of fibers. Consequently the exceptional fibers, if any, will be components
of the singular set~$S$.

\begin{lemma}
Let $P$ be a Heegaard torus in $L$ which is a union of fibers, bounding
solid tori $V$ and $W$. Suppose that a loop in $P$ is a longitude for $V$
and for $W$. Then $q=1$ and the loop is isotopic in $P$ to a fiber.
\label{lem:bilongitude}
\end{lemma}

\begin{proof}
Let $a$ and $b$ be loops in $P$ which are respectively a longitude and a
meridian of $V$, and with $a$ determined by the condition that $ma+qb$ is a
meridian of $W$. Let $c$ be a loop in $P$ which is a longitude for both $V$
and $W$. Since $c$ is a longitude of $V$, it has (for one of its two
orientations) the form $a + kb$ in $H_1(P)$ for some $k$. The intersection
number of $c$ with $ma+qb$ is $q-km$, which must be $\pm1$ since $c$ is a
longitude of $W$. Since $1\leq q\leq m/2$ and $m>2$, this implies that
$k=0$ and $q=1$. Since $k=0$, $c$ is uniquely determined and $c=a$.  Since
$q=1$, the Hopf fibering is nonsingular, so the fiber is a longitude of
both $V$ and $W$ and hence is isotopic in $P$ to~$c$.
\end{proof}

\begin{lemma}
Let $h\colon L\to L$ be a diffeomorphism isotopic to the identity, with
$h(P_s)=P_t$. Then the image of a fiber of $P_s$ is isotopic in $P_t$ to a
fiber.\par
\label{lem:coincident levels}
\end{lemma}

\begin{proof}
Composing $f$ with a fiber-preserving diffeomorphism of $L$ that moves
$P_s$ to $P_t$, we may assume that $s=t$. Write $P$, $V$, and $W$ for
$P_t$, $V_t$, and $W_t$. Let $a$ and $b$ be loops in $P$ selected as in the
proof of lemma~\ref{lem:bilongitude}, and write $h_*\colon H_1(P)\to
H_1(P)$ for the induced isomorphism.

Suppose first that $h(V)=V$. Since the meridian disk of $V$ is unique up to
isotopy, we have $h_*(b)=\pm b$. Since $h$ is isotopic to the identity on
$L$ and $m>2$, $h$ is orientation-preserving and induces the identity on
$\pi_1(V)$. This implies that $h_*(b)=b$. Similar considerations for $W$
show that $h_*(ma+qb)=ma+qb$, so $h_*(a)=a$. Thus $h_*$ is the identity on
$H_1(P)$ and the lemma follows for this case.

Suppose now that $h(V)=W$. Then $h$ is orientation-reversing on $P$. Since
$h$ must take a meridian of $V$ to one of $W$, we have
$h_*(b)=\epsilon_1(ma+qb)$ where $\epsilon_1=\pm1$. Writing $h_*(a)=ua+vb$,
we find that $1=a\cdot b=-h_*(a)\cdot h_*(b)=-\epsilon_1(qu-mv)$. The facts
that $h$ is isotopic to the identity on $L$, $a$ generates $\pi_1(L)$, and
$b$ is $0$ in $\pi_1(V)$ imply that $u\equiv 1\pmod m$, so modulo $m$ we
have $1\equiv -\epsilon_1q$. Since $1\leq q\leq m/2$, this forces $q=1$,
$\epsilon_1=-1$, and $h_*(b)=-ma-b$.  Since $a$ has intersection number
$-1$ with the meridian $-ma-b$ of $W$, it is also a longitude of $W$.
Since $h$ is a homeomorphism interchanging $V$ and $W$, $h(a)$ is a
longitude of $V$ and of $W$, and an application of
lemma~\ref{lem:bilongitude} completes the proof.
\end{proof}

We now give several lemmas which allow the deformation of diffeomorphisms
and imbeddings to make them fiber-preserving or level-preserving. For
$Y\subset X$, $\imb(Y,X)$ means the connected component of the inclusion in
the space of all imbeddings of $Y$ in $X$. When $X$ is a fibered object,
$\Diff_f(X)$ means the space of diffeomorphisms of $X$ that take fibers to
fibers, and $\diff_f(X)$ is the connected component of the identity
~$\Diff_f(X)$.

\begin{lemma}
Let $X$ be either a solid torus or $S^1\times S^1\times I$, with a
fixed Seifert fibering.  Then the inclusion $\diff_f(X)\to \diff(X)$ is a
homotopy equivalence.\par
\label{lem:fiber-preserving on X}
\end{lemma}

\begin{proof}
Here is a very brief sketch of the proof; for detailed arguments of this
kind, see the final section of \cite{MR}. All needed results on fibrations
of spaces of diffeomorphisms appear in \cite{KM}.

Results from surface theory imply that $\diff(S^1\times S^1)\simeq
S^1\times S^1$. Using \cite{KM}, $\diff_f(X)$ is homotopy equivalent to
$S^1\times S^1$, and if $T$ is a boundary component, the restriction map
$\diff_f(X)\to \diff(T)$ is a weak homotopy equivalence. Using \cite{H},
one can show that the fiber of $\diff(X)\to \diff(T)$ is contractible, so
in $\diff_f(X)\to \diff(X)\to \diff(T)$, the composition and the second map
are homotopy equivalences, hence the first is as well.
\end{proof}

Lemma~\ref{lem:fiber-preserving on X} guarantees that if $g\colon\Delta\to
\diff(X)$ is a continuous map from an $n$-simplex, $n\geq 1$, with
$g(\partial \Delta)\subset \diff_f(X)$, then $g$ is homotopic relative to
$\partial \Delta$ to a map with image in $\diff_f(X)$. Analogous
observations hold for the next four lemmas as well.

When $X$ is fibered or Seifert-fibered and $Y\subset X$ is a union of
fibers, we write $\imb_f(Y,X)$ for the connected component of the inclusion
in the subspace of $\imb(Y,X)$ consisting of all imbeddings that take
fibers to fibers. The next two lemmas were proven in \cite{MR}, using
results on fibrations of spaces of mappings from \cite{KM}.

\begin{lemma} Let $T$ be a torus with a fixed
$S^1$-fibering, and let $C_n$ be a union of $n$ distinct fibers.  Then
$\imb_f(C_n,T)\to\imb(C_n,T)$ is a weak homotopy equivalence.
\label{lem:fiberpreserving1}
\end{lemma}

\begin{lemma}
Let $\Sigma$ be a compact $3$-manifold with nonempty boundary and having
a fixed Seifert fibering.  Let $F$ be a compact $2$-manifold properly
imbedded in $\Sigma$, such that $F$ is a union of fibers. Let
$\imb_{\partial f}(F,\Sigma)$ be the connected component of the
inclusion in the space of (proper) imbeddings for which the image of
$\partial F$ is a union of fibers. Then $\imb_f(F,\Sigma)\to
\imb_{\partial f}(F,\Sigma)$ is a weak homotopy equivalence.
\label{lem:fiberpreserving2}
\end{lemma}

The proof of the next lemma uses surface theory, somewhat along the lines
of the proof of theorem~\ref{thm:reduce to fiber-preserving}, and we do not
include the details.
\begin{lemma} Let $T$ be a torus with a fixed
$S^1$-fibering. Let $\Diff_h(T)$ be the subspace of $\Diff(T)$ consisting
of the diffeomorphisms that take some fiber to a loop isotopic to a fiber.
Then the inclusion $\Diff_f(T)\to \Diff_h(T)$ is a weak homotopy equivalence.
\label{lem:fiberpreserving0}
\end{lemma}

For $e\in (0,1)$ we let $eD^2$ denote the concentric disk of radius $e$ in
the standard disk $D^2\subset \R^2$. Let $X$ be either a solid torus
$D^2\times S^1$, or $T\times I$ where $T$ is a torus. Let $F=\cup F_i$ be a
disjoint union of finitely many tori. Fix an inclusion of $F$ into $X$ such
that each $F_i$ is of the form $\partial (e_iD^2\times S^1)$, in the solid
torus case, or of the form $T\times \set{e_i}$, in the $T^2\times I$ case,
for distinct numbers $e_i$ in $(0,1)$. Let $\imb_{\textit{int}}(F,X)$ be
the connected component of the inclusion in the space of all imbeddings of
$F$ into the interior of $X$, and let $\imb_{\textit{conc}}(F,X)$ be the
connected component of the inclusion in the set of imbeddings for which
each $F_i$ is of the form $\partial (eD^2)\times S^1$ or $T\times \set{e}$
for some $e\in (0,1)$. The next lemma is essentially the uniqueness of
collars of a boundary component.
\begin{lemma}
Let $X$ be a Seifert-fibered solid torus or $S^1\times S^1\times I$.  Then
the inclusion $\imb_{\text{conc}}(F,X)\to
\imb_{\text{int}}(F,X)$ is a homotopy equivalence.\par
\label{lem:level-preserving on F}
\end{lemma}
\noindent

\newpage
\section[Deforming to fiber-preserving families]
{Deforming to fiber-preserving families}
\label{sec:fiber-preserving families}
\longpage

\begin{theorem}
Let $L=L(m,q)$ with $m>2$ and let $f\colon S^d\to \diff(L)$.  Then
$f$ is homotopic to a map into $\diff_f(L)$.
\label{thm:make fiber-preserving}
\end{theorem}

\begin{proof}
Applying theorems~\ref{thm:generalposition}, \ref{thm:finding good levels},
and~\ref{thm:from good to very good}, we may assume that $f$ satisfies the
conclusion of theorem~\ref{thm:from good to very good}. That is, there are
pairs $(s_i,t_i)$ and an open cover $\set{U(t_i)}$ of $S^d$ with the
property that for every $u\in U(t_i)$, $Q_{s_i}(u)$ and $P_{t_i}$ meet in
very good position, and $Q_{s_i}(u)$ meets every $P_{t_j}$ transversely,
with no discal intersections. The $U(t_i)$ are selected to be connected, so
the intersection $Q_{s_i}(u)\cap P_{t_j}$ is independent, up to isotopy in
$P_{t_j}$, of the parameter~$u$. We remind the reader of our convention
that assertions about $Q_{s_i}$ implicitly mean ``for every $u\in
U(t_i)$.''  We can and always will assume that each $U(t_i)$ is connected,
and that conditions stated for parameters in $U(t_i)$ actually hold for all
parameters in~$\overline{U(t_i)}$.

Since the $t_j$ are distinct, we may select notation so that
$t_1<t_2<\cdots <t_m$. The corresponding $s_i$ typically are not in
ascending order. Figure~\ref{fig:si_chaos} shows a schematic picture of a
block of three levels for which the image levels $Q_{s_1}$, $Q_{s_2}$, and
$Q_{s_3}$ have $s_1<s_3<s_2$.
\begin{figure}
\includegraphics[width=30 ex]{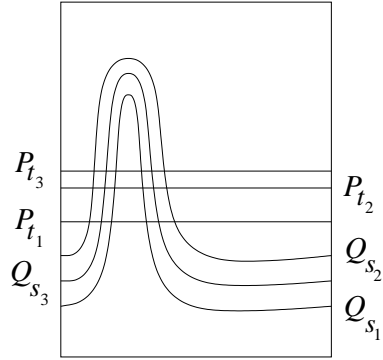}
\caption{A block of level tori with the $Q_{s_i}$ out of order.}
\label{fig:si_chaos}
\end{figure}

The basic idea of the proof is to make the $f_u$ fiber-preserving on the
$P_{s_i}$, then use lemma~\ref{lem:fiber-preserving on X} to make the $f_u$
fiber-preserving on the complementary $S^1\times S^1\times I$ or solid tori
of the $P_{s_i}$-levels. We must be very careful that none of the isotopic
adjustments to a $Q_{s_i}$ destroys any condition that must be preserved on
the other~$Q_{s_j}$.

Before listing the steps in the proof of theorem~\ref{thm:make
fiber-preserving}, a definition is needed. For each $i$, the intersection
circles of $Q_{s_i}\cap P_{t_i}$ cannot be meridians in both $V_{t_i}$ and
$W_{t_i}$, so $Q_{s_i}$ must satisfy exactly one of the following:
\begin{enumerate}
\item
The circles of $Q_{s_i}\cap P_{t_i}$ are not longitudes or meridians for
$V_{t_i}$, so the annuli of $Q_{s_i}\cap V_{t_i}$ are uniquely boundary
parallel in $V_{t_i}$.
\item
The circles of $Q_{s_i}\cap P_{t_i}$ are longitudes or meridians for
$V_{t_i}$, but are not longitudes or meridians for $W_{t_i}$, so the annuli
of $Q_{s_i}\cap W_{t_i}$ are uniquely boundary parallel in $W_{t_i}$.
\item
The circles of $Q_{s_i}\cap P_{t_i}$ are longitudes both for $V_{t_i}$ and
for $W_{t_i}$.
\end{enumerate}
In the first case, we say that $Q_{s_i}$ and $P_{t_i}$ are
\emph{$V$-cored,} in the second that they are \emph{$W$-cored,} and in the
third that they are \emph{bilongitudinal.} If they are either $V$-cored or
$W$-cored, we say they are \emph{cored.} Lemma~\ref{lem:bilongitude} shows
that the bilongitudinal case can occur only when $q=1$, and then only when
the intersection circles are isotopic in $P_{t_i}$ to fibers of the Hopf
fibering.

We can now list the steps in the procedure. In this list, and in the
ensuing details, ``push $Q_{s_i}$'' means perform a deformation of $f$ that
moves $Q_{s_i}$ as stated, and preserves all other conditions needed.
Making $Q_{s_i}$ ``vertical'' (at a parameter $u$) means making the
restriction of $f_u$ to $P_{s_i}$ fiber-preserving.  When we say that
something is done ``at all parameters of $U(t_i)$,'' we mean that a
deformation of $f$ will be performed, and that $U(t_i)$ is replaced by a
smaller set, so that the result is achieved for all parameters in the new
$\overline{U(t_i)}$, while retaining all other needed properties (such as
that $\set{U(t_i)}$ is an open covering of $S^d$).
\begin{enumerate}
\item[1.]  Push the $Q_{s_i}$ that meet $P_{t_j}$ out of $V_{t_j}$, for all
the $V$-cored $P_{t_j}$, at all parameters in $U(t_j)$. At the end of this
step, each $Q_{s_i}$ that was $V$-cored is parallel to $P_{t_i}$.
\vspace*{0.5 ex}
\item[2.]  Push the $Q_{s_i}$ that meet $P_{t_j}$ out of $W_{t_j}$, for all
the $W$-cored $P_{t_j}$, at all parameters in $U(t_j)$. At the end of this
step, each $Q_{s_i}$ that was $W$-cored is parallel to $P_{t_i}$.
\end{enumerate}
\noindent 
These first two steps are performed using a method of Hatcher like that of
the proof of section~\ref{sec:from good to very good}, although simpler.
After they are completed, a triangulation of $S^d$ is fixed with mesh
smaller than a Lebesgue number for the open cover by the $U(t_i)$. Each of
the remaining steps is performed by inductive procedures that move up the
skeleta of the triangulation, achieving the objective for $Q_{s_i}$ at all
parameters that lie in a simplex completely contained in $U(t_i)$.
\begin{enumerate}
\item[3.]  Push the $Q_{s_i}$ that originally were cored so that each one
equals some level torus. These level tori may vary from parameter to
parameter.\par
\vspace*{0.5 ex}
\item[4.]  
Push the $Q_{s_i}$ that originally were cored to be vertical.
\vspace*{0.5 ex}
\item[5.]
Push the bilongitudinal $Q_{s_i}$ to be vertical.
\vspace*{0.5 ex}
\item[6.]  
Use lemma~\ref{lem:fiber-preserving on X} to make $f_u$ fiber-preserving
on the complementary $S^1\times S^1\times I$ or solid tori of the
$P_{s_i}$-levels.
\end{enumerate}

The underlying fact that allows all of this pushing to be carried out
without undoing the results of the previous work is lemma~\ref{lem:hitting
levels}. Its use involves the concepts of \emph{compatibility} and
\emph{blocks,} which we will now define.

Recall that $R(t_i,t_j)$ means the closure of the region between $P_{t_i}$
and $P_{t_j}$. For a connected subset $Z$ of $S^d$, which in practice will
be either a single parameter or a simplex of a triangulation, denote by
$B_Z$ the set of $t_i$ such that $Z\subset U(t_i)$. Elements $t_i$ and
$t_j$ of $B_Z$ are called \emph{$Z$-compatible} when $Q_{s_i}(u)\cap
P_{t_i}$ and $Q_{s_k}(u)\cap P_{t_k}$ are homotopic in $R(t_i,t_k)$ for
every $t_k\in B_Z$ with $t_i<t_k\leq t_j$. Whether or not $t_i$ and $t_j$
are $u$-compatible typically varies as $u$ varies over $U(t_i)\cap U(t_j)$,
since $B_u$ depends on the other $U(t_k)$ that contain~$u$.

Because our family $f$ satisfies the conclusion of theorem~\ref{thm:from
good to very good}, lemma~\ref{lem:hitting levels} has the following
consequence: if $t_i$ and $t_j$ are $u$-compatible for any $u$,
then $P_{t_i}$ and $P_{t_j}$ are both $V$-cored, or both $W$-cored, or both
bilongitudinal. The next proposition is also immediate from
lemma~\ref{lem:hitting levels}.
\begin{proposition} Suppose that $t_i,t_j,t_k\in B_Z$. Then at parameters
in $Z$, $Q_{s_k}$ can meet both $P_{t_i}$ and $P_{t_j}$ only if $t_i$
and $t_j$ are $Z$-compatible.
\label{prop:hitting levels}
\end{proposition}

For a simplex $\Delta$, write $B_{\Delta}=\set{b_1,\ldots,b_m}$ with each
$b_i<b_{i+1}$, and for each $i\leq m$ define $a_i$ to be the $s_j$ for
which $b_i=t_j$. Decompose $B_{\Delta}$ into maximal $\Delta$-compatible
blocks $C_1=\set{b_1,b_2,\ldots,b_{\ell_1}}$,
$C_2=\set{b_{\ell_1+1},\allowbreak\ldots,\allowbreak b_{\ell_2}},\ldots\,$,
$C_r=\set{b_{\ell_{r-1}+1},\ldots,b_{\ell_r}}$, with $\ell_r=m$.  Since the
blocks are maximal, proposition~\ref{prop:hitting levels} shows that
$Q_{a_i}$ is disjoint from $P_{b_j}$ if $b_i$ and $b_j$ are not in the same
block. In steps~3-6, this disjointness will ensure that isotopies of these
$Q_{a_i}$ do not disturb the results of previous work.

Note that if $b_i$ and $b_j$ lie in the same block, then either both
$P_{b_i}$ and $P_{b_j}$ are $V$-cored, or both are $W$-cored, or both are
bilongitudinal. Thus we can speak of $V$-cored blocks, and so on.

When $\delta$ is a face of $\Delta$, $B_\Delta\subseteq
B_\delta$. Therefore if $b_i$ and $b_j$ in $B_\Delta$ are
$\delta$-compatible, then they are $\Delta$-compatible. So for each block
$C$ of $B_\delta$, $C\cap B_\Delta$ is contained in a block of $B_\Delta$.
However, levels that are not compatible in $B_\delta$ may become compatible
in $B_\Delta$, since the $t_i$ for intervening levels in $B_\delta$ may
fail to be in $B_\Delta$. Typically, the intersections of blocks of
$B_\delta$ with $B_\Delta$ will combine into larger blocks in~$B_\Delta$.

We should emphasize that the blocks of $B_Z$, and whether a level $P_{t_i}$
is $V$-cored, $W$-cored, or bilongitudinal, are defined with respect to the
original configuration, not the new positioning after the procedure
begins. Indeed, after steps~1 and~2, many of the $Q_{s_i}$ will be disjoint
from their~$P_{t_j}$.

We now fill in the details of these procedures.

\vspace*{0.5 ex}\noindent \textsl{Step 1: Push the $Q_{s_i}$ that meet
$P_{t_j}$ out of $V_{t_j}$, for all the $V$-cored $P_{t_j}$, at all
parameters in $U(t_j)$.}\vspace*{0.5 ex}

We perform this in order of increasing $t_j$ for the $V$-cored image
levels. Begin with $t_1$. If $Q_{s_1}$ is $W$-cored or bilongitudinal, do
nothing. Suppose it is $V$-cored. Then for each $u$ in $U(t_1)$, the
$Q_{s_j}(u)$ that meet $P_{t_1}$ intersect $V_{t_1}$ in a union of
incompressible uniquely boundary-parallel annuli.  Since any such $Q_{s_j}$
are transverse to $P_{t_1}$ at each point of $U(t_j)$, the set of
intersection annuli $Q_{s_j}\cap V_{t_1}$ falls into finitely many isotopic
families, with each family a copy of the connected set $U(t_j)$.  For each
$j$ with $U(t_1)\cap U(t_j)$ nonempty, let $\mathcal{A}_j$ be the
collection of the annuli $Q_{s_j}\cap V_{t_1}$, over all parameters in
$U(t_j)$, and let $\mathcal{A}$ be the union of these~$\mathcal{A}_j$. The
nonempty intersection of $U(t_1)$ and $U(t_j)$ ensures that the loops of
$Q_{s_j}\cap P_{t_1}$ and $Q_{s_1}\cap P_{t_1}$ are all in the same isotopy
class in $P_{t_1}$.

One might hope to push these families of annuli out of $V_{t_1}$ one at a
time, beginning with an outermost one, but an outermost family might not
exist. There could be a sequence $U(t_{j_1}),\ldots\,$, $U(t_{j_k})$ such
that $U(t_{j_i})\cap U(t_{j_{i+1}})$ is nonempty for each $i$,
$U(t_{j_k})\cap U(t_{j_1})$ is nonempty, and for some parameters $u_{j_i}$
in $U(t_{j_i})$, an annulus $Q_{s_{j_{i+1}}}(u_{j_i})\cap V_{t_1}$ lies
outside one of $Q_{s_{j_i}}(u_{j_i})\cap V_{t_1}$ for each $i$, and an
annulus of $Q_{s_{j_1}}(u_{j_k})\cap V_{t_1}$ lies outside one of
$Q_{s_{j_k}}(u_{j_k})\cap V_{t_1}$. Since an outermost family might not
exist, we will need to utilize the method of Hatcher as in the proof of
theorem~\ref{thm:from good to very good}, but only a simple version of it.

Shrink the $U(t_i)$ slightly, obtaining a new open cover by sets $U(t_i)'$
with $\overline{U(t_i)'}\subset U(t_i)$.
%Consider the $t_j$ for
%which $Q_{s_i}$ meets $V_{t_1}$ (at the parameters of $U(t_1)\cap U(t_j)$).
%Recalling that $Q_{s_j}$ is transverse to $P_{t_1}$ throughout $U(t_j)$,
%the annuli $\mathcal{A}_j$ of $Q_{s_j}\cap V_{t_1}$ fall into finitely many
%families as one varies over $U(t_j)$. 
We will use a function $\Psi\colon \mathcal{A}\to (0,2)$, so that at each
parameter $u$, the restriction $\psi_u$ of $\Psi$ to the annuli at that
parameter has the property that $\psi_u(A_1)<\psi_u(A_2)$ whenever $A_1,
A_2\in \mathcal{A}_i$ and $A_1$ lies in the region of parallelism between
$A_2$ and $\partial V_{t_1}$. Moreover, we will have $\psi_u(A)<1$ whenever
$A\in \mathcal{A}_i$ and $u\in \overline{U(t_i)'}$, while $\psi_u(A)>1$ for
$u$ near the boundary of $U(t_i)$. We construct $\Psi$ by letting
$\Psi_0(A)$ be the volume of the region of parallelism between $A$ and an
annulus in $\partial V_{t_1}$ (assuming that the volume of $L$ has been
normalized to $1$ to ensure that $\Psi_0(A)<1$), then adding on auxiliary
values $\alpha_{t_i}(u)$ as in the proof of theorem~\ref{thm:from good to
very good}.

Form the union $G_0\subset S^d\times (0,2)$ of the $(u,\psi_u(A))$ as in
the proof of theorem~\ref{thm:from good to very good}, and thicken each of
its sheets as was done there, obtaining an interval for each parameter.
These intervals tell the supports of the isotopies that push the annuli of
$Q_{s_j}\cap V_{t_1}$ out of $V_{t_1}$. If two sheets of $\mathcal{A}$
cross in $S^d\times (0,2)$, then the corresponding regions of parallelism
have the same volume, so must be disjoint and the isotopies can be
performed simultaneously without interference. At each individual
parameter $u$, each annulus is outermost during the time it is being
pushed out of $V_{t_1}$, but the times need to be different since there may
be no outermost family.

After the process is completed, $Q_{s_j}$ will lie outside of $V_{t_1}$ at
all parameters in $\overline{U(t_j)'}$, whenever $U(t_j)$ had nonempty
intersection with $U(t_1)$. Replacing each $U(t_j)$ by $U(t_j)'$, we have
$Q_{s_j}$ pushed out of $V_{t_1}$ at all parameters in these $U(t_j)$.
Moreover, lemma~\ref{lem:pushout}(2) shows that $V_{t_1}$ is concentric in
either $X_{s_1}$ or $Y_{s_1}$ at all parameters in~$U(t_1)$.

Some of the $Q_{s_k}$ for which $U(t_k)$ did not meet $U(t_1)$ may be moved
by the isotopies of the $Q_{s_j}$ at parameters in $U(t_j)\cap U(t_k)$. The
condition that these $Q_{s_k}$ meet $P_{t_1}$ transversely may be lost, but
this will not matter, because these intersections never matter when
$U(t_k)$ does not meet~$U(t_1)$.

Now consider $t_2$. Again, we do nothing if $Q_{s_2}$ is $W$-cored or
bilongitudinal, so suppose that it is $V$-cored. Use the Hatcher process as
before, to push annuli in the $Q_{s_j}$ out of $V_{t_2}$, when $Q_{s_j}$
meets $P_{t_2}$ and $U(t_j)$ meets $U(t_2)$.  Notice that these $Q_{s_j}$
cannot meet $V_{t_1}$ at parameters in $U(t_1)$.  For if $t_2$ is not
$u$-compatible with $t_1$ at some parameters in $U(t_1)$, then (by
lemma~\ref{lem:hitting levels}) $Q_{s_j}$ cannot meet both $P_{t_2}$ and
$P_{t_1}$, while if it is $u$-compatible at some parameter in $U(t_1)$,
then it has already been pushed out of $V_{t_1}$.  And $V_{t_1}$ cannot lie
in any of the regions of parallelism for the pushouts from $V_{t_2}$, since
the intersection circles of the $Q_{s_j}$ with $P_{t_2}$ are not longitudes
in $V_{t_2}$.

After these pushouts are completed, if $i=1$ or $i=2$ and $Q_{s_i}$ was
$V$-cored, then $V_{t_i}$ is concentric in either $X_{s_i}$ or $Y_{s_i}$ at
all parameters in~$U(t_i)$.

We continue working up the increasing $t_i$ in this way.  At the end of
this process, $V_{t_i}$ is concentric in either $X_{s_i}$ or $Y_{s_i}$ for
all $i$ such that $Q_{s_i}$ was $V$-cored, and at all parameters in
$U(t_i)$. For $Q_{s_i}$ that were $W$-cored or bilongitudinal, the
intersections $Q_{s_i}\cap P_{t_i}$ have not been disturbed at parameters
in $U(t_i)$. We have not introduced any new intersections of $Q_{s_i}$ with
$P_{t_j}$, so we still have the property that at any parameter $u$ in
$U(t_i)\cap U(t_j)$, $Q_{s_j}$ can meet $P_{t_i}$ only if $t_i$ and $t_j$
were originally $u$-compatible.

\vspace*{0.5 ex}\noindent \textsl{Step 2: Push the $Q_{s_i}$ that meet
$P_{t_j}$ out of $W_{t_j}$ for all the $Q_{s_j}$ that are $W$-cored, at all
parameters in $U(t_j)$.}\vspace*{0.5 ex}

The entire process is repeated with $W$-cored levels, except that we start
with $t_m$ and proceed in order of decreasing $t_i$. Each $W$-cored
$Q_{s_i}$ is pushed out of $W_{t_i}$, and at the end of the process
$W_{t_i}$ is concentric in either $X_{s_i}$ or $Y_{s_i}$ at all parameters
in $U(t_i)$, whenever $Q_{s_i}$ was $W$-cored. No intersection of a
$Q_{s_j}$ with a $V$-cored or bilongitudinal level $P_{t_i}$ is changed at
any parameter in~$U(t_i)$.

For the remaining steps, we fix a triangulation of $S^d$ with mesh smaller
than a Lebesgue number for $\set{U(t_i)}$, which will ensure that
$B_\Delta$ is nonempty for every simplex~$\Delta$. We will no longer
proceed up or down all $t_i$-levels, working on the sets $U(t_i)$, but
instead will work inductively up the skeleta of the triangulation. Recall
that each $B_\Delta$ is decomposed into blocks, according to the original
intersections of the $Q_{s_i}$ and $P_{t_i}$ before steps~1 and~2 were
performed.

\vspace*{0.5 ex}\noindent \textsl{Step 3: Push the $Q_{s_i}$ that were
originally cored so that each one equals some level torus.}\vspace*{0.5 ex}

We will proceed inductively up the skeleta of the triangulation, moving
cored $Q_{s_i}$ to level tori, without changing $Q_{s_k}\cap P_{s_k}$ for
the bilongitudinal $Q_{s_k}$. We want to use the fact that $V_{t_i}$ (or
$W_{t_i}$) is concentric with $X_{s_i}$ or $Y_{s_i}$ to push $Q_{s_i}$ onto
a level torus, but when moving multiple levels at a given parameter, there
is a consistency condition needed. As shown in figure~\ref{fig:inconsistent
levels}, it might happen that $V_{t_i}$ is concentric in $X_{s_i}$ while
$V_{t_j}$ is concentric in $Y_{s_j}$. Then, we might not be able to push
$Q_{s_i}$ and $Q_{s_j}$ onto level tori without disrupting other
levels. The following lemma rules out this bad configuration.
\begin{figure}
\includegraphics[width=35 ex]{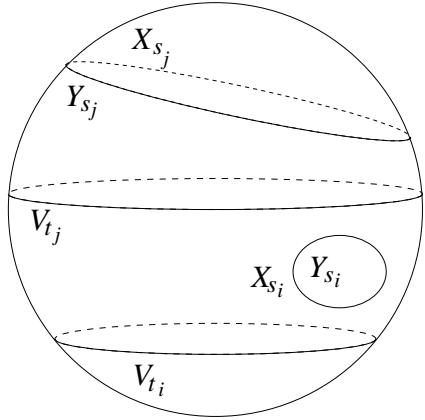}
\caption{Hypothetical inconsistent nesting: $V_{t_i}\subset X_{s_i}$ and
$V_{t_j}\subset Y_{s_j}$.}
\label{fig:inconsistent levels}
\end{figure}

\begin{lemma} Suppose, after steps~1 and~2 have been completed,
that $u\in U(t_i)\cap U(t_j)$, $t_i<t_j$, and that
$Q_{s_i}$ is $V$-cored.
\begin{enumerate}
\item The region between $Q_{s_i}$ and $Q_{s_j}$ does not contain a core
circle of $V_{t_i}$.
\item Suppose that $t_i$ and $t_j$ are $u$-compatible, and $V_{t_i}$ is
concentric in $Z_{s_i}$ where $Z$ is $X$ or $Z$ is $Y$. Then $V_{t_j}$ is
concentric in $Z_{s_j}$.
\item If $t_i$ and $t_j$ are $u$-incompatible, then $Q_{s_i}$ is parallel
to $P_{t_i}$ in $R(t_i,t_j)$.\par
\end{enumerate}\par
\noindent The analogous statement holds when $Q_{s_j}$ is $W$-cored
and $W_{t_j}$ is concentric in $Z_{s_j}$.
\label{lem:nesting}
\end{lemma}
\noindent 

\begin{proof}
It suffices to consider the case when $Q_{s_i}$ is $V$-cored. In the
situation at the start of Step~1 above, when annuli in the $Q_{s_k}$
were being pushed out of $V_{t_i}$, the intersection of $Q_{s_i}\cup
Q_{s_j}$ with $V_{t_i}$ was a union $F$ of incompressible nonlongitudinal
annuli. Since $Q_{s_i}$ met $P_{t_i}$, $F$ was nonempty. By
proposition~\ref{prop:core circles}, exactly one complementary region of
$F$ in $V_{t_i}$ contained a core circle $C$ of $V_{t_i}$. For at least one
of $s_i$ and $s_j$, say for $s_k$, $Q_{s_k}$ met this complementary region.

Since the annuli of $F$ are nonlongitudinal, there is an imbedded circle
$C'$ in $Q_{s_k}$ that is homotopic in the core region to a proper multiple
of $C$. If $C$ were in the region $R=f_u(R(s_i,s_j))$ between $Q_{s_i}$ and
$Q_{s_j}$, then the imbedded circle $C'$ in $\partial R$ would be a proper
multiple in $\pi_1(R)$, which is impossible since $R$ is homeomorphic to
$S^1\times S^1\times I$. This proves~(1).

Assume that $t_i$ and $t_j$ are $u$-compatible and suppose that
$V_{t_i}\subset X_{s_i}$ and $V_{t_j}\subset Y_{s_j}$.  Then $C$ is
contained in $X_{s_i}\cap Y_{s_j}$, forcing $s_i>s_j$ and $C$ in the region
between $Q_{s_i}$ and $Q_{s_j}$, contradicting (1). The case of
$V_{t_i}\subset Y_{s_i}$ and $V_{t_j}\subset X_{s_j}$ is similar, so (2)
holds.

For (3), if $t_i$ and $t_j$ are not $u$-compatible, then $Q_{s_j}$ was
disjoint from $P_{t_i}$ before the pushouts, so $i=k$ and $Q_{s_i}$ is
disjoint from $P_{t_j}$. If $Q_{s_i}$ is not parallel to $P_{t_i}$ in
$R(t_i,t_j)$, then $Q_{s_i}$ does not separate $P_{t_j}$ from
$P_{t_i}$. This implies that before pushouts from $V_{t_i}$ it did not
separate $Q_{s_j}\cap P_{t_j}$ from $P_{t_i}$, so after pushouts it could
not separate $Q_{s_j}$ from $P_{t_i}$, contradicting~(1).
\end{proof}

It will be convenient to extend our previous notation $R(s,t)$ for the
closure of the region between $P_s$ and $P_t$, by putting $R(0,t)=V_t$,
$R(t,1)=W_t$, and $R(0,1)=L$.

We will now define target regions. The isotopies that we will use in the
rest of our process will only change values within a single target region,
ensuring that the necessary positioning of the $Q_{s_i}$ is retained.  Let
$\Delta$ be a simplex of the triangulation, and recall the decomposition of
$B_\Delta=\set{b_1,\ldots,b_m}$ into maximal $\Delta$-compatible blocks
$C_1=\set{b_1,b_2,\ldots,b_{\ell_1}}$,
$C_2=\set{b_{\ell_1+1},\allowbreak\ldots,\allowbreak b_{\ell_2}},\ldots\,$,
$C_r=\set{b_{\ell_{r-1}+1},\ldots,b_{\ell_r}}$. Define the \emph{target
region} of a block $C_n$ to be the submanifold $T_\Delta(C_n)$ of $L$
defined as follows. Put $\ell_0=0$, $b_0=0$, and $b_{\ell_r+1}=1$.
\begin{enumerate}
\item If $C_n$ is $V$-cored, then 
$T_\Delta(C_n)=R(b_{\ell_{n-1}+1},b_{\ell_n+1})$.
\item If $C_n$ is $W$-cored, then 
$T_\Delta(C_n)=R(b_{\ell_{n-1}},b_{\ell_n})$.
\item If $C_n$ is bilongitudinal, then
$T_\Delta(C_n)=R(b_{\ell_{n-1}},b_{\ell_n+1})$.
\end{enumerate}
We remark that $T_\Delta(C_n)$ is all of $L$ when $B_\Delta$ consists of a
single bilongitudinal block, otherwise is of the form $V_t$ when $n=1$ and
$C_1$ is $W$-cored or bilongitudinal and of the form $W_t$ when $n=r$ and
$C_n$ is $V$-cored or bilongitudinal, and in all other cases it is a region
$R(s,t)$ diffeomorphic to $S^1\times S^1\times I$.

As noted in the next lemma, the interior of the target region of a block
contains the $Q_{a_i}$ for the $b_i$ in the block, at this point of our
argument.
\begin{lemma} Target regions satisfy the following.
\begin{enumerate}
\item If $b_i\in C_n$ and $u\in \Delta$, then $Q_{a_i}(u)$ is in the
interior of $T_\Delta(C_n)$.
\item If $\delta$ is a face of $\Delta$, and $C'_1,\ldots\,$, $C'_{r'}$ are
the blocks of $B_\delta$, then for each $i$, there exists a $j$ such that
$T_\delta(C'_i)\subseteq T_\Delta(C_j)$.
\end{enumerate}
\label{lem:target}
\end{lemma}

\begin{proof}
Property (1) is a consequence of proposition~\ref{prop:hitting levels} and
the fact that Steps~1 and~2 do not create new intersections of the
$Q_{s_i}(u)$ with the $P_{t_j}$. For part~(2), the proof is direct
from the definitions, dividing into various subcases.
\end{proof}

Target regions can overlap in the following ways: the target region for a
$V$-cored block $C_n$ will overlap the target region of a succeeding
$W$-cored block $C_{n+1}$, and the target region of a bilongitudinal block
will overlap the target region of a preceeding $V$-cored block or of a
succeeding $W$-cored block (note that by lemma~\ref{lem:bilongitude},
successive blocks cannot both be bilongitudinal). The latter cause no
difficulties, but the conjunctions of a $V$-cored block and a succeeding
$W$-cored block will necessitate some care during the ensuing argument.

We can now begin the process that will complete Step~3. We will start at
the parameters that are vertices of the triangulation and move the
$Q_{a_i}$ for each $V$-cored or $W$-cored block to be level, that is, so
that each $Q_{a_i}(u)$ equals some $P_t$. The isotopies will be fixed on
each $P_{b_i}$ for which $Q_{a_i}$ is bilongitudinal, and these unchanged
$Q_{a_i}\cap P_{b_i}$ will be used to work with the bilongitudinal levels
in a later step.  For each cored block, the isotopy that levels the
$Q_{a_i}$ will move points only in the interior of the target region of the
block. As we move to higher-dimensional simplices, the $Q_{a_i}$ will
already be level at parameters on the boundary, and the deformation will be
fixed at those parameters. Each deformation for the parameters in a simplex
$\delta_0$ of dimension less than $d$ must be extended to a deformation of
$f$. The extension will change an $f_u$ only when $u$ is in the open star
of $\delta_0$, by a deformation that performs some initial portion of the
deformation of $f_{u_0}$ at a parameter $u_0$ of $\delta_0$--- the
parameter that is the $\delta_0$-coordinate of $u$ when the simplex that
contains it is written as a join $\delta_0*\delta_1$. We will see that
because the target regions can overlap, the deformation of an $f_u$ might
not preserve all target regions, but enough positioning of the image
levels~$Q_{a_i}$ will be retained to continue the inductive process.

Fix a vertex $\delta_0$ of the triangulation, and consider the first block
$C_1$ of $B_{\delta_0}$. If it is bilongitudinal, we do nothing. Suppose
that it is $V$-cored.  All of the $Q_{a_1},\ldots\,$, $Q_{a_{\ell_1}}$ lie
in the interior of the target region
$T_{\delta_0}(C_1)$. Lemma~\ref{lem:nesting}(2) shows that for either $Z=X$
or $Z=Y$, $V_{b_i}$ is concentric in $Z_{a_i}$ for $b_i\in C_1$. We claim
that there is an isotopy, supported on $T_{\delta_0}(C_1)$, that moves each
$Q_{a_i}$ to be level. If $C_1$ is the only block, then
$T_{\delta_0}(C_1)=L$ and the isotopy exists by the definition of
concentric. If there is a second block, then lemma~\ref{lem:nesting}(3)
shows that the $Q_{a_i}$ for $b_i\in C_1$ are parallel to $P_{b_1}$ in
$T_{\delta_0}(C_1)=R(b_1,b_{\ell_1+1})$, and again the isotopy exists.
After performing the isotopy, we may assume that the $Q_{a_i}(u_0)$ are
level.

To extend this deformation of $f_{\delta_0}$ to a deformation of the
parameterized family $f$, we regard each simplex $\Delta$ of the closed
star of $\delta_0$ in the triangulation as the join $\delta_0*\delta_1$,
where $\delta_1$ is the face of $\Delta$ spanned by the vertices of
$\Delta$ other than $\delta_0$. Each point of $\Delta$ is uniquely of the
form $u=s\delta_0+(1-s)u_1$ with $u_1\in \delta_1$. Write the isotopy of
$f_{\delta_0}$ as $j_t\circ f_{\delta_0}$, with $j_0$ the identity map of
$L$. Then, at $u$ the isotopy at time $t$ is $j_t\circ f_u$ for $0\leq
t\leq s$ and $j_s\circ f_u$ for $s\leq t\leq 1$.  For any two simplices
containing $\delta_0$, this deformation agrees on their intersection, so it
defines a deformation of~$f$.

The target region $T_{\delta_0}(C_1)$ will overlap $T_{\delta_0}(C_2)$ if
$C_2$ is bilongitudinal or $W$-cored. When $C_2$ is bilongitudinal, this
does not affect any of our necessary positioning. If it is $W$-cored, then
$Q_{a_i}$ with $b_i\in C_2$ may be moved into $T_{\delta_0}(C_1)$. At
$\delta_0$, they will end up somewhere between the now-level
$Q_{a_{\ell_1}}$ and $P_{b_{\ell_2}}$, and at other parameters in the star
of $\delta_0$ they will lie somewhere in $R(b_1,b_{\ell_2})$. This will
require only a bit of attention in the later argument.

In case $C_1$ was $W$-cored, we use lemmas~\ref{lem:nesting}(2)
and~\ref{lem:level-preserving on F}, producing a deformation of
$f_{\delta_0}$ supported on the interior of the solid torus
$T_\Delta(C_1)=V_{b_{\ell_1}}$, which does not meet any other target
region. This is extended to a deformation of $f$ just as before.

We move on to consider $C_2$ in analogous fashion, doing nothing if $C_2$
is bilongitudinal, and moving the $Q_{a_i}$ to be level at the parameter
$\delta_0$. If $C_1$ was $V$-cored and $C_2$ is $W$-cored, then instead of
the initial target region $T_{\delta_0}(C_2)$ we must use the region
between the now-level $Q_{a_{\ell_1}}(u)$ and $P_{b_{\ell_2}}$, but
otherwise the argument is the same. Proceed in the same way through the
remaining blocks $C_n$ of~$B_{\delta_0}$, ending with all the cored
$Q_{a_i}(u_0)$ moved to be level. This process for $u_0$ is repeated for
each $0$-simplex of the triangulation.

Now, consider a simplex $\delta$ of positive dimension. Inductively, we may
assume that at each $u$ in $\partial \delta$, each cored $Q_{a_i}$ has been
moved to a level torus, and $Q_{a_i}\cap P_{a_i}$ is unchanged for each
bilongitudinal $Q_{a_i}$. Moreover, if $a_i$ is contained in a target block
$T_\delta(C_i)$, then $Q_{a_i}$ lies in its target region, or else lies in
the union of the target regions for a $V$-cored block and a succeeding
$W$-cored block. Note that we are using lemma~\ref{lem:nesting}(2) here.

We apply lemma~\ref{lem:level-preserving on F} to each cored block of
$B_\Delta$, sequentially up the cored blocks. We obtain a sequence of
deformations of $f$ on $\delta$, constant at parameters in $\partial
\delta$. There is no interference between different blocks, except when a
$W$-cored block $C_{n+1}$ succeeds a $V$-cored block $C_n$. First, the
$Q_{a_i}$ for the $V$-cored block are moved to be level. Then, at each
parameter in $\delta$, the $Q_{a_i}(u)$ for the $W$-cored block lie between
the now-level $Q_{a_{\ell_n}}(u)$ and $P_{b_{\ell_{n+1}}}$. We regard the
union of these regions over the parameters of $\delta$ as a product
$\delta\times S^1\times S^1\times I$, and apply
lemma~\ref{lem:level-preserving on F}. Thus the isotopy that levels the
$Q_{a_i}$ from the $W$-cored block need not move any of the $Q_{a_i}$ from
the $V$-cored block. In other cases, the successive isotopies take place in
disjoint regions. To extend this to a deformation of $f$, we adapt the join
method from above (of course when $\delta_0$ is $d$-dimensional, no
extension is necessary). Regard each simplex $\Delta$ of the closed star of
$\delta_0$ in the triangulation as the join $\delta_0*\delta_1$, where
$\delta_1$ is the face of $\Delta$ spanned by the vertices of $\Delta$ not
in $\delta_0$. Each point of $\Delta$ is uniquely of the form $u=s
u_0+(1-s)u_1$ with $u_i\in \delta_i$. Write the isotopy of $f_{u_0}$ as
$j_t\circ f_{u_0}$, with $j_0$ the identity map of $L$. Then, at $u$ the
isotopy at time $t$ is $j_t\circ f_u$ for $0\leq t\leq s$ and $j_s\circ
f_u$ for $s\leq t\leq 1$. For any two simplices containing $\delta_0$, this
deformation agrees on their intersection, so it defines a deformation
of~$f$.

At the completion of this process, each cored $Q_{s_i}$ is level at all
parameters in $\Delta$, whenever $\Delta\subset U(t_i)$. The bilongitudinal
$Q_{s_i}$ may have been moved around some, but their intersections
$Q_{s_i}\cap P_{t_i}$ will not be altered at parameters for which $t_i\in
B_\Delta$ since these intersections will not lie in the interior of any
target region for a cored level.

\vspace*{0.5 ex}\noindent \textsl{Step 4: Push all cored $Q_{s_i}$ to be
vertical, that is, make each image of a fiber of $P_{s_i}$ a fiber in $L$.
}\vspace*{0.5 ex}

Again we work our way up the simplices of the triangulation. Start at a
$0$-simplex $\delta_0$. Each cored $Q_{a_i}(\delta_0)$ for $b_i\in
B_{\delta_0}$ is now level.  By lemma~\ref{lem:coincident levels}, the
image fibers in $Q_{a_i}(\delta_0)$ are isotopic in that level torus to
fibers of $L$. Using lemma~\ref{lem:fiberpreserving0}, there is an isotopy
of $f_{\delta_0}$ that preserves the level tori and makes
$Q_{a_i}(\delta_0)$ vertical. This isotopy can be chosen to fix all points
in other $Q_{a_j}(\delta_0)$, and is extended to a deformation of $f$ by
using the method of Step~3. We work our way up the skeleta; if
$\delta\subset U(b_i)$, then for every $u$ in $\delta$, each $Q_{a_i}(u)$
is level torus, and at parameters $u\in\partial \delta$, $Q_{a_i}(u)$ is
vertical.  Using lemma~\ref{lem:fiberpreserving0}, we make the $Q_{a_i}(u)$
vertical at all $u\in \delta$, and extend to a deformation of $f$ as
before. We repeat this for all levels of cored blocks.

\vspace*{0.5 ex}\noindent \textsl{Step 5: Push all bilongitudinal $Q_{s_i}$
to be vertical.}\vspace*{0.5 ex}

Now, we examine the bilongitudinal levels. For a bilongitudinal level
$Q_{a_i}$ at a vertex $\delta_0$, corollary~\ref{coro:comeridian} shows
that the intersection circles are longitudes for $X_{a_i}$ and
$Y_{a_i}$. Lemma~\ref{lem:bilongitude} then shows that the circles of
$Q_{a_i}\cap P_{b_j}$ are isotopic in $Q_{a_i}$ and in $P_{b_j}$ to fibers.
First, use lemma~\ref{lem:fiberpreserving1} to find an isotopy preserving
levels, such that postcomposing $f_{\delta_0}$ by the isotopy makes the
intersection circles fibers of the $P_{b_j}$. Then, use
lemma~\ref{lem:fiberpreserving1} applied to $f_{\delta_0}^{-1}$ to find an
isotopy preserving levels of the domain, such that precomposing
$f_{\delta_0}$ by the isotopy makes the intersection circles the images of
fibers of $P_{s_i}$. After this process has been complete for the
bilongitudinal $Q_{a_i}$, the preimage (in their union $\cup Q_{a_i}$) of
each region $R(b_j,b_{j+1})$ with $b_j$ or $b_{j+1}$ in a bilongitudinal
block is a collection of fibered annuli which map into $R(b_j,b_{j+1})$ by
imbeddings that are fiber-preserving on their boundaries. We use
lemma~\ref{lem:fiberpreserving2} to find an isotopy that makes the
$Q_{a_i}$ vertical.  Again, we extend to a deformation of $f$ and work our
way up the skeleta, to assume that $Q_{s_i}(u)$ is vertical whenever $u\in
\Delta$ and $\Delta\subset U(t_i)$.

\vspace*{0.5 ex}\noindent \textsl{Step 6: Make $f$ fiber-preserving on
the complementary $S^1\times S^1\times I$ or solid tori of the
$P_{s_i}$-levels}\vspace*{0.5 ex}

We work our way up the skeleta one last time, using
lemma~\ref{lem:fiber-preserving on X} to make $f$ fiber-preserving on the
complementary $S^1\times S^1\times I$ or solid tori of the~$P_{a_i}$.
\end{proof}

\newpage
\section[Parameters in $D^d$]
{Parameters in $D^d$}
\label{sec:Dk}

Regard $D^d$ as the unit ball in $d$-dimensional Euclidean space, with
boundary the unit sphere $S^{d-1}$. As mentioned in
section~\ref{sec:reduction}, to prove that $\diff_f(L)\to \diff(L)$ is a
homotopy equivalence, we actually need to work with a family of
diffeomorphisms $f$ of $L$ parameterized by $D^d$, $d\geq 1$, for which
$f(u)$ is fiber-preserving whenever $u$ lies in the boundary
$S^{d-1}$. We must deform $f$ so that each $f(u)$ is
fiber-preserving, by a deformation that keeps $f(u)$ fiber-preserving at
all times when $u\in S^{d-1}$.

We now present a trick that allows us to gain good control of what happens
on $S^{d-1}$. The Hopf fibering we are using on $L$ can be described as a
Seifert fibering of $L$ over the round $2$-sphere $S$, in such a way that
each isometry of $L$ projects to an isometry of $S$ (details appear
in~\cite{M}, see also~\cite{MR}). By conjugating $\pi_1(L)$ in $\SO(4)$, we
may assume that the singular fibers, when $q>1$, are the preimages of the
poles. We choose our sweepout so that the level tori are the preimages of
latitude circles. Denote by $p_t$ the latitude circle that is the image of
the level torus~$P_t$.

There is an isotopy $J_t$ with $J_0$ the identity map of $L$ and each $J_t$
fiber-preserving, so that the images of the level tori $P_s$ under $J_1$
project to circles in the $2$-sphere as indicated in figure~\ref{fig:Dk
trick}. Denote the image of $J_1(P_s)$ in $S$ by~$q_s$. Their key property
is that when moved by any orthogonal rotation of $S$, each $p_t$ meets the
image of some $q_s$ transversely in two or four points.
\begin{figure}
\includegraphics[width=25 ex]{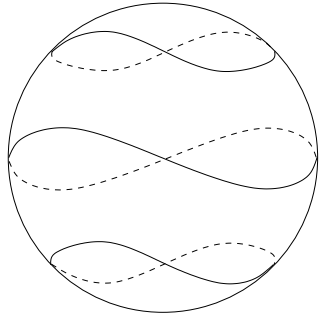}
\caption{Projections of the $J_1(P_t)$ into the $2$-sphere.}
\label{fig:Dk trick}
\end{figure}

Using theorem~\ref{thm:reduce to fiber-preserving}, we may assume that
$f_u$ is actually an isometry of $L$ for each $u\in S^{d-1}$.  Denote the
isometry that $f_u$ induces on $S$ by $\overline{f_u}$.  Now, deform the
entire family $f$ by precomposing each $f_u$ with $J_t$. At points in
$S^{d-1}$, each $f_u\circ J_t$ is fiber-preserving, so this is an allowable
deformation of $f$.  At the end of the deformation, for each $u\in
S^{d-1}$, $f_u\circ J_1(P_s)$ is a fibered torus $Q_s$ that projects to
$\overline{f_u}(q_s)$. Since $\overline{f_u}$ is an isometry of $S$, it
follows that for any latitude circle $p_t$, some $\overline{f_u}(q_s)$
meets $p_t$ transversely, in either two or four points. So $P_t$ and this
$Q_s$ meet transversely in either two or four circles which are fibers of
$L$. In particular, they are in very good position. We call such a pair
$P_t$ and $Q_s$ at $u$ an \emph{instant pair.}

Cover $S^{d-1}$ by finitely many open sets $Z_i'$ such that for each $i$,
there is an $(x_i,y_i)$ such that $Q_{x_i}$ and $P_{y_i}$ are an instant
pair at every point of $\overline{Z_i'}$. We may assume that there are open
sets $Z_i$ in $D^d$ such that $\overline{Z_i}\cap S^{d-1}=\overline{Z_i'}$
and $Q_{x_i}$ and $P_{y_i}$ meet in very good position at each point
of~$\overline{Z_i}$. For any sufficiently small deformation of $f$,
$Q_{x_i}$ and $P_{y_i}$ will still meet in very good position at all points
of~$\overline{Z_i}$. Let $V$ be a neighborhood of $S^{d-1}$ in $D^d$ such
that $\overline{V}$ is contained in the union of the~$Z_i$.

Now, we apply to $D^d$ the entire process used for the case when the
parameters lie in $S^d$, using appropriate fiber-preserving deformations at
parameters in $S^{d-1}$. Here are the steps:
\begin{enumerate}
\item By theorem~\ref{thm:generalposition}, there are arbitrarily small
deformations of $f$ that put it in general position with respect to the
sweepout. Select the deformation sufficiently small so that the $Q_{s_i}$
and $P_{t_i}$ still meet in very good position at every point of
$\overline{Z_i}$. Within $V$, we taper the deformation off to the identity,
so that no change has taken place at parameters in $S^{d-1}$. At every
parameter, either there is already a pair in very good position, or $f_u$
satisfies the conditions (GP1), (GP2), and~(GP3) of a general position
family.
\item
Theorem~\ref{thm:finding good levels} guarantees that at each of the
parameters in $D^d-V$, there is a pair $Q_s$ and $P_t$ meeting in good
position.
\item
Applying theorem~\ref{thm:from good to very good} to $D^d$, with $S^{d-1}$
in the role of $W_0$, we find a deformation of $f$, fixed on $S^{d-1}$, and
a covering $U(t_i)$ of $D^d$ and associated values $s_i$ so that for every
$u\in U(t_i)$, $Q_{s_i}$ and $P_{t_i}$ meet in very good position, and
$Q_{s_i}$ has no discal intersection with any~$P_{t_j}$.
\item
In the pushout step of the proof of theorem~\ref{thm:make
fiber-preserving}, we may assume that all the $U(t_i)$ that meet $S^{d-1}$
are the open sets $Z_i$. At parameters $u$ in $S^{d-1}$, the annuli to be
pushed out of each $V_{t_i}$ will be vertical annuli. So the pushouts may be
performed using fiber-preserving isotopies at these parameters, because the
necessary deformations can be taken as lifts of deformations of circles in
the quotient sphere $S$, and \cite{MR} provides fiber-preserving lifts of
any such deformations.
\item
After the triangulation of $D^d$ is chosen, the deformation that move the
$Q_{s_i}$ onto level tori can be performed using fiber-preserving isotopies
at parameters in $S^{d-1}$, again because the necessary deformations cover
deformations of circles in the quotient surface $S$. No further deformation
will be needed on simplices in $S^{d-1}$, since the $f_u$ are already
fiber-preserving there.
\end{enumerate}

This completes the discussion of the case of parameters in $D^d$, and the
proof of the Smale Conjecture for lens spaces.

\bibliographystyle{amsplain}

\end{document}